\documentclass[11pt,reqno]{amsart} 

\evensidemargin=\oddsidemargin
\usepackage{amsmath,amssymb}
\usepackage{amsfonts}
\usepackage{times}
\usepackage{graphicx}
\usepackage{color}
\definecolor{webred}{rgb}{0.75,0,0}
\definecolor{webgreen}{rgb}{0,0.75,0}
\usepackage[citecolor=webgreen,colorlinks=true,linkcolor=webred]{hyperref}

\begin{document}
\title[Primal and Shadow functions, Dual and Dual-Shadow functions]
{Primal and Shadow functions, Dual and Dual-Shadow functions\\ for
a circular crack and a circular $90^\circ$ V-notch\\ with Neumann boundary conditions}


\author{Samuel SHANNON, Zohar YOSIBASH, Monique DAUGE and Martin COSTABEL}

\address{{\bf Z.Y. \& S.S.}: Pearlstone Center for Aeronautical Eng.\ Studies, Dept.\ of Mechanical Engineering,
        Ben-Gurion University of the Negev, Beer-Sheva, 84105, Israel, \ \
              Tel.: +972-8-6477103, \ \ Fax: +972-8-6477101 } 
\email{zohary@bgu.ac.il} 
              
\address{{\bf M.D. \& M.C.}: IRMAR, Universit\'{e} de Rennes 1,
Campus de Beaulieu,
35042 Rennes Cedex, France}
\email{monique.dauge@univ-rennes1.fr}
\email{martin.costabel@univ-rennes1.fr}

\begin{abstract}
This report presents explicit analytical expressions for the primal, primal shadows, dual and dual shadows functions for the Laplace equation in the vicinity of a circular singular edge with Neumann boundary conditions on the faces that intersect at the singular edge. Two configurations are investigated: a penny-shaped crack and a $90^\circ$ V-notch.
\end{abstract}

\maketitle
\thispagestyle{empty}

\tableofcontents
\section*{Introduction}
The edges that we consider are in an axisymmetric configuration. Cylindrical coordinates are denoted by $(r,\theta, x_3)$ with the distance $r$ to the axis, the rotation angle $\theta\in[0,2\pi]$ and the coordinate $x_3$ along the axis. The circular edges are generated by rotating a point $P$ in the $r$-$x_3$ half-plane
around the $x_3$ axis. Let $R$ be $r(P)$. We consider the polar coordinate system $(\rho,\varphi)$ centered at $P$ such that the following relations hold
\[
   r = R + \rho\cos\varphi \quad\mbox{and}\quad x_3 = \rho\sin\varphi\,.
\]
In the vicinity of the edge the domain coincides with the set
\[
   \{(\rho,\varphi,\theta),\quad 
   0<\rho<\rho_0,\ \ \varphi_1<\varphi<\varphi_2,\ \ \theta\in[0,2\pi]\}\,.
\]

\newpage
In this report we address the cases of
\begin{itemize}
\item the penny-shaped crack, for which $\varphi_1=-\pi$ and $\varphi_2=\pi$
\item the $90^\circ$ V-notch that corresponds to $\varphi_1=-\pi$ and $\varphi_2=\pi/2$.
\end{itemize}

Solutions $\tau$ of the homogeneous Laplace equation $\Delta\tau=0$ with zero Neumann boundary conditions on the faces abutting the edge
can be expressed as follows \cite{YoSh11} in coordinates $(\rho,\varphi,\theta)$:
\begin{equation} \label{e.tau_gen}
\tau(\rho,\varphi,\theta) =\sum_{j \ge 0} \sum_{h=0,2,4\cdots} \partial_\theta^h
    A_j(\theta)\rho ^{\alpha _j}\sum _{f \ge 0} \left(\frac{\rho }{R}\right)^{h+f}
    \phi _{h, j, f}(\varphi),
\end{equation}
The functions $A_j$ are coefficients which depend on $\tau$, whereas the functions $\phi _{h, j, f}$ are completely determined by the geometry of the domain around the edge, I.E. by the limiting angles $\varphi_1$ and $\varphi_2$. For each natural number $j$, the term $\phi _{0, j, 0}$ is an eigenfunction of a 1D Laplace operator and for $h+j>0$, the ``shadows'' $\phi _{h, j, f}(\varphi)$ are obtained by solving the following recursive set of equations, cf
\cite[eq.\ (15-17) and (36)]{YoSh11}.

\vspace{0.1in}
\noindent{For $h= 0$}
\begin{eqnarray}
\alpha_j^2 \phi_{0,j,0} + \phi''_{0,j,0} & = & 0
\label{e.200} \\
(\alpha_j+1)^2 \phi_{0,j,1} + \phi''_{0,j,1} & = & - \left( \alpha_j \cos \varphi \phi_{0,j,0} 
                    - \sin\varphi \phi'_{0,j,0} \right)
\label{e.201} \\
(\alpha_j+i)^2 \phi_{0,j,f} + \phi''_{0,j,f} & = &  - \left[ (\alpha_j+f) (\alpha_j+f-1) \cos \varphi \phi_{0,j,f-1}
                            - \sin\varphi \phi'_{0,j,f-1} + \cos\varphi \phi''_{0,j,f-1} \right]
\nonumber \\
& & \quad \quad \quad \quad \quad \quad f\geq 2
\label{e.202}
\end{eqnarray}

\noindent{For $h= 2,4,6\cdots,$}
\begin{eqnarray}
(\alpha_j+f+h)^2\phi_{h,j,f} + \phi_{h,j,f}'' &=&
  -\ (h+f+\alpha_j-1)\left[ 2 (h+f+\alpha_j)-1 \right] \cos\varphi \phi_{h,j,(f-1)}
  \hskip -10mm \nonumber \\
&& 
+\ \sin\varphi \phi_{h,j,(f-1)}' - 2 \cos\varphi \phi_{h,j,(f-1)}''
\nonumber \\
& & 
            -\  (h+\alpha_j+f-2)(h+\alpha_j+f-1)\cos^2\varphi \phi_{h,j,(f-2)}
 \nonumber\\
&&    +\ \cos\varphi \sin\varphi \phi_{h,j,(f-2)}' - \cos^2\varphi \phi_{h,j,(f-2)}'' -
             \partial_\varphi \phi_{(h-2),j,f}
\label{e.450}
\end{eqnarray}

\noindent These equations have to be completed by the Neumann boundary conditions:
\begin{equation}
\partial_\varphi \phi_{h,j,f} = 0, \quad \text{on}\quad \varphi=\varphi_1 
\quad\mbox{and}\quad\varphi=\varphi_2
\label{e.Lap_Neumann_BCs}
\end{equation}

\medskip
Equations (\ref{e.200}-\ref{e.202}) with $h=0$ are associated with the axisymmetric case, 
and for $h=2,4,6,\cdots$ equations \eqref{e.450} are associated with the non-axisymmetric case.

\vspace{0.2in}
Each eigenfunction and shadow, has a dual and dual shadow counterpart, that can be computed by the following
recursive system.

For the dual eigenfunctions and their shadows  $\psi_{0,j,f}(\varphi)$ (for $h=0$ 
associated with the axisymmetric case), the recursive equation to be solved is \cite{YoSh11}:

\vspace{0.1in}
\noindent{For $h= 0$}
\vspace{-0.08in}
\begin{eqnarray} 
\alpha_j^2 \psi_{0,j,0} + \psi''_{0,j,0} & = & 0,
\label{e.200n} \\
(-\alpha_j+1)^2 \psi_{0,j,1} + \psi''_{0,j,1} & = &
    - \left( -\alpha_j \cos \varphi\  \psi_{0,j,0} - \sin\varphi\  \psi'_{0,j,0} \right),
\label{e.201n} \\
(-\alpha_j+i)^2 \psi_{0,j,f} + \psi''_{0,j,f} & = &
    - \left[ (-\alpha_j+f) (-\alpha_j+f-1) \cos \varphi\  \psi_{0,j,f-1}\right.
\label{e.202n} \\
& & \left.  - \sin\varphi\  \psi'_{0,j,f-1} + \cos\varphi\  \psi''_{0,j,f-1} \right],\quad f\geq 2,
\nonumber
\end{eqnarray}

\noindent{For $h= 2,4,6\cdots,$}
\vspace{-0.08in}
\begin{eqnarray}
& &
  (-\alpha_j+f+h)^2\psi_{h,j,f} + \psi_{h,j,f}'' = 
\nonumber  \\
& & \hskip 8ex
  - (h+f-\alpha_j-1)\left[ 2 (h+f-\alpha_j)-1 \right] \cos\varphi\ \psi_{h,j,(f-1)}
   +\sin\varphi\ \psi_{h,j,(f-1)}'
\nonumber\\
& & \hskip 8ex
   - 2 \cos\varphi\ \psi_{h,f,(f-1)}''  - (h-\alpha_j+f-2)(h-\alpha_j+f-1)\cos^2\varphi\ \psi_{h,j,(f-2)}
\nonumber \\
& & \hskip 8ex
   +\cos\varphi\ \sin\varphi\ \psi_{h,j,(f-2)}' - \cos^2\varphi\ \psi_{h,j,(f-2)}'' -\psi_{(h-2),j,f}
\end{eqnarray}
with homogeneous Neumann boundary conditions
\begin{equation}
\partial_\varphi \psi_{h,j,f}  =  0 \quad \text{on}\quad \varphi=\varphi_1 
\quad\mbox{and}\quad\varphi=\varphi_2.
\label{e.Neumann_BCs_psi}
\end{equation}

The above systems of equations have been solved using the formal calculus software Mathematica.

\newpage
\section{Primal functions and shadows for the crack $-\pi \le \varphi \le \pi$}
\label{s.fnc_shadows_axi_crack}
\subsection{First singular exponent ($j=1$)}

\begin{eqnarray}
\phi _{0,1,0}& =& \sin \frac{\varphi }{2}
\nonumber\\
\phi _{0,1,1}& =&\frac{1}{4} \sin \frac{\varphi }{2}
\nonumber\\
\phi _{0,1,2}& =&\frac{1}{12} \sin \frac{\varphi }{2}-\frac{3}{32} \sin \frac{3 \varphi }{2}
\nonumber\\
\phi _{0,1,3}& =&\frac{1}{16} \sin \frac{\varphi }{2}-\frac{1}{30} \sin \frac{3 \varphi }{2}+\frac{5}{128} \sin \frac{5 \varphi }{2}
\nonumber\\
\phi _{0,1,4}& =&\frac{3}{80} \sin \frac{\varphi }{2}-\frac{5}{128} \sin \frac{3 \varphi }{2}+\frac{1}{70} \sin \frac{5 \varphi }{2}-\frac{35}{2048} \sin \frac{7 \varphi }{2}
\nonumber\\
\phi _{0,1,5}& =&\frac{1}{32}\sin \frac{\varphi }{2}-\frac{3}{140}\sin \frac{3 \varphi }{2}+\frac{35}{1536}\sin \frac{5 \varphi }{2}-\frac{2}{315}\sin \frac{7 \varphi }{2}+\frac{63}{8192}\sin \frac{9 \varphi }{2}
\nonumber\\
\phi _{0,1,6}& =&\frac{5}{224}\sin \frac{\varphi }{2}-\frac{35}{1536}\sin \frac{3 \varphi }{2}+\frac{1}{84}\sin \frac{5 \varphi }{2}-\frac{105}{8192}\sin \frac{7 \varphi }{2}+\frac{2}{693}\sin \frac{9 \varphi }{2}-\frac{231}{65536}\sin \frac{11 \varphi }{2}
\nonumber\\
\phi _{0,1,7}& =&\frac{5}{256}\sin \frac{\varphi }{2}-\frac{5}{336}\sin \frac{3 \varphi }{2}+\frac{63}{4096}\sin \frac{5 \varphi }{2}-\frac{1}{154}\sin \frac{7 \varphi }{2}+\frac{231}{32768}\sin \frac{9 \varphi }{2}-\frac{4}{3003}\sin \frac{11 \varphi }{2}
\nonumber\\& + &\frac{429}{262144}\sin \frac{13 \varphi }{2}
\nonumber\\
\phi _{0,1,8}& =&\frac{35}{2304}\sin \frac{\varphi }{2}-\frac{63}{4096}\sin \frac{3 \varphi }{2}+\frac{5}{528}\sin \frac{5 \varphi }{2}-\frac{1617}{163840}\sin \frac{7 \varphi }{2}+\frac{1}{286}\sin \frac{9 \varphi }{2}-\frac{1001}{262144}\sin \frac{11 \varphi }{2}
\nonumber\\
& +& \frac{4}{6435}\sin \frac{13 \varphi }{2}-\frac{6435}{8388608}\sin \frac{15 \varphi }{2}
\nonumber\\
\phi_{0,1,9}& =&\frac{7}{512}\sin \frac{\varphi }{2}-\frac{35}{3168}\sin \frac{3 \varphi }{2}+\frac{231}{20480}\sin \frac{5 \varphi }{2}-\frac{5}{858}\sin \frac{7 \varphi }{2}+\frac{1001}{163840}\sin \frac{9 \varphi }{2}
\nonumber\\
& -& \frac{4}{2145}\sin \frac{11 \varphi }{2}+\frac{2145}{1048576}\sin \frac{13 \varphi }{2}-\frac{32}{109395}\sin \frac{15 \varphi }{2}+\frac{12155}{33554432}\sin \frac{17 \varphi }{2}
\nonumber\\
\phi _{0,1,10}& =&\frac{63}{5632}\sin \frac{\varphi }{2}-\frac{231}{20480}\sin \frac{3 \varphi }{2}+\frac{35}{4576}\sin \frac{5 \varphi }{2}-\frac{1287}{163840}\sin \frac{7 \varphi }{2}+\frac{1}{286}\sin \frac{9 \varphi }{2}
\nonumber\\
& -& \frac{3861}{1048576}\sin \frac{11 \varphi }{2}+\frac{12}{12155}\sin \frac{13 \varphi }{2}-\frac{36465}{33554432}\sin \frac{15 \varphi }{2}+\frac{32}{230945}\sin \frac{17 \varphi }{2}
\nonumber\\
& -& \frac{46189}{268435456}\sin \frac{19 \varphi }{2}
\nonumber
\end{eqnarray}

\newpage

\begin{eqnarray}
\phi _{2,1,0}& = &-\frac{1}{6}\sin \frac{\varphi }{2}
\nonumber\\
\phi _{2,1,1}& = &-\frac{1}{8}\sin \frac{\varphi }{2}+\frac{7}{60}\sin \frac{3 \varphi }{2}
\nonumber\\
\phi _{2,1,2}& = &-\frac{31}{240}\sin \frac{\varphi }{2}+\frac{5}{64}\sin \frac{3 \varphi }{2}-\frac{19}{280}\sin \frac{5 \varphi }{2}
\nonumber\\
\phi _{2,1,3}& = &-\frac{41}{384}\sin \frac{\varphi }{2}+\frac{111}{1120} \sin \frac{3 \varphi }{2}-\frac{35}{768} \sin \frac{5 \varphi }{2}+\frac{187}{5040} \sin \frac{7 \varphi }{2}
\nonumber\\
\phi _{2,1,4}& = &-\frac{8209}{80640}\sin \frac{\varphi }{2}+\frac{119}{1536} \sin \frac{3 \varphi }{2}-\frac{271}{4032} \sin \frac{5 \varphi }{2}+\frac{105}{4096}\sin \frac{7 \varphi }{2}-\frac{437}{22176}\sin \frac{9 \varphi }{2}
\nonumber\\
\phi _{2,1,5}& = &-\frac{4069}{46080} \sin \frac{\varphi }{2}+\frac{39899}{483840} \sin \frac{3 \varphi }{2}-\frac{427}{8192} \sin \frac{5 \varphi }{2}+\frac{1259}{29568} \sin \frac{7 \varphi }{2}-\frac{231}{16384} \sin \frac{9 \varphi }{2}
\nonumber\\
&+&\frac{1979}{192192}\sin \frac{11 \varphi }{2}
\nonumber\\
\phi _{2,1,6}& = &-\frac{483047}{5806080} \sin \frac{\varphi }{2}+\frac{5663}{81920} \sin \frac{3 \varphi }{2}-\frac{18467}{304128} \sin \frac{5 \varphi }{2}+\frac{5467}{163840} \sin \frac{7 \varphi }{2}-\frac{2837}{109824} \sin \frac{9 \varphi }{2}
\nonumber\\
&+&\frac{1001}{131072}\sin \frac{11 \varphi }{2}-\frac{4387}{823680}\sin \frac{13 \varphi }{2}
\nonumber\\
\phi _{2,1,7}& = &-\frac{153325}{2064384} \sin \frac{\varphi }{2}+\frac{8920073}{127733760} \sin \frac{3 \varphi }{2}-\frac{49621}{983040} \sin \frac{5 \varphi }{2}+\frac{331871}{7907328} \sin \frac{7 \varphi }{2}
\nonumber\\
&-&\frac{27027}{1310720} \sin \frac{9 \varphi }{2}+\frac{50111}{3294720}\sin \frac{11 \varphi }{2}-\frac{2145}{524288}\sin \frac{13 \varphi }{2}+\frac{76627}{28005120}\sin \frac{15 \varphi }{2}
\nonumber\\
\phi _{2,1,8}& = &-\frac{79481719}{1135411200} \sin \frac{\varphi }{2}+\frac{2514853}{41287680} \sin \frac{3 \varphi }{2}-\frac{6664667}{123002880} \sin \frac{5 \varphi }{2}+\frac{275561}{7864320} \sin \frac{7 \varphi }{2}
\nonumber\\
&-&\frac{81201}{2928640} \sin \frac{9 \varphi }{2}+\frac{13013}{1048576}\sin \frac{11 \varphi }{2}-\frac{327121}{37340160}\sin \frac{13 \varphi }{2}+\frac{36465}{16777216}\sin \frac{15 \varphi }{2}
\nonumber\\
&-&\frac{165409}{118243840}\sin \frac{17 \varphi }{2}
\nonumber
\end{eqnarray}

\bigskip
\begin{eqnarray}
\phi _{4,1,0}& = &\frac{1}{120}\sin \frac{\varphi }{2}
\nonumber\\
\phi _{4,1,1}& = &\frac{1}{96}\sin \frac{\varphi }{2}-\frac{3}{280}\sin \frac{3 \varphi }{2}
\nonumber\\
\phi _{4,1,2}& = &\frac{17}{1008}\sin \frac{\varphi }{2}-\frac{7}{768}\sin \frac{3 \varphi }{2}+\frac{19}{2016}\sin \frac{5 \varphi }{2}
\nonumber\\
\phi _{4,1,3}& = &\frac{41}{2304}\sin \frac{\varphi }{2}-\frac{433}{24192}\sin \frac{3 \varphi }{2}+\frac{7}{1024}\sin \frac{5 \varphi }{2}-\frac{13}{1848}\sin \frac{7 \varphi }{2}
\nonumber\\
\phi _{4,1,4}& = &\frac{1019}{46080}\sin \frac{\varphi }{2}-\frac{1}{64}\sin \frac{3 \varphi }{2}+\frac{595}{38016}\sin \frac{5 \varphi }{2}-\frac{77}{16384}\sin \frac{7 \varphi }{2}+\frac{263}{54912}\sin \frac{9 \varphi }{2}
\nonumber\\
\phi _{4,1,5}& = &\frac{1367 \sin \frac{\varphi }{2}}{61440}-\frac{5617 \sin \frac{3 \varphi }{2}}{253440}+\frac{605 \sin \frac{5 \varphi }{2}}{49152}-\frac{24155 \sin \frac{7 \varphi }{2}}{1976832}+\frac{1001 \sin \frac{9 \varphi }{2}}{327680}-\frac{2533 \sin \frac{11 \varphi }{2}}{823680}
\nonumber\\
\phi _{4,1,6}& = &\frac{3211069 \sin \frac{\varphi }{2}}{127733760}-\frac{58267 \sin \frac{3 \varphi }{2}}{2949120}+\frac{514259 \sin \frac{5 \varphi }{2}}{26357760}-\frac{4433 \sin \frac{7 \varphi }{2}}{491520}+\frac{81 \sin \frac{9 \varphi }{2}}{9152}
\nonumber\\
&-&\frac{1001 \sin \frac{11 \varphi }{2}}{524288}+\frac{35401 \sin \frac{13 \varphi }{2}}{18670080}
\nonumber
\end{eqnarray}
\newpage

\begin{eqnarray}
\phi _{6,1,0}& = &-\frac{1}{5040}\sin \frac{\varphi }{2}
\nonumber\\
\phi _{6,1,1}& = &-\frac{1}{2880}\sin \frac{\varphi }{2}+\frac{11}{30240}\sin \frac{3 \varphi }{2}
\nonumber\\
\phi _{6,1,2}& = &-\frac{13}{17280}\sin \frac{\varphi }{2}+\frac{1}{2560}\sin \frac{3 \varphi }{2}-\frac{1}{2376}\sin \frac{5 \varphi }{2}
\nonumber\\
\phi _{6,1,3}& = &-\frac{1}{1024}\sin \frac{\varphi }{2}+\frac{191}{190080}\sin \frac{3 \varphi }{2}-\frac{11}{30720}\sin \frac{5 \varphi }{2}+\frac{97}{247104}\sin \frac{7 \varphi }{2}
\nonumber\\
\phi _{6,1,4}& = &-\frac{511}{345600}\sin \frac{\varphi }{2}+\frac{187}{184320}\sin \frac{3 \varphi }{2}-\frac{1751}{1647360}\sin \frac{5 \varphi }{2}+\frac{143}{491520}\sin \frac{7 \varphi }{2}-\frac{59}{183040}\sin \frac{9 \varphi }{2}
\nonumber
\end{eqnarray}

\bigskip
\begin{eqnarray}
\phi _{8,1,0}& = &\frac{1}{362880}\sin \frac{\varphi }{2}
\nonumber\\
\phi _{8,1,1}& = &\frac{\sin \frac{\varphi }{2}}{161280}-\frac{13 \sin \frac{3 \varphi }{2}}{1995840}
\nonumber\\
\phi _{8,1,2}& = &\frac{89 \sin \frac{\varphi }{2}}{5322240}-\frac{11 \sin \frac{3 \varphi }{2}}{1290240}+\frac{107 \sin \frac{5 \varphi }{2}}{11531520}
\nonumber
\end{eqnarray}

\bigskip
\begin{eqnarray}
\phi _{10,1,0}& = &-\frac{1}{39916800}\sin \frac{\varphi }{2}
\nonumber
\end{eqnarray}

\newpage
\subsection{Second singular exponent ($j=3$)}

\begin{eqnarray}
\phi _{0,3,0}& =&\sin \frac{3\varphi }{2}
\nonumber\\
\phi _{0,3,1}& =&-\frac{1}{4} \sin \frac{\varphi }{2}
\nonumber\\
\phi _{0,3,2}& =&-\frac{3}{32} \sin \frac{\varphi }{2}+\frac{1}{10} \sin \frac{3 \varphi }{2}
\nonumber\\
\phi _{0,3,3}& =&-\frac{3}{40} \sin \frac{\varphi }{2}+\frac{5}{128} \sin \frac{3 \varphi }{2}-\frac{3}{70} \sin \frac{5 \varphi }{2}
\nonumber\\
\phi _{0,3,4}& =&-\frac{3}{64} \sin \frac{\varphi }{2}+\frac{27}{560} \sin \frac{3 \varphi }{2}-\frac{35}{2048} \sin \frac{5 \varphi }{2}+\frac{2}{105} \sin \frac{7 \varphi }{2}
\nonumber\\
\phi _{0,3,5}& =&-\frac{9}{224}\sin \frac{\varphi }{2}+\frac{7}{256} \sin \frac{3 \varphi }{2}-\frac{1}{35}\sin \frac{5 \varphi }{2}+\frac{63}{8192}\sin \frac{7 \varphi }{2}-\frac{2}{231}\sin \frac{9 \varphi }{2}
\nonumber\\
\phi _{0,3,6}& =&-\frac{15}{512}\sin \frac{\varphi }{2}+\frac{5}{168} \sin \frac{3 \varphi }{2}-\frac{63}{4096}\sin \frac{5 \varphi }{2}+\frac{5}{308}\sin \frac{7 \varphi }{2}-\frac{231}{65536}\sin \frac{9 \varphi }{2}
\nonumber\\
& +& \frac{4}{1001}\sin \frac{11 \varphi }{2}
\nonumber\\
\phi _{0,3,7}& =&-\frac{5}{192}\sin \frac{\varphi }{2}+\frac{81}{4096} \sin \frac{3 \varphi }{2}-\frac{25}{1232}\sin \frac{5 \varphi }{2}+\frac{693}{81920}\sin \frac{7 \varphi }{2}-\frac{9}{1001}\sin \frac{9 \varphi }{2}
\nonumber\\
&+& \frac{429}{262144}\sin \frac{11 \varphi }{2}-\frac{4}{2145} \sin \frac{13 \varphi }{2}
\nonumber\\
\phi _{0,3,8}& =&-\frac{21}{1024}\sin \frac{\varphi }{2}+\frac{175}{8448} \sin \frac{3 \varphi }{2}-\frac{2079}{163840}\sin \frac{5 \varphi }{2}+\frac{15}{1144}\sin \frac{7 \varphi }{2}-\frac{3003}{655360}\sin \frac{9 \varphi }{2}
\nonumber\\
& +& \frac{7}{1430}\sin \frac{11 \varphi }{2}-\frac{6435}{8388608} \sin \frac{13 \varphi }{2}+\frac{32}{36465}\sin \frac{15 \varphi }{2}
\nonumber\\
\phi _{0,3,9}& =&-\frac{105}{5632}\sin \frac{\varphi }{2}+\frac{77}{5120} \sin \frac{3 \varphi }{2}-\frac{35}{2288}\sin \frac{5 \varphi }{2}+\frac{1287}{163840}\sin \frac{7 \varphi }{2}-\frac{7}{858}\sin \frac{9 \varphi }{2}
\nonumber\\
&+&\frac{1287}{524288}\sin \frac{11 \varphi }{2}-\frac{32}{12155} \sin \frac{13 \varphi }{2}+\frac{12155}{33554432}\sin \frac{15 \varphi }{2}-\frac{96}{230945}\sin \frac{17 \varphi }{2}
\nonumber
\end{eqnarray}

\begin{eqnarray}
\phi _{2,3,0}& =&-\frac{1}{10} \sin \frac{3 \varphi }{2}
\nonumber\\
\phi _{2,3,1}& =&\frac{3}{40} \sin \frac{\varphi }{2}+\frac{11}{140} \sin \frac{5 \varphi }{2}
\nonumber\\
\phi _{2,3,2}& =&\frac{3}{64} \sin \frac{\varphi }{2}-\frac{7}{80} \sin \frac{3 \varphi }{2}-\frac{41}{840} \sin \frac{7 \varphi }{2}
\nonumber\\
\phi _{2,3,3}& =&\frac{65}{896} \sin \frac{\varphi }{2}-\frac{7}{256} \sin \frac{3 \varphi }{2}+\frac{145 \sin \frac{5 \varphi }{2}}{2016}+\frac{103 \sin \frac{9 \varphi }{2}}{3696}
\nonumber\\
\phi _{2,3,4}& =&\frac{27}{512} \sin \frac{\varphi }{2}-\frac{17959 \sin \frac{3 \varphi }{2}}{241920}+\frac{63 \sin \frac{5 \varphi }{2}}{4096}-\frac{2263 \sin \frac{7 \varphi }{2}}{44352}-\frac{489 \sin \frac{11 \varphi }{2}}{32032}
\nonumber\\
\phi _{2,3,5}& =&\frac{8929 \sin \frac{\varphi }{2}}{138240}-\frac{291 \sin \frac{3 \varphi }{2}}{8192}+\frac{2483 \sin \frac{5 \varphi }{2}}{39424}-\frac{693 \sin \frac{7 \varphi }{2}}{81920}+\frac{12829 \sin \frac{9 \varphi }{2}}{384384}+\frac{1121 \sin \frac{13 \varphi }{2}}{137280}
\nonumber\\
\phi _{2,3,6}& =&\frac{4151 \sin \frac{\varphi }{2}}{81920}-\frac{1358089 \sin \frac{3 \varphi }{2}}{21288960}+\frac{3729 \sin \frac{5 \varphi }{2}}{163840}-\frac{63161 \sin \frac{7 \varphi }{2}}{1317888}+\frac{3003 \sin \frac{9 \varphi }{2}}{655360}
\nonumber\\
&-&\frac{34169 \sin \frac{11 \varphi }{2}}{1647360}-\frac{20081 \sin \frac{15 \varphi }{2}}{4667520}
\nonumber\\
\phi _{2,3,7}& =&\frac{2161361 \sin \frac{\varphi }{2}}{37847040}-\frac{21857 \sin \frac{3 \varphi }{2}}{589824}+\frac{10216499 \sin \frac{5 \varphi }{2}}{184504320}-\frac{18447 \sin \frac{7 \varphi }{2}}{1310720}+\frac{1342447 \sin \frac{9 \varphi }{2}}{39536640}
\nonumber\\
&-&\frac{1287 \sin \frac{11 \varphi }{2}}{524288}+\frac{63367 \sin \frac{13 \varphi }{2}}{5091840}+\frac{132641 \sin \frac{17 \varphi }{2}}{59121920}
\nonumber
\end{eqnarray}

\begin{eqnarray}
\phi _{4,3,0}& =&\frac{1}{280} \sin \frac{3 \varphi }{2}
\nonumber\\
\phi _{4,3,1}& =&-\frac{1}{224} \sin \frac{\varphi }{2}-\frac{13 \sin \frac{5 \varphi }{2}}{2520}
\nonumber\\
\phi _{4,3,2}& =&-\frac{1}{256} \sin \frac{\varphi }{2}+\frac{7}{864} \sin \frac{3 \varphi }{2}+\frac{109 \sin \frac{7 \varphi }{2}}{22176}
\nonumber\\
\phi _{4,3,3}& =&-\frac{59 \sin \frac{\varphi }{2}}{6912}+\frac{3 \sin \frac{3 \varphi }{2}}{1024}-\frac{25 \sin \frac{5 \varphi }{2}}{2688}-\frac{17 \sin \frac{9 \varphi }{2}}{4368}
\nonumber\\
\phi _{4,3,4}& =&-\frac{23 \sin \frac{\varphi }{2}}{3072}+\frac{17431 \sin \frac{3 \varphi }{2}}{1520640}-\frac{33 \sin \frac{5 \varphi }{2}}{16384}+\frac{355 \sin \frac{7 \varphi }{2}}{41184}+\frac{457 \sin \frac{11 \varphi }{2}}{164736}
\nonumber\\
\phi _{4,3,5}& =&-\frac{7781 \sin \frac{\varphi }{2}}{675840}+\frac{869 \sin \frac{3 \varphi }{2}}{147456}-\frac{120277 \sin \frac{5 \varphi }{2}}{9884160}+\frac{429 \sin \frac{7 \varphi }{2}}{327680}-\frac{13919 \sin \frac{9 \varphi }{2}}{1976832}
\nonumber\\
&-&\frac{25879 \sin \frac{13 \varphi }{2}}{14002560}
\nonumber
\end{eqnarray}

\begin{eqnarray}
\phi _{6,3,0}& =&-\frac{\sin \frac{3 \varphi }{2}}{15120}
\nonumber\\
\phi _{6,3,1}& =&\frac{\sin \frac{\varphi }{2}}{8640}+\frac{\sin \frac{5 \varphi }{2}}{7392}
\nonumber\\
\phi _{6,3,2}& =&\frac{\sin \frac{\varphi }{2}}{7680}-\frac{53 \sin \frac{3 \varphi }{2}}{190080}-\frac{7 \sin \frac{7 \varphi }{2}}{41184}
\nonumber\\
\phi _{6,3,3}& =&\frac{61 \sin \frac{\varphi }{2}}{168960}-\frac{11 \sin \frac{3 \varphi }{2}}{92160}+\frac{997 \sin \frac{5 \varphi }{2}}{2471040}+\frac{19 \sin \frac{9 \varphi }{2}}{112320}
\nonumber
\end{eqnarray}

\begin{eqnarray}
\phi _{8,3,0}& =&\frac{\sin \frac{3 \varphi }{2}}{1330560}
\nonumber\\
\phi _{8,3,1}& =&-\frac{\sin \frac{\varphi }{2}}{591360}-\frac{17 \sin \frac{5 \varphi }{2}}{8648640}
\nonumber
\end{eqnarray}

\newpage
\subsection{Third singular exponent ($j=5$)}

\begin{eqnarray}
\phi _{0,5,0}& =&\sin \frac{5\varphi }{2}
\nonumber\\
\phi _{0,5,1}& =&-\frac{1}{4} \sin \frac{3 \varphi }{2}
\nonumber\\
\phi _{0,5,2}& =&\frac{3}{32} \sin \frac{\varphi }{2}+\frac{3}{28} \sin \frac{5 \varphi }{2}
\nonumber\\
\phi _{0,5,3}& =&\frac{5}{128} \sin \frac{\varphi }{2}-\frac{9}{112} \sin \frac{3 \varphi }{2}-\frac{1}{21} \sin \frac{7 \varphi }{2}
\nonumber\\
\phi _{0,5,4}& =&\frac{45}{896} \sin \frac{\varphi }{2}-\frac{35}{2048} \sin \frac{3 \varphi }{2}+\frac{3}{56} \sin \frac{5 \varphi }{2}+\frac{5}{231} \sin \frac{9 \varphi }{2}
\nonumber\\
\phi _{0,5,5}& =&\frac{15}{512}\sin \frac{\varphi }{2}-\frac{5}{112}\sin \frac{3 \varphi }{2}+\frac{63}{8192}\sin \frac{5 \varphi }{2}-\frac{5}{154}\sin \frac{7 \varphi }{2}-\frac{10}{1001}\sin \frac{11 \varphi }{2}
\nonumber\\
\phi _{0,5,6}& =&\frac{25}{768}\sin \frac{\varphi }{2}-\frac{135}{8192}\sin \frac{3 \varphi }{2}+\frac{125}{3696}\sin \frac{5 \varphi }{2}-\frac{231}{65536}\sin \frac{7 \varphi }{2}+\frac{75}{4004}\sin \frac{9 \varphi }{2}+\frac{2}{429}\sin \frac{13 \varphi }{2}
\nonumber\\
\phi _{0,5,7}& =&\frac{45}{2048}\sin \frac{\varphi }{2}-\frac{125}{4224}\sin \frac{3 \varphi }{2}+\frac{297}{32768}\sin \frac{5 \varphi }{2}-\frac{375}{16016}\sin \frac{7 \varphi }{2}+\frac{429}{262144}\sin \frac{9 \varphi }{2}
\nonumber\\
&-&\frac{3}{286}\sin \frac{11 \varphi }{2}-\frac{16}{7293}\sin \frac{15 \varphi }{2}
\nonumber\\
\phi _{0,5,8}& =&\frac{525}{22528}\sin \frac{\varphi }{2}-\frac{231}{16384}\sin \frac{3 \varphi }{2}+\frac{875}{36608}\sin \frac{5 \varphi }{2}-\frac{1287}{262144}\sin \frac{7 \varphi }{2}+\frac{35}{2288}\sin \frac{9 \varphi }{2}
\nonumber\\
&-&\frac{6435}{8388608}\sin \frac{11 \varphi }{2}+\frac{14}{2431}\sin \frac{13 \varphi }{2}+\frac{48}{46189}\sin \frac{17 \varphi }{2}
\nonumber
\end{eqnarray}

\begin{eqnarray}
\phi _{2,5,0}&=&-\frac{1}{14} \sin \frac{5 \varphi }{2}
\nonumber\\
\phi _{2,5,1}&=&\frac{3}{56} \sin \frac{3 \varphi }{2}+\frac{5}{84} \sin \frac{7 \varphi }{2}
\nonumber\\
\phi _{2,5,2}&=&-\frac{15}{448} \sin \frac{\varphi }{2}-\frac{67 \sin \frac{5 \varphi }{2}}{1008}-\frac{71 \sin \frac{9 \varphi }{2}}{1848}
\nonumber\\
\phi _{2,5,3}&=&-\frac{5}{256} \sin \frac{\varphi }{2}+\frac{445 \sin \frac{3 \varphi }{2}}{8064}+\frac{1261 \sin \frac{7 \varphi }{2}}{22176}+\frac{361 \sin \frac{11 \varphi }{2}}{16016}
\nonumber\\
\phi _{2,5,4}&=&-\frac{185 \sin \frac{\varphi }{2}}{4608}+\frac{45 \sin \frac{3 \varphi }{2}}{4096}-\frac{31309 \sin \frac{5 \varphi }{2}}{532224}-\frac{725 \sin \frac{9 \varphi }{2}}{17472}-\frac{173 \sin \frac{13 \varphi }{2}}{13728}
\nonumber\\
\phi _{2,5,5}&=&-\frac{665 \sin \frac{\varphi }{2}}{24576}+\frac{15589 \sin \frac{3 \varphi }{2}}{304128}-\frac{99 \sin \frac{5 \varphi }{2}}{16384}+\frac{33853 \sin \frac{7 \varphi }{2}}{658944}+\frac{4567 \sin \frac{11 \varphi }{2}}{164736}+\frac{3197 \sin \frac{15 \varphi }{2}}{466752}
\nonumber\\
\phi _{2,5,6}&=&-\frac{21763 \sin \frac{\varphi }{2}}{540672}+\frac{1705 \sin \frac{3 \varphi }{2}}{98304}-\frac{961937 \sin \frac{5 \varphi }{2}}{18450432}+\frac{429 \sin \frac{7 \varphi }{2}}{131072}-\frac{263249 \sin \frac{9 \varphi }{2}}{6589440}
\nonumber\\
&-&\frac{5771 \sin \frac{13 \varphi }{2}}{329472}-\frac{21613 \sin \frac{17 \varphi }{2}}{5912192}
\nonumber
\end{eqnarray}

\begin{eqnarray}
\phi _{4,5,0}&=&\frac{1}{504} \sin \frac{5 \varphi }{2}
\nonumber\\
\phi _{4,5,1}&=&-\frac{5 \sin \frac{3 \varphi }{2}}{2016}-\frac{17 \sin \frac{7 \varphi }{2}}{5544}
\nonumber\\
\phi _{4,5,2}&=&\frac{5 \sin \frac{\varphi }{2}}{2304}+\frac{10 \sin \frac{5 \varphi }{2}}{2079}+\frac{295 \sin \frac{9 \varphi }{2}}{96096}
\nonumber\\
\phi _{4,5,3}&=&\frac{5 \sin \frac{\varphi }{2}}{3072}-\frac{35 \sin \frac{3 \varphi }{2}}{6912}-\frac{955 \sin \frac{7 \varphi }{2}}{164736}-\frac{1}{396} \sin \frac{11 \varphi }{2}
\nonumber\\
\phi _{4,5,4}&=&\frac{25 \sin \frac{\varphi }{2}}{5632}-\frac{55 \sin \frac{3 \varphi }{2}}{49152}+\frac{28211 \sin \frac{5 \varphi }{2}}{3953664}+\frac{919 \sin \frac{9 \varphi }{2}}{164736}+\frac{5197 \sin \frac{13 \varphi }{2}}{2800512}
\nonumber
\end{eqnarray}

\begin{eqnarray}
\phi _{6,5,0}&=&-\frac{\sin \frac{5 \varphi }{2}}{33264}
\nonumber\\
\phi _{6,5,1}&=&\frac{\sin \frac{3 \varphi }{2}}{19008}+\frac{19 \sin \frac{7 \varphi }{2}}{288288}
\nonumber\\
\phi _{6,5,2}&=&-\frac{\sin \frac{\varphi }{2}}{16896}-\frac{67 \sin \frac{5 \varphi }{2}}{494208}-\frac{\sin \frac{9 \varphi }{2}}{11440}
\nonumber
\end{eqnarray}

\begin{eqnarray}
\phi _{8,5,0}&=&\frac{\sin \frac{5 \varphi }{2}}{3459456}
\nonumber
\end{eqnarray}

\subsection{Fourth singular exponent ($j=7$)}

\begin{eqnarray}
\phi _{0,7,0}& =&\sin \frac{7\varphi }{2}
\nonumber\\
\phi _{0,7,1}& =&-\frac{1}{4} \sin \frac{5 \varphi }{2}
\nonumber\\
\phi _{0,7,2}& =&\frac{3}{32} \sin \frac{3 \varphi }{2}+\frac{1}{9} \sin \frac{7 \varphi }{2}
\nonumber\\
\phi _{0,7,3}& =&-\frac{5}{128} \sin \frac{\varphi }{2}-\frac{1}{12} \sin \frac{5 \varphi }{2}-\frac{5}{99} \sin \frac{9 \varphi }{2}
\nonumber\\
\phi _{0,7,4}& =&-\frac{35 \sin \frac{\varphi }{2}}{2048}+\frac{5}{96} \sin \frac{3 \varphi }{2}+\frac{5}{88} \sin \frac{7 \varphi }{2}+\frac{10}{429} \sin \frac{11 \varphi }{2}
\nonumber\\
\phi _{0,7,5}& =&-\frac{35 \sin \frac{\varphi }{2}}{1152}+\frac{63 \sin \frac{3 \varphi }{2}}{8192}-\frac{25}{528} \sin \frac{5 \varphi }{2}-\frac{5}{143} \sin \frac{9 \varphi }{2}-\frac{14 \sin \frac{13 \varphi }{2}}{1287}
\nonumber\\
\phi _{0,7,6}& =&-\frac{35 \sin \frac{\varphi }{2}}{2048}+\frac{875 \sin \frac{3 \varphi }{2}}{25344}-\frac{231 \sin \frac{5 \varphi }{2}}{65536}+\frac{125 \sin \frac{7 \varphi }{2}}{3432}+\frac{35 \sin \frac{11 \varphi }{2}}{1716}+\frac{112 \sin \frac{15 \varphi }{2}}{21879}
\nonumber\\
\phi _{0,7,7}& =&-\frac{525 \sin \frac{\varphi }{2}}{22528}+\frac{77 \sin \frac{3 \varphi }{2}}{8192}-\frac{875 \sin \frac{5 \varphi }{2}}{27456}+\frac{429 \sin \frac{7 \varphi }{2}}{262144}-\frac{175 \sin \frac{9 \varphi }{2}}{6864}-\frac{28 \sin \frac{13 \varphi }{2}}{2431}-\frac{112 \sin \frac{17 \varphi }{2}}{46189}
\nonumber
\end{eqnarray}

\begin{eqnarray}
\phi _{2,7,0}& =&-\frac{1}{18} \sin \frac{7 \varphi }{2}
\nonumber\\
\phi _{2,7,1}& =&\frac{1}{24} \sin \frac{5 \varphi }{2}+\frac{19}{396} \sin \frac{9 \varphi }{2}
\nonumber\\
\phi _{2,7,2}& =&-\frac{5}{192} \sin \frac{3 \varphi }{2}-\frac{85 \sin \frac{7 \varphi }{2}}{1584}-\frac{109 \sin \frac{11 \varphi }{2}}{3432}
\nonumber\\
\phi _{2,7,3}& =&\frac{35 \sin \frac{\varphi }{2}}{2304}+\frac{565 \sin \frac{5 \varphi }{2}}{12672}+\frac{647 \sin \frac{9 \varphi }{2}}{13728}+\frac{391 \sin \frac{13 \varphi }{2}}{20592}
\nonumber\\
\phi _{2,7,4}& =&\frac{35 \sin \frac{\varphi }{2}}{4096}-\frac{1645 \sin \frac{3 \varphi }{2}}{50688}-\frac{48259 \sin \frac{7 \varphi }{2}}{988416}-\frac{321 \sin \frac{11 \varphi }{2}}{9152}-\frac{7543 \sin \frac{15 \varphi }{2}}{700128}
\nonumber\\
\phi _{2,7,5}& =&\frac{5915 \sin \frac{\varphi }{2}}{270336}-\frac{77 \sin \frac{3 \varphi }{2}}{16384}+\frac{168343 \sin \frac{5 \varphi }{2}}{3953664}+\frac{430007 \sin \frac{9 \varphi }{2}}{9884160}+\frac{66637 \sin \frac{13 \varphi }{2}}{2800512}+\frac{17517 \sin \frac{17 \varphi }{2}}{2956096}
\nonumber
\end{eqnarray}

\begin{eqnarray}
\phi _{4,7,0}& =&\frac{1}{792} \sin \frac{7 \varphi }{2}
\nonumber\\
\phi _{4,7,1}& =&-\frac{5 \sin \frac{5 \varphi }{2}}{3168}-\frac{7 \sin \frac{9 \varphi }{2}}{3432}
\nonumber\\
\phi _{4,7,2}& =&\frac{35 \sin \frac{3 \varphi }{2}}{25344}+\frac{395 \sin \frac{7 \varphi }{2}}{123552}+\frac{29 \sin \frac{11 \varphi }{2}}{13728}
\nonumber\\
\phi _{4,7,3}& =&-\frac{35 \sin \frac{\varphi }{2}}{33792}-\frac{3325 \sin \frac{5 \varphi }{2}}{988416}-\frac{179 \sin \frac{9 \varphi }{2}}{44928}-\frac{625 \sin \frac{13 \varphi }{2}}{350064}
\nonumber
\end{eqnarray}

\begin{eqnarray}
\phi _{6,7,0}& =&-\frac{\sin \frac{7 \varphi }{2}}{61776}
\nonumber\\
\phi _{6,7,1}& =&\frac{7 \sin \frac{5 \varphi }{2}}{247104}+\frac{23 \sin \frac{9 \varphi }{2}}{617760}
\nonumber
\end{eqnarray}

\subsection{Fifth singular exponent ($j=9$)}

\begin{eqnarray}
\phi _{0,9,0}& =&\sin \frac{9\varphi }{2}
\nonumber\\
\phi _{0,9,1}& =&-\frac{1}{4} \sin \frac{7 \varphi }{2}
\nonumber\\
\phi _{0,9,2}& =&\frac{3}{32} \sin \frac{5 \varphi }{2}+\frac{5}{44} \sin \frac{9 \varphi }{2}
\nonumber\\
\phi _{0,9,3}& =&-\frac{5}{128} \sin \frac{3 \varphi }{2}-\frac{15}{176} \sin \frac{7 \varphi }{2}-\frac{15}{286} \sin \frac{11 \varphi }{2}
\nonumber\\
\phi _{0,9,4}& =&\frac{35 \sin \frac{\varphi }{2}}{2048}+\frac{75 \sin \frac{5 \varphi }{2}}{1408}+\frac{135 \sin \frac{9 \varphi }{2}}{2288}+\frac{7}{286} \sin \frac{13 \varphi }{2}
\nonumber\\
\phi _{0,9,5}& =&\frac{63 \sin \frac{\varphi }{2}}{8192}-\frac{175 \sin \frac{3 \varphi }{2}}{5632}-\frac{225 \sin \frac{7 \varphi }{2}}{4576}-\frac{21}{572} \sin \frac{11 \varphi }{2}-\frac{28 \sin \frac{15 \varphi }{2}}{2431}
\nonumber\\
\phi _{0,9,6}& =&\frac{1575 \sin \frac{\varphi }{2}}{90112}-\frac{231 \sin \frac{3 \varphi }{2}}{65536}+\frac{2625 \sin \frac{5 \varphi }{2}}{73216}+\frac{175 \sin \frac{9 \varphi }{2}}{4576}+\frac{105 \sin \frac{13 \varphi }{2}}{4862}+\frac{252 \sin \frac{17 \varphi }{2}}{46189}
\nonumber
\end{eqnarray}

\begin{eqnarray}
\phi _{2,9,0}& =&-\frac{1}{22} \sin \frac{9 \varphi }{2}
\nonumber\\
\phi _{2,9,1}& =&\frac{3}{88} \sin \frac{7 \varphi }{2}+\frac{23}{572} \sin \frac{11 \varphi }{2}
\nonumber\\
\phi _{2,9,2}& =&-\frac{15}{704} \sin \frac{5 \varphi }{2}-\frac{103 \sin \frac{9 \varphi }{2}}{2288}-\frac{31 \sin \frac{13 \varphi }{2}}{1144}
\nonumber\\
\phi _{2,9,3}& =&\frac{35 \sin \frac{3 \varphi }{2}}{2816}+\frac{685 \sin \frac{7 \varphi }{2}}{18304}+\frac{553 \sin \frac{11 \varphi }{2}}{13728}+\frac{639 \sin \frac{15 \varphi }{2}}{38896}
\nonumber\\
\phi _{2,9,4}& =&-\frac{315 \sin \frac{\varphi }{2}}{45056}-\frac{1995 \sin \frac{5 \varphi }{2}}{73216}-\frac{5293 \sin \frac{9 \varphi }{2}}{126720}-\frac{14197 \sin \frac{13 \varphi }{2}}{466752}-\frac{13933 \sin \frac{17 \varphi }{2}}{1478048}
\nonumber
\end{eqnarray}

\begin{eqnarray}
\phi _{4,9,0}& =&\frac{\sin \frac{9 \varphi }{2}}{1144}
\nonumber\\
\phi _{4,9,1}& =&-\frac{5 \sin \frac{7 \varphi }{2}}{4576}-\frac{5 \sin \frac{11 \varphi }{2}}{3432}
\nonumber\\
\phi _{4,9,2}& =&\frac{35 \sin \frac{5 \varphi }{2}}{36608}+\frac{47 \sin \frac{9 \varphi }{2}}{20592}+\frac{361 \sin \frac{13 \varphi }{2}}{233376}
\nonumber
\end{eqnarray}

\begin{eqnarray}
\phi _{6,9,0}& =&-\frac{\sin \frac{9 \varphi }{2}}{102960}
\nonumber
\end{eqnarray}

\subsection{Sixth singular exponent ($j=11$)}
\begin{eqnarray}
\phi _{0,11,0}& =&\sin \frac{11\varphi }{2}
\nonumber\\
\phi _{0,11,1}& =&-\frac{1}{4} \sin \frac{9 \varphi }{2}
\nonumber\\
\phi _{0,11,2}& =&\frac{3}{32} \sin \frac{7 \varphi }{2}+\frac{3}{26} \sin \frac{11 \varphi }{2}
\nonumber\\
\phi _{0,11,3}& =&-\frac{5}{128} \sin \frac{5 \varphi }{2}-\frac{9}{104} \sin \frac{9 \varphi }{2}-\frac{7}{130} \sin \frac{13 \varphi }{2}
\nonumber\\
\phi _{0,11,4}& =&\frac{35 \sin \frac{3 \varphi }{2}}{2048}+\frac{45}{832} \sin \frac{7 \varphi }{2}+\frac{63 \sin \frac{11 \varphi }{2}}{1040}+\frac{28 \sin \frac{15 \varphi }{2}}{1105}
\nonumber\\
\phi _{0,11,5}& =&-\frac{63 \sin \frac{\varphi }{2}}{8192}-\frac{105 \sin \frac{5 \varphi }{2}}{3328}-\frac{21}{416} \sin \frac{9 \varphi }{2}-\frac{42 \sin \frac{13 \varphi }{2}}{1105}-\frac{252 \sin \frac{17 \varphi }{2}}{20995}
\nonumber
\end{eqnarray}

\begin{eqnarray}
\phi _{2,11,0}& =&-\frac{1}{26} \sin \frac{11 \varphi }{2}
\nonumber\\
\phi _{2,11,1}& =&\frac{3}{104} \sin \frac{9 \varphi }{2}+\frac{9}{260} \sin \frac{13 \varphi }{2}
\nonumber\\
\phi _{2,11,2}& =&-\frac{15}{832} \sin \frac{7 \varphi }{2}-\frac{121 \sin \frac{11 \varphi }{2}}{3120}-\frac{209 \sin \frac{15 \varphi }{2}}{8840}
\nonumber\\
\phi _{2,11,3}& =&\frac{35 \sin \frac{5 \varphi }{2}}{3328}+\frac{161 \sin \frac{9 \varphi }{2}}{4992}+\frac{3733 \sin \frac{13 \varphi }{2}}{106080}+\frac{4867 \sin \frac{17 \varphi }{2}}{335920}
\nonumber
\end{eqnarray}

\begin{eqnarray}
\phi _{4,11,0}& =&\frac{\sin \frac{11 \varphi }{2}}{1560}
\nonumber\\
\phi _{4,11,1}& =&-\frac{\sin \frac{9 \varphi }{2}}{1248}-\frac{29 \sin \frac{13 \varphi }{2}}{26520}
\nonumber
\end{eqnarray}

\subsection{Seventh singular exponent ($j=13$)}
\begin{eqnarray}
\phi _{0,13,0}& =&\sin \frac{13\varphi }{2}
\nonumber\\
\phi _{0,13,1}& =&-\frac{1}{4} \sin \frac{11 \varphi }{2}
\nonumber\\
\phi _{0,13,2}& =&\frac{3}{32} \sin \frac{9 \varphi }{2}+\frac{7}{60} \sin \frac{13 \varphi }{2}
\nonumber\\
\phi _{0,13,3}& =&-\frac{5}{128} \sin \frac{7 \varphi }{2}-\frac{7}{80} \sin \frac{11 \varphi }{2}-\frac{14}{255} \sin \frac{15 \varphi }{2}
\nonumber\\
\phi _{0,13,4}& =&\frac{35 \sin \frac{5 \varphi }{2}}{2048}+\frac{7}{128} \sin \frac{9 \varphi }{2}+\frac{21}{340} \sin \frac{13 \varphi }{2}+\frac{42 \sin \frac{17 \varphi }{2}}{1615}
\nonumber
\end{eqnarray}

\begin{eqnarray}
\phi _{2,13,0}& =&-\frac{1}{30} \sin \frac{13 \varphi }{2}
\nonumber\\
\phi _{2,13,1}& =&\frac{1}{40} \sin \frac{11 \varphi }{2}+\frac{31 \sin \frac{15 \varphi }{2}}{1020}
\nonumber\\
\phi _{2,13,2}& =&-\frac{1}{64} \sin \frac{9 \varphi }{2}-\frac{139 \sin \frac{13 \varphi }{2}}{4080}-\frac{271 \sin \frac{17 \varphi }{2}}{12920}
\nonumber
\end{eqnarray}

\begin{eqnarray}
\phi _{4,13,0}& =&\frac{\sin \frac{13 \varphi }{2}}{2040}
\nonumber
\end{eqnarray}


\subsection{Higher order singular exponent ($j=15,17,19,21$)}
\begin{eqnarray}
\phi _{0,15,0}& =&\sin \frac{15\varphi }{2}
\nonumber\\
\phi _{0,15,1}& =&-\frac{1}{4} \sin \frac{13 \varphi }{2}
\nonumber\\
\phi _{0,15,2}& =&\frac{3}{32} \sin \frac{11 \varphi }{2}+\frac{2}{17} \sin \frac{15 \varphi }{2}
\nonumber\\
\phi _{0,15,3}& =&-\frac{5}{128} \sin \frac{9 \varphi }{2}-\frac{3}{34} \sin \frac{13 \varphi }{2}-\frac{18}{323} \sin \frac{17 \varphi }{2}
\nonumber
\end{eqnarray}

\begin{eqnarray}
\phi _{2,15,0}& =&-\frac{1}{34} \sin \frac{15 \varphi }{2}
\nonumber\\
\phi _{2,15,1}& =&\frac{3}{136} \sin \frac{13 \varphi }{2}+\frac{35 \sin \frac{17 \varphi }{2}}{1292}
\nonumber
\end{eqnarray}


\begin{eqnarray}
\phi _{0,17,0}& =&\sin \frac{17\varphi }{2}
\nonumber\\
\phi _{0,17,1}& =&-\frac{1}{4} \sin \frac{15 \varphi }{2}
\nonumber\\
\phi _{0,17,2}& =&\frac{3}{32} \sin \frac{13 \varphi }{2}+\frac{9}{76} \sin \frac{17 \varphi }{2}
\nonumber
\end{eqnarray}

\begin{eqnarray}
\phi _{2,17,0}& =&\frac{1}{38} \sin \frac{17 \varphi }{2}
\nonumber
\end{eqnarray}


\begin{eqnarray}
\phi _{0,19,0}& =&\sin \frac{19\varphi }{2}
\nonumber\\
\phi _{0,19,1}& =&-\frac{1}{4} \sin \frac{17 \varphi }{2}
\nonumber
\end{eqnarray}


\begin{eqnarray}
\phi _{0,21,0}& =&\sin \frac{21\varphi }{2}
\nonumber
\end{eqnarray}

\newpage
\section{Dual singular functions and Dual-shadows for the crack $-\pi \le \varphi \le \pi$} \label{s.fnc dual_shadows_axi_crack}

\subsection{First singular exponent ($j=1$)}

\begin{eqnarray}
\psi _{0,1,0}& = &\sin \frac{\varphi }{2}
\nonumber\\
\psi _{0,1,1}& = &-\frac{1}{4} \sin \frac{3 \varphi }{2}
\nonumber\\
\psi _{0,1,2}& = &\frac{3}{32} \sin \frac{5 \varphi }{2}
\nonumber\\
\psi _{0,1,3}& = &-\frac{5}{128}\sin \frac{7 \varphi }{2}
\nonumber\\
\psi _{0,1,4}& = &\frac{35}{2048} \sin \frac{9 \varphi }{2}
\nonumber
\end{eqnarray}

\begin{eqnarray}
\psi _{2,1,0}& = &-\frac{1}{2} \sin \frac{\varphi }{2}
\nonumber\\
\psi _{2,1,1}& = &-\frac{1}{4} \sin \frac{\varphi }{2}+\frac{3}{8} \sin \frac{3 \varphi }{2}
\nonumber\\
\psi _{2,1,2}& = &-\frac{13}{48} \sin \frac{\varphi }{2}+\frac{1}{8} \sin \frac{3 \varphi }{2}-\frac{15}{64} \sin \frac{5 \varphi }{2}
\nonumber\\
\psi _{2,1,3}& = &-\frac{17}{96} \sin \frac{\varphi }{2}+\frac{85}{384} \sin \frac{3 \varphi }{2}-\frac{1}{16} \sin \frac{5 \varphi }{2}+\frac{35}{256} \sin \frac{7 \varphi }{2}
\nonumber\\
\psi _{2,1,4}& = &-\frac{2059 \sin \frac{\varphi }{2}}{11520}+\frac{7}{64} \sin \frac{3 \varphi }{2}-\frac{245 \sin \frac{5 \varphi }{2}}{1536}+\frac{1}{32} \sin \frac{7 \varphi }{2}-\frac{315 \sin \frac{9 \varphi }{2}}{4096}
\nonumber
\end{eqnarray}

\begin{eqnarray}
\psi _{4,1,0}& = &\frac{1}{24} \sin \frac{\varphi }{2}
\nonumber\\
\psi _{4,1,1}& = &\frac{1}{24} \sin \frac{\varphi }{2}-\frac{5}{96} \sin \frac{3 \varphi }{2}
\nonumber\\
\psi _{4,1,2}& = &\frac{19}{288} \sin \frac{\varphi }{2}-\frac{1}{32} \sin \frac{3 \varphi }{2}+\frac{35}{768} \sin \frac{5 \varphi }{2}
\nonumber\\
\psi _{4,1,3}& = &\frac{23}{384} \sin \frac{\varphi }{2}-\frac{161 \sin \frac{3 \varphi }{2}}{2304}+\frac{1}{48} \sin \frac{5 \varphi }{2}-\frac{35 \sin \frac{7 \varphi }{2}}{1024}
\nonumber\\
\psi _{4,1,4}& = &\frac{3431 \sin \frac{\varphi }{2}}{46080}-\frac{3}{64} \sin \frac{3 \varphi }{2}+\frac{63 \sin \frac{5 \varphi }{2}}{1024}-\frac{5}{384} \sin \frac{7 \varphi }{2}+\frac{385 \sin \frac{9 \varphi }{2}}{16384}
\nonumber
\end{eqnarray}

\newpage
\subsection{Second singular exponent ($j=3$)}
\begin{eqnarray}
\psi _{0,3,0}&=&\sin \left(\frac{3\varphi }{2}\right)
\nonumber\\
\psi _{0,3,1}&=&-\frac{1}{4}\sin \left(\frac{5 \varphi }{2}\right)
\nonumber\\
\psi _{0,3,2}&=&\frac{1}{4}\sin \left(\frac{3 \varphi }{2}\right)+\frac{3}{32}\sin \left(\frac{7 \varphi }{2}\right)
\nonumber\\
\psi _{0,3,3}&=&-\frac{3}{16}\sin \left(\frac{5 \varphi }{2}\right)-\frac{5}{128}\sin \left(\frac{9 \varphi }{2}\right)
\nonumber\\
\psi _{0,3,4}&=&\frac{15}{128}\sin \left(\frac{7 \varphi }{2}\right)+\frac{35}{2048}\sin \left(\frac{11 \varphi }{2}\right)
\nonumber
\end{eqnarray}

\begin{eqnarray}
\psi _{2,3,0}&=&\frac{1}{2}\sin \left(\frac{3 \varphi }{2}\right)
\nonumber\\
\psi _{2,3,1}&=&\frac{1}{4} \sin \left(\frac{\varphi }{2}\right)-\frac{3}{8} \sin \left(\frac{5 \varphi }{2}\right)
\nonumber\\
\psi _{2,3,2}&=&\frac{1}{8} \sin \left(\frac{\varphi }{2}\right)-\frac{5}{16} \sin \left(\frac{3 \varphi }{2}\right)+\frac{15}{64} \sin \left(\frac{7 \varphi }{2}\right)
\nonumber\\
\psi _{2,3,3}&=&\frac{19}{96} \sin \left(\frac{\varphi }{2}\right)-\frac{1}{16} \sin \left(\frac{3 \varphi }{2}\right)+\frac{35}{128} \sin \left(\frac{5 \varphi }{2}\right)-\frac{35}{256} \sin \left(\frac{9 \varphi }{2}\right)
\nonumber\\
\psi _{2,3,4}&=&\frac{23}{192} \sin \left(\frac{\varphi }{2}\right)-\frac{491 \sin \left(\frac{3 \varphi }{2}\right)}{2304}+\frac{1}{32} \sin \left(\frac{5 \varphi }{2}\right)-\frac{105}{512} \sin \left(\frac{7 \varphi }{2}\right)+\frac{315 \sin \left(\frac{11 \varphi }{2}\right)}{4096}
\nonumber
\end{eqnarray}

\begin{eqnarray}
\psi _{4,3,0}&=&-\frac{1}{8} \sin \left(\frac{3 \varphi }{2}\right)
\nonumber\\
\psi _{4,3,1}&=&\frac{1}{24} \sin \left(\frac{\varphi }{2}\right)+\frac{5}{32} \sin \left(\frac{5 \varphi }{2}\right)
\nonumber\\
\psi _{4,3,2}&=&\frac{1}{96} \sin \left(\frac{\varphi }{2}\right)-\frac{5}{72} \sin \left(\frac{3 \varphi }{2}\right)-\frac{35}{256} \sin \left(\frac{7 \varphi }{2}\right)
\nonumber\\
\psi _{4,3,3}&=&\frac{3}{128} \sin \left(\frac{\varphi }{2}\right)+\frac{175 \sin \left(\frac{5 \varphi }{2}\right)}{2304}+\frac{105 \sin \left(\frac{9 \varphi }{2}\right)}{1024}
\nonumber\\
\psi _{4,3,4}&=&\frac{1}{384} \sin \left(\frac{\varphi }{2}\right)-\frac{1477 \sin \left(\frac{3 \varphi }{2}\right)}{46080}-\frac{1}{384} \sin \left(\frac{5 \varphi }{2}\right)-\frac{35}{512} \sin \left(\frac{7 \varphi }{2}\right)-\frac{1155 \sin \left(\frac{11 \varphi }{2}\right)}{16384}
\nonumber
\end{eqnarray}

\newpage
\subsection{Third singular exponent ($j=5$)}
\begin{eqnarray}
\psi _{0,5,0}&=&\sin \left(\frac{5\varphi }{2}\right)
\nonumber\\
\psi _{0,5,1}&=&-\frac{1}{4}\sin \left(\frac{7 \varphi }{2}\right)
\nonumber\\
\psi _{0,5,2}&=& \frac{1}{6}\sin \left(\frac{5 \varphi }{2}\right)+\frac{3}{32}\sin \left(\frac{9 \varphi }{2}\right)
\nonumber\\
\psi _{0,5,3}&=&-\frac{1}{6}\sin \left(\frac{3 \varphi }{2}\right)-\frac{1}{8}\sin \left(\frac{7 \varphi }{2}\right)-\frac{5}{128}\sin \left(\frac{11 \varphi }{2}\right)
\nonumber\\
\psi _{0,5,4}&=&\frac{3}{16}\sin \left(\frac{5 \varphi }{2}\right)+\frac{5}{64}\sin \left(\frac{9 \varphi }{2}\right)+\frac{35}{2048}\sin \left(\frac{13 \varphi }{2}\right)
\nonumber
\end{eqnarray}

\begin{eqnarray}
\psi _{2,5,0}&=&\frac{1}{6} \sin \left(\frac{5 \varphi }{2}\right)
\nonumber\\
\psi _{2,5,1}&=&-\frac{5}{12} \sin \left(\frac{3 \varphi }{2}\right)-\frac{1}{8} \sin \left(\frac{7 \varphi }{2}\right)
\nonumber\\
\psi _{2,5,2}&=&-\frac{1}{8} \sin \left(\frac{\varphi }{2}\right)+\frac{23}{48} \sin \left(\frac{5 \varphi }{2}\right)+\frac{5}{64} \sin \left(\frac{9 \varphi }{2}\right)
\nonumber\\
\psi _{2,5,3}&=&-\frac{1}{16} \sin \left(\frac{\varphi }{2}\right)+\frac{7}{32} \sin \left(\frac{3 \varphi }{2}\right)-\frac{155}{384} \sin \left(\frac{7 \varphi }{2}\right)-\frac{35}{768} \sin \left(\frac{11 \varphi }{2}\right)
\nonumber\\
\psi _{2,5,4}&=&-\frac{25}{192} \sin \left(\frac{\varphi }{2}\right)+\frac{1}{32} \sin \left(\frac{3 \varphi }{2}\right)-\frac{559 \sin \left(\frac{5 \varphi }{2}\right)}{2304}+\frac{455 \sin \left(\frac{9 \varphi }{2}\right)}{1536}+\frac{105 \sin \left(\frac{13 \varphi }{2}\right)}{4096}
\nonumber
\end{eqnarray}

\begin{eqnarray}
\psi _{4,5,0}&=&\frac{1}{24} \sin \left(\frac{5 \varphi }{2}\right)
\nonumber\\
\psi _{4,5,1}&=&\frac{1}{8} \sin \left(\frac{3 \varphi }{2}\right)-\frac{5}{96} \sin \left(\frac{7 \varphi }{2}\right)
\nonumber\\
\psi _{4,5,2}&=&-\frac{5}{96} \sin \left(\frac{\varphi }{2}\right)-\frac{55}{288} \sin \left(\frac{5 \varphi }{2}\right)+\frac{35}{768} \sin \left(\frac{9 \varphi }{2}\right)
\nonumber\\
\psi _{4,5,3}&=&-\frac{1}{48} \sin \left(\frac{\varphi }{2}\right)+\frac{115 \sin \left(\frac{3 \varphi }{2}\right)}{1152}+\frac{455 \sin \left(\frac{7 \varphi }{2}\right)}{2304}-\frac{35 \sin \left(\frac{11 \varphi }{2}\right)}{1024}
\nonumber\\
\psi _{4,5,4}&=&-\frac{7}{144} \sin \left(\frac{\varphi }{2}\right)+\frac{1}{128} \sin \left(\frac{3 \varphi }{2}\right)-\frac{1141 \sin \left(\frac{5 \varphi }{2}\right)}{9216}-\frac{175 \sin \left(\frac{9 \varphi }{2}\right)}{1024}+\frac{385 \sin \left(\frac{13 \varphi }{2}\right)}{16384}
\nonumber
\end{eqnarray}
\normalsize

\clearpage
\section{Primal functions and shadows for $90^\circ$ V-notch $-\pi \le \varphi \le \pi/2$}
\subsection{First singular exponent ($j=1$)}\mbox{}

\tiny
\begin{eqnarray}
\phi_{0,1,0}& = &\sin \frac{2\varphi }{3}+\frac{1}{\sqrt{3}}\cos \frac{2\varphi }{3}
\nonumber\\
\phi _{0,1,1}& = &\frac{1}{4}\sin \frac{\varphi }{3}-\frac{1}{4 \sqrt{3}}\cos \frac{\varphi }{3}+\frac{1}{20} \sin \frac{5\varphi }{3}+\frac{1}{20 \sqrt{3}} \cos \frac{5\varphi }{3}
\nonumber\\
\phi _{0,1,2}& = &\frac{\sqrt{3}}{40}\cos \frac{2 \varphi }{3}+\frac{3}{40}\sin \frac{2 \varphi }{3}+\frac{\sqrt{3}}{32}\cos \frac{4 \varphi }{3}-\frac{3}{32}\sin \frac{4 \varphi }{3}
\nonumber\\
\phi _{0,1,3}& = &-\frac{13 \sqrt{3}}{640}\cos \frac{\varphi }{3}+\frac{39}{640}\sin \frac{\varphi }{3}-\frac{13 \sqrt{3}}{1280}\cos \frac{5 \varphi }{3}-\frac{39}{1280}\sin \frac{5 \varphi }{3}-\frac{5}{128 \sqrt{3}}\cos \frac{7 \varphi }{3}
\nonumber\\ & + &
\frac{5}{128}\sin \frac{7 \varphi }{3}+\frac{7}{1280} \sin \frac{11\varphi }{3}+\frac{7}{1280 \sqrt{3}} \cos \frac{11\varphi }{3}
\nonumber\\
\phi _{0,1,4}& = &\frac{247 \sqrt{3}}{20480}\cos \frac{2 \varphi }{3}+\frac{741}{20480}\sin \frac{2 \varphi }{3}+\frac{5}{128 \sqrt{3}}\cos \frac{4 \varphi }{3}-\frac{5}{128}\sin \frac{4 \varphi }{3}+
\frac{83}{7040 \sqrt{3}}\cos \frac{8 \varphi }{3}
\nonumber\\ & + &
\frac{83}{7040}\sin \frac{8 \varphi }{3}+\frac{35}{2048 \sqrt{3}}\cos \frac{10 \varphi }{3}-\frac{35}{2048}\sin \frac{10 \varphi }{3}
\nonumber\\
\phi _{0,1,5}& = &-\frac{169 \sqrt{3}}{16384}\cos \frac{\varphi }{3}+\frac{507}{16384}\sin \frac{\varphi }{3}-\frac{6021 \sqrt{3}}{901120}\cos \frac{5 \varphi }{3}-\frac{18063}{901120}\sin \frac{5 \varphi }{3}-\frac{63 \sqrt{3}}{8192}\cos \frac{7 \varphi }{3}
\nonumber\\ & + &
\frac{189}{8192}\sin \frac{7 \varphi }{3}-\frac{415}{78848 \sqrt{3}}\cos \frac{11 \varphi }{3}-\frac{415}{78848}\sin \frac{11 \varphi }{3}-\frac{21 \sqrt{3}}{8192}\cos \frac{13 \varphi }{3}+\frac{63}{8192}\sin \frac{13 \varphi }{3}
\nonumber\\ & + &
\frac{109}{157696 \sqrt{3}} \cos \frac{17\varphi }{3}+\frac{109}{157696} \sin \frac{17\varphi }{3}
\nonumber\\
\phi _{0,1,6}& = &\frac{2609 \sqrt{3}}{360448} \cos \frac{2 \varphi }{3}+\frac{1001 \sqrt{3}}{131072}\cos \frac{4 \varphi }{3}+\frac{36745 \sqrt{3}}{10092544} \cos \frac{8 \varphi }{3} + \frac{357 \sqrt{3}}{81920} \cos \frac{10 \varphi }{3}
\nonumber\\ & + &
\frac{1137 \sqrt{3}}{1531904}\cos \frac{14 \varphi }{3}+\frac{77 \sqrt{3}}{65536}\cos \frac{16 \varphi }{3}+\frac{7827}{360448} \sin \frac{2 \varphi }{3}-\frac{3003}{131072}\sin \frac{4 \varphi }{3}
\nonumber\\ & + &
\frac{110235}{10092544} \sin \frac{8 \varphi }{3}-\frac{1071}{81920}\sin \frac{10 \varphi }{3}+\frac{3411}{1531904} \sin \frac{14 \varphi }{3}-\frac{231}{65536} \sin \frac{16 \varphi }{3}
\nonumber\\
\phi _{0,1,7}& = &-\frac{18641 \sqrt{3}}{2883584} \cos \frac{\varphi }{3}-\frac{575075}{40370176 \sqrt{3}}\cos \frac{5 \varphi }{3}-\frac{1365 \sqrt{3}}{262144}\cos \frac{7 \varphi }{3}-\frac{1330111 \sqrt{3}}{686292992}\cos \frac{11 \varphi }{3}
\nonumber\\ & - &
\frac{3157 \sqrt{3}}{1310720}\cos \frac{13 \varphi }{3}-\frac{42069 \sqrt{3}}{122552320}\cos \frac{17 \varphi }{3}-\frac{143 \sqrt{3}}{262144} \cos \frac{19 \varphi }{3}+\frac{55923}{2883584} \sin \frac{\varphi }{3}
\nonumber\\ & - &
\frac{575075}{40370176}\sin \frac{5 \varphi }{3}+\frac{4095}{262144} \sin \frac{7 \varphi }{3}-\frac{3990333}{686292992} \sin \frac{11 \varphi }{3}+\frac{9471}{1310720} \sin \frac{13 \varphi }{3}
\nonumber\\ & - &
\frac{126207}{122552320} \sin \frac{17 \varphi }{3}+\frac{429}{262144} \sin \frac{19 \varphi }{3}+\frac{87217}{857866240} \sin \frac{23\varphi }{3}+\frac{87217}{857866240 \sqrt{3}} \cos \frac{23\varphi }{3}
\nonumber\\
\phi _{0,1,8}& = &\frac{1373975}{92274688 \sqrt{3}}\cos \frac{2 \varphi }{3}+\frac{74361 \sqrt{3}}{14417920}\cos \frac{4 \varphi }{3}+\frac{437105}{49020928 \sqrt{3}}\cos \frac{8 \varphi }{3}+\frac{141141 \sqrt{3}}{41943040}\cos \frac{10 \varphi }{3}
\nonumber\\ & + &
\frac{56281699 \sqrt{3}}{54903439360}\cos \frac{14 \varphi }{3}+\frac{429 \sqrt{3}}{327680}\cos \frac{16 \varphi }{3}+\frac{8995579}{19730923520 \sqrt{3}} \cos \frac{20 \varphi }{3}+\frac{2145 \sqrt{3}}{8388608}\cos \frac{22 \varphi }{3}
\nonumber\\ & + &
\frac{1373975}{92274688} \sin \frac{2 \varphi }{3}-\frac{223083}{14417920} \sin \frac{4 \varphi }{3}+\frac{437105}{49020928} \sin \frac{8 \varphi }{3}-\frac{423423}{41943040}\sin \frac{10 \varphi }{3}
\nonumber\\ & + &
\frac{168845097}{54903439360} \sin \frac{14 \varphi }{3}-\frac{1287}{327680} \sin \frac{16 \varphi }{3}+\frac{8995579}{19730923520} \sin \frac{20 \varphi }{3}-\frac{6435}{8388608} \sin \frac{22 \varphi }{3}
\nonumber\\
\phi _{0,1,9}& = &-\frac{1674225 \sqrt{3}}{369098752}\cos \frac{\varphi }{3}-\frac{22357035 \sqrt{3}}{6274678784} \cos \frac{5 \varphi }{3}-\frac{640381 \sqrt{3}}{167772160} \cos \frac{7 \varphi }{3}-\frac{79139337 \sqrt{3}}{43922751488} \cos \frac{11 \varphi }{3}
\nonumber\\ & - &
\frac{351351 \sqrt{3}}{167772160}\cos \frac{13 \varphi }{3}-\frac{2140689 \sqrt{3}}{3992977408} \cos \frac{17 \varphi }{3}-\frac{23595 \sqrt{3}}{33554432}\cos \frac{19 \varphi }{3}
\nonumber\\ & - &
\frac{62969053}{293145149440 \sqrt{3}}\cos \frac{23 \varphi }{3}-\frac{12155}{33554432 \sqrt{3}}\cos \frac{25 \varphi }{3}+\frac{5022675}{369098752} \sin \frac{\varphi }{3}-\frac{67071105}{6274678784} \sin \frac{5 \varphi }{3}
\nonumber\\ & + &
\frac{1921143}{167772160} \sin \frac{7 \varphi }{3}-\frac{237418011}{43922751488} \sin \frac{11 \varphi }{3}+\frac{1054053}{167772160} \sin \frac{13 \varphi }{3}-\frac{6422067}{3992977408} \sin \frac{17 \varphi }{3}
\nonumber\\ & + &
\frac{70785}{33554432} \sin \frac{19 \varphi }{3}-\frac{62969053}{293145149440} \sin \frac{23 \varphi }{3}+\frac{12155}{33554432} \sin \frac{25 \varphi }{3}+\frac{187981}{11152261120} \sin \frac{29 \varphi }{3}+\frac{187981}{11152261120 \sqrt{3}} \cos \frac{29 \varphi }{3}
\nonumber\\
\phi_{0,1,10}& = &\frac{2710107 \sqrt{3}}{738197504} \cos \frac{2 \varphi }{3}+\frac{507745 \sqrt{3}}{134217728}\cos \frac{4 \varphi }{3}+\frac{61145349 \sqrt{3}}{25098715136} \cos \frac{8 \varphi }{3}+\frac{898079 \sqrt{3}}{335544320} \cos \frac{10 \varphi }{3}
\nonumber\\ & + &
\frac{2152824343 \sqrt{3}}{2020446568448}\cos \frac{14 \varphi }{3}+\frac{340197 \sqrt{3}}{268435456} \cos \frac{16 \varphi }{3}+\frac{20970810791 \sqrt{3}}{75045158256640}\cos \frac{20 \varphi }{3}+\frac{75361}{67108864 \sqrt{3}} \cos \frac{22 \varphi }{3}
\nonumber\\ & + &
\frac{1787693753}{18310297026560 \sqrt{3}} \cos \frac{26 \varphi }{3}+\frac{46189}{268435456 \sqrt{3}} \cos \frac{28 \varphi }{3}+\frac{8130321}{738197504} \sin \frac{2 \varphi }{3}-\frac{1523235}{134217728} \sin \frac{4 \varphi }{3}+\frac{183436047}{25098715136} \sin \frac{8 \varphi }{3}
\nonumber\\ & - &
\frac{2694237}{335544320} \sin \frac{10 \varphi }{3}+\frac{6458473029}{2020446568448} \sin \frac{14 \varphi }{3}-\frac{1020591}{268435456} \sin \frac{16 \varphi }{3}+\frac{62912432373}{75045158256640} \sin \frac{20 \varphi }{3}-\frac{75361}{67108864} \sin \frac{22 \varphi }{3}
\nonumber\\ & + &
\frac{1787693753}{18310297026560} \sin \frac{26 \varphi }{3}-\frac{46189}{268435456} \sin \frac{28 \varphi }{3}
\nonumber
\end{eqnarray}
\newpage

\begin{eqnarray}
\nonumber\\
\phi _{2,1,0}& =&\frac{1}{20} \left(-\sqrt{3} \cos \frac{2 \varphi }{3}-3 \sin \frac{2 \varphi }{3}\right)
\nonumber\\
\phi _{2,1,1}& =&\frac{3}{80} \sqrt{3} \cos \frac{\varphi }{3}+\frac{11}{320} \sqrt{3} \cos \frac{5 \varphi }{3}+\frac{67 \sqrt{3} \cos \frac{11 \varphi }{3}}{3520}-\frac{9}{80} \sin \frac{\varphi }{3}+\frac{33}{320} \sin \frac{5 \varphi }{3}+\frac{201 \sin \frac{11 \varphi }{3}}{3520}
\nonumber\\
\phi _{2,1,2}& =&-\frac{99 \sqrt{3} \cos \frac{2 \varphi }{3}}{2560}-\frac{3}{128} \sqrt{3} \cos \frac{4 \varphi }{3}-\frac{87 \sqrt{3} \cos \frac{8 \varphi }{3}}{3520}-\frac{297 \sin \frac{2 \varphi }{3}}{2560}+\frac{9}{128} \sin \frac{4 \varphi }{3}-\frac{261 \sin \frac{8 \varphi }{3}}{3520}
\nonumber\\
\phi _{2,1,3}& =&\frac{33 \sqrt{3} \cos \frac{\varphi }{3}}{1024}+\frac{177 \sqrt{3} \cos \frac{5 \varphi }{3}}{5632}+\frac{7}{512} \sqrt{3} \cos \frac{7 \varphi }{3}+\frac{10207 \sqrt{3} \cos \frac{11 \varphi }{3}}{788480}+\frac{38427 \sqrt{3} \cos \frac{17 \varphi }{3}}{13404160}
\nonumber\\
&-&\frac{99 \sin \frac{\varphi }{3}}{1024}+\frac{531 \sin \frac{5 \varphi }{3}}{5632}-\frac{21}{512} \sin \frac{7 \varphi }{3}+\frac{30621 \sin \frac{11 \varphi }{3}}{788480}+\frac{115281 \sin \frac{17 \varphi }{3}}{13404160}
\nonumber\\
\phi _{2,1,4}& =&-\frac{28387 \sqrt{3} \cos \frac{2 \varphi }{3}}{901120}-\frac{385 \sqrt{3} \cos \frac{4 \varphi }{3}}{16384}-\frac{136329 \sqrt{3} \cos \frac{8 \varphi }{3}}{6307840}-\frac{63 \sqrt{3} \cos \frac{10 \varphi }{3}}{8192}-\frac{709 \sqrt{3} \cos \frac{14 \varphi }{3}}{95744}
\nonumber\\
&-&\frac{85161 \sin \frac{2 \varphi }{3}}{901120}+\frac{1155 \sin \frac{4 \varphi }{3}}{16384}-\frac{408987 \sin \frac{8 \varphi }{3}}{6307840}+\frac{189 \sin \frac{10 \varphi }{3}}{8192}-\frac{2127 \sin \frac{14 \varphi }{3}}{95744}
\nonumber\\
\phi _{2,1,5}& =&\frac{195993 \sqrt{3} \cos \frac{\varphi }{3}}{7208960}+\frac{163131 \sqrt{3} \cos \frac{5 \varphi }{3}}{6307840}+\frac{2079 \sqrt{3} \cos \frac{7 \varphi }{3}}{131072}+\frac{1200027 \sqrt{3} \cos \frac{11 \varphi }{3}}{85786624}+\frac{693 \sqrt{3} \cos \frac{13 \varphi }{3}}{163840}+\frac{4022021 \sqrt{3} \cos \frac{17 \varphi }{3}}{1072332800}
\nonumber\\
&+&\frac{11088127 \sqrt{3} \cos \frac{23 \varphi }{3}}{24663654400}-\frac{587979 \sin \frac{\varphi }{3}}{7208960}+\frac{489393 \sin \frac{5 \varphi }{3}}{6307840}-\frac{6237 \sin \frac{7 \varphi }{3}}{131072}+\frac{3600081 \sin \frac{11 \varphi }{3}}{85786624}-\frac{2079 \sin \frac{13 \varphi }{3}}{163840}+\frac{12066063 \sin \frac{17 \varphi }{3}}{1072332800}
\nonumber\\
&+&\frac{33264381 \sin \frac{23 \varphi }{3}}{24663654400}
\nonumber\\
\phi _{2,1,6}& =&-\frac{2611551 \sqrt{3} \cos \frac{2 \varphi }{3}}{100925440}-\frac{613557 \sqrt{3} \cos \frac{4 \varphi }{3}}{28835840}-\frac{9448749 \sqrt{3} \cos \frac{8 \varphi }{3}}{490209280}-\frac{53361 \sqrt{3} \cos \frac{10 \varphi }{3}}{5242880}-\frac{293357739 \sqrt{3} \cos \frac{14 \varphi }{3}}{34314649600}-\frac{3003 \sqrt{3} \cos \frac{16 \varphi }{3}}{1310720}
\nonumber\\
&-&\frac{1164771 \sqrt{3} \cos \frac{20 \varphi }{3}}{580321280}-\frac{7834653 \sin \frac{2 \varphi }{3}}{100925440}+\frac{1840671 \sin \frac{4 \varphi }{3}}{28835840}-\frac{28346247 \sin \frac{8 \varphi }{3}}{490209280}+\frac{160083 \sin \frac{10 \varphi }{3}}{5242880}-\frac{880073217 \sin \frac{14 \varphi }{3}}{34314649600}+\frac{9009 \sin \frac{16 \varphi }{3}}{1310720}
\nonumber\\
&-&\frac{3494313 \sin \frac{20 \varphi }{3}}{580321280}
\nonumber\\
\phi _{2,1,7}& =&\frac{37211697 \sqrt{3} \cos \frac{\varphi }{3}}{1614807040}+\frac{120244237 \sqrt{3} \cos \frac{5 \varphi }{3}}{5490343936}+\frac{81543 \sqrt{3} \cos \frac{7 \varphi }{3}}{5242880}+\frac{1843285679 \sqrt{3} \cos \frac{11 \varphi }{3}}{137258598400}+\frac{33033 \sqrt{3} \cos \frac{13 \varphi }{3}}{5242880}
\nonumber\\
&+&\frac{16043532597 \sqrt{3} \cos \frac{17 \varphi }{3}}{3156947763200}+\frac{1287 \sqrt{3} \cos \frac{19 \varphi }{3}}{1048576}+\frac{2067082063 \sqrt{3} \cos \frac{23 \varphi }{3}}{2052016046080}+\frac{22910746871 \sqrt{3} \cos \frac{29 \varphi }{3}}{297542326681600}-\frac{111635091 \sin \frac{\varphi }{3}}{1614807040}
\nonumber\\
&+&\frac{360732711 \sin \frac{5 \varphi }{3}}{5490343936}-\frac{244629 \sin \frac{7 \varphi }{3}}{5242880}+\frac{5529857037 \sin \frac{11 \varphi }{3}}{137258598400}-\frac{99099 \sin \frac{13 \varphi }{3}}{5242880}+\frac{48130597791 \sin \frac{17 \varphi }{3}}{3156947763200}
\nonumber\\
&-&\frac{3861 \sin \frac{19 \varphi }{3}}{1048576}+\frac{6201246189 \sin \frac{23 \varphi }{3}}{2052016046080}+\frac{68732240613 \sin \frac{29 \varphi }{3}}{297542326681600}
\nonumber\\
\phi _{2,1,8}& =&-\frac{23999498037 \sqrt{3} \cos \frac{2 \varphi }{3}}{1098068787200}-\frac{22191111 \sqrt{3} \cos \frac{4 \varphi }{3}}{1174405120}-\frac{18777125183 \sqrt{3} \cos \frac{8 \varphi }{3}}{1098068787200}-\frac{1813851 \sqrt{3} \cos \frac{10 \varphi }{3}}{167772160}-\frac{225549913639 \sqrt{3} \cos \frac{14 \varphi }{3}}{25255582105600}
\nonumber\\
&-&\frac{127413 \sqrt{3} \cos \frac{16 \varphi }{3}}{33554432}-\frac{964379391513 \sqrt{3} \cos \frac{20 \varphi }{3}}{328322567372800}-\frac{21879 \sqrt{3} \cos \frac{22 \varphi }{3}}{33554432}-\frac{260981373 \sqrt{3} \cos \frac{26 \varphi }{3}}{497562419200}-\frac{71998494111 \sin \frac{2 \varphi }{3}}{1098068787200}
\nonumber\\
&+&\frac{66573333 \sin \frac{4 \varphi }{3}}{1174405120}-\frac{56331375549 \sin \frac{8 \varphi }{3}}{1098068787200}+\frac{5441553 \sin \frac{10 \varphi }{3}}{167772160}-\frac{676649740917 \sin \frac{14 \varphi }{3}}{25255582105600}+\frac{382239 \sin \frac{16 \varphi }{3}}{33554432}
\nonumber\\
&-&\frac{2893138174539 \sin \frac{20 \varphi }{3}}{328322567372800}+\frac{65637 \sin \frac{22 \varphi }{3}}{33554432}-\frac{782944119 \sin \frac{26 \varphi }{3}}{497562419200}
\nonumber
\end{eqnarray}
\newpage

\begin{eqnarray}
\phi_{4,1,0}&=&\frac{3 \sqrt{3} \cos \frac{2 \varphi }{3}}{1280}+\frac{9 \sin \frac{2 \varphi }{3}}{1280}
\nonumber\\
\phi_{4,1,1}&=&-\frac{3 \sqrt{3} \cos \frac{\varphi }{3}}{1024}-\frac{171 \sqrt{3} \cos \frac{5 \varphi }{3}}{56320}-\frac{201 \sqrt{3} \cos \frac{11 \varphi }{3}}{197120}+\frac{1359 \sqrt{3} \cos \frac{17 \varphi }{3}}{3351040}+\frac{9 \sin \frac{\varphi }{3}}{1024}-\frac{513 \sin \frac{5 \varphi }{3}}{56320}-\frac{603 \sin \frac{11 \varphi }{3}}{197120}+\frac{4077 \sin \frac{17 \varphi }{3}}{3351040}
\nonumber\\
\phi _{4,1,2}&=&\frac{27 \sqrt{3} \cos \frac{2 \varphi }{3}}{5632}+\frac{21 \sqrt{3} \cos \frac{4 \varphi }{3}}{8192}+\frac{10881 \sqrt{3} \cos \frac{8 \varphi }{3}}{3153920}+\frac{75 \sqrt{3} \cos \frac{14 \varphi }{3}}{95744}+\frac{81 \sin \frac{2 \varphi }{3}}{5632}-\frac{63 \sin \frac{4 \varphi }{3}}{8192}+\frac{32643 \sin \frac{8 \varphi }{3}}{3153920}+\frac{225 \sin \frac{14 \varphi }{3}}{95744}
\nonumber\\
\phi _{4,1,3}&=&-\frac{1827 \sqrt{3} \cos \frac{\varphi }{3}}{360448}-\frac{14127 \sqrt{3} \cos \frac{5 \varphi }{3}}{2523136}-\frac{63 \sqrt{3} \cos \frac{7 \varphi }{3}}{32768}-\frac{127053 \sqrt{3} \cos \frac{11 \varphi }{3}}{42893312}-\frac{695133 \sqrt{3} \cos \frac{17 \varphi }{3}}{1072332800}+\frac{2098839 \sqrt{3} \cos \frac{23 \varphi }{3}}{24663654400}+\frac{5481 \sin \frac{\varphi }{3}}{360448}
\nonumber\\
&-&\frac{42381 \sin \frac{5 \varphi }{3}}{2523136}+\frac{189 \sin \frac{7 \varphi }{3}}{32768}-\frac{381159 \sin \frac{11 \varphi }{3}}{42893312}-\frac{2085399 \sin \frac{17 \varphi }{3}}{1072332800}+\frac{6296517 \sin \frac{23 \varphi }{3}}{24663654400}
\nonumber\\
\phi _{4,1,4}&=&\frac{47757 \sqrt{3} \cos \frac{2 \varphi }{3}}{7208960}+\frac{3213 \sqrt{3} \cos \frac{4 \varphi }{3}}{720896}+\frac{5907 \sqrt{3} \cos \frac{8 \varphi }{3}}{1114112}+\frac{693 \sqrt{3} \cos \frac{10 \varphi }{3}}{524288}+\frac{39272811 \sqrt{3} \cos \frac{14 \varphi }{3}}{17157324800}+\frac{100143 \sqrt{3} \cos \frac{20 \varphi }{3}}{246636544}+\frac{143271 \sin \frac{2 \varphi }{3}}{7208960}
\nonumber\\
&-&\frac{9639 \sin \frac{4 \varphi }{3}}{720896}+\frac{17721 \sin \frac{8 \varphi }{3}}{1114112}-\frac{2079 \sin \frac{10 \varphi }{3}}{524288}+\frac{117818433 \sin \frac{14 \varphi }{3}}{17157324800}+\frac{300429 \sin \frac{20 \varphi }{3}}{246636544}
\nonumber\\
\phi _{4,1,5}&=&-\frac{34371 \sqrt{3} \cos \frac{\varphi }{3}}{5242880}-\frac{6843897 \sqrt{3} \cos \frac{5 \varphi }{3}}{980418560}-\frac{7371 \sqrt{3} \cos \frac{7 \varphi }{3}}{2097152}-\frac{305405949 \sqrt{3} \cos \frac{11 \varphi }{3}}{68629299200}-\frac{9009 \sqrt{3} \cos \frac{13 \varphi }{3}}{10485760}-\frac{362650311 \sqrt{3} \cos \frac{17 \varphi }{3}}{225496268800}
\nonumber\\
&-&\frac{7307871 \sqrt{3} \cos \frac{23 \varphi }{3}}{28187033600}+\frac{4039929 \sqrt{3} \cos \frac{29 \varphi }{3}}{260089446400}+\frac{103113 \sin \frac{\varphi }{3}}{5242880}-\frac{20531691 \sin \frac{5 \varphi }{3}}{980418560}+\frac{22113 \sin \frac{7 \varphi }{3}}{2097152}-\frac{916217847 \sin \frac{11 \varphi }{3}}{68629299200}+\frac{27027 \sin \frac{13 \varphi }{3}}{10485760}
\nonumber\\
&-&\frac{1087950933 \sin \frac{17 \varphi }{3}}{225496268800}-\frac{21923613 \sin \frac{23 \varphi }{3}}{28187033600}+\frac{12119787 \sin \frac{29 \varphi }{3}}{260089446400}
\nonumber\\
\phi _{4,1,6}&=&\frac{10523907 \sqrt{3} \cos \frac{2 \varphi }{3}}{1372585984}+\frac{484713 \sqrt{3} \cos \frac{4 \varphi }{3}}{83886080}+\frac{101102037 \sqrt{3} \cos \frac{8 \varphi }{3}}{15686696960}+\frac{27027 \sqrt{3} \cos \frac{10 \varphi }{3}}{10485760}+\frac{491571873 \sqrt{3} \cos \frac{14 \varphi }{3}}{143497625600}+\frac{9009 \sqrt{3} \cos \frac{16 \varphi }{3}}{16777216}
\nonumber\\
&+&\frac{15993341307 \sqrt{3} \cos \frac{20 \varphi }{3}}{14923753062400}+\frac{3421914141 \sqrt{3} \cos \frac{26 \varphi }{3}}{22887871283200}+\frac{31571721 \sin \frac{2 \varphi }{3}}{1372585984}-\frac{1454139 \sin \frac{4 \varphi }{3}}{83886080}+\frac{303306111 \sin \frac{8 \varphi }{3}}{15686696960}-\frac{81081 \sin \frac{10 \varphi }{3}}{10485760}
\nonumber\\
&+&\frac{1474715619 \sin \frac{14 \varphi }{3}}{143497625600}-\frac{27027 \sin \frac{16 \varphi }{3}}{16777216}+\frac{47980023921 \sin \frac{20 \varphi }{3}}{14923753062400}+\frac{10265742423 \sin \frac{26 \varphi }{3}}{22887871283200}
\nonumber
\end{eqnarray}

\begin{eqnarray}
\phi _{6,1,0}&=&-\frac{3 \sqrt{3} \cos \frac{2 \varphi }{3}}{56320}-\frac{9 \sin \frac{2 \varphi }{3}}{56320}
\nonumber\\
\phi _{6,1,1}&=&\frac{21 \sqrt{3} \cos \frac{\varphi }{3}}{225280}+\frac{39 \sqrt{3} \cos \frac{5 \varphi }{3}}{394240}+\frac{603 \sqrt{3} \cos \frac{11 \varphi }{3}}{26808320}
\nonumber\\
&-&\frac{4077 \sqrt{3} \cos \frac{17 \varphi }{3}}{268083200}+\frac{21951 \sqrt{3} \cos \frac{23 \varphi }{3}}{6165913600}-\frac{63 \sin \frac{\varphi }{3}}{225280}+\frac{117 \sin \frac{5 \varphi }{3}}{394240}+\frac{1809 \sin \frac{11 \varphi }{3}}{26808320}-\frac{12231 \sin \frac{17 \varphi }{3}}{268083200}+\frac{65853 \sin \frac{23 \varphi }{3}}{6165913600}
\nonumber\\
\phi _{6,1,2}&=&-\frac{369 \sqrt{3} \cos \frac{2 \varphi }{3}}{1802240}-\frac{189 \sqrt{3} \cos \frac{4 \varphi }{3}}{1802240}-\frac{4401 \sqrt{3} \cos \frac{8 \varphi }{3}}{30638080}
\nonumber\\
&-&\frac{963 \sqrt{3} \cos \frac{14 \varphi }{3}}{38297600}+\frac{7911 \sqrt{3} \cos \frac{20 \varphi }{3}}{616591360}-\frac{1107 \sin \frac{2 \varphi }{3}}{1802240}+\frac{567 \sin \frac{4 \varphi }{3}}{1802240}-\frac{13203 \sin \frac{8 \varphi }{3}}{30638080}-\frac{2889 \sin \frac{14 \varphi }{3}}{38297600}+\frac{23733 \sin \frac{20 \varphi }{3}}{616591360}
\nonumber\\
\phi _{6,1,3}&=&\frac{1917 \sqrt{3} \cos \frac{\varphi }{3}}{7208960}+\frac{7371 \sqrt{3} \cos \frac{5 \varphi }{3}}{24510464}+\frac{63 \sqrt{3} \cos \frac{7 \varphi }{3}}{655360}+\frac{194433 \sqrt{3} \cos \frac{11 \varphi }{3}}{1225523200}+\frac{4882743 \sqrt{3} \cos \frac{17 \varphi }{3}}{197309235200}-\frac{4186251 \sqrt{3} \cos \frac{23 \varphi }{3}}{366431436800}
\nonumber\\
&-&\frac{17914173 \sqrt{3} \cos \frac{29 \varphi }{3}}{14877116334080}-\frac{5751 \sin \frac{\varphi }{3}}{7208960}+\frac{22113 \sin \frac{5 \varphi }{3}}{24510464}-\frac{189 \sin \frac{7 \varphi }{3}}{655360}+\frac{583299 \sin \frac{11 \varphi }{3}}{1225523200}+\frac{14648229 \sin \frac{17 \varphi }{3}}{197309235200}-\frac{12558753 \sin \frac{23 \varphi }{3}}{366431436800}
\nonumber\\
&-&\frac{53742519 \sin \frac{29 \varphi }{3}}{14877116334080}
\nonumber\\
\phi _{6,1,4}&=&-\frac{1041957 \sqrt{3} \cos \frac{2 \varphi }{3}}{2451046400}-\frac{1449 \sqrt{3} \cos \frac{4 \varphi }{3}}{5242880}-\frac{1717143 \sqrt{3} \cos \frac{8 \varphi }{3}}{4902092800}-\frac{819 \sqrt{3} \cos \frac{10 \varphi }{3}}{10485760}-\frac{238063431 \sqrt{3} \cos \frac{14 \varphi }{3}}{1578473881600}-\frac{50690457 \sqrt{3} \cos \frac{20 \varphi }{3}}{2565020057600}
\nonumber\\
&+&\frac{12531 \sqrt{3} \cos \frac{26 \varphi }{3}}{1541478400}-\frac{3125871 \sin \frac{2 \varphi }{3}}{2451046400}+\frac{4347 \sin \frac{4 \varphi }{3}}{5242880}-\frac{5151429 \sin \frac{8 \varphi }{3}}{4902092800}+\frac{2457 \sin \frac{10 \varphi }{3}}{10485760}-\frac{714190293 \sin \frac{14 \varphi }{3}}{1578473881600}-\frac{152071371 \sin \frac{20 \varphi }{3}}{2565020057600}
\nonumber\\
&+&\frac{37593 \sin \frac{26 \varphi }{3}}{1541478400}
\nonumber
\end{eqnarray}

\begin{eqnarray}
\phi _{8,1,0}&=&\frac{9 \sqrt{3} \cos \frac{2 \varphi }{3}}{12615680}+\frac{27 \sin \frac{2 \varphi }{3}}{12615680}
\nonumber\\
\phi _{8,1,1}&=&-\frac{81 \sqrt{3} \cos \frac{\varphi }{3}}{50462720}-\frac{1467 \sqrt{3} \cos \frac{5 \varphi }{3}}{857866240}-\frac{603 \sqrt{3} \cos \frac{11 \varphi }{3}}{2144665600}+\frac{12231 \sqrt{3} \cos \frac{17 \varphi }{3}}{49327308800}-\frac{65853 \sqrt{3} \cos \frac{23 \varphi }{3}}{641255014400}+\frac{44469 \sqrt{3} \cos \frac{29 \varphi }{3}}{2656627916800}+\frac{243 \sin \frac{\varphi }{3}}{50462720}
\nonumber\\
&-&\frac{4401 \sin \frac{5 \varphi }{3}}{857866240}-\frac{1809 \sin \frac{11 \varphi }{3}}{2144665600}+\frac{36693 \sin \frac{17 \varphi }{3}}{49327308800}-\frac{197559 \sin \frac{23 \varphi }{3}}{641255014400}+\frac{133407 \sin \frac{29 \varphi }{3}}{2656627916800}
\nonumber\\
\phi _{8,1,2}&=&\frac{7533 \sqrt{3} \cos \frac{2 \varphi }{3}}{1715732480}+\frac{81 \sqrt{3} \cos \frac{4 \varphi }{3}}{36700160}+\frac{20277 \sqrt{3} \cos \frac{8 \varphi }{3}}{6862929920}+\frac{2403 \sqrt{3} \cos \frac{14 \varphi }{3}}{7046758400}-\frac{405 \sqrt{3} \cos \frac{20 \varphi }{3}}{1165918208}+\frac{129519 \sqrt{3} \cos \frac{26 \varphi }{3}}{1430491955200}
\nonumber\\
&+&\frac{22599 \sin \frac{2 \varphi }{3}}{1715732480}-\frac{243 \sin \frac{4 \varphi }{3}}{36700160}+\frac{60831 \sin \frac{8 \varphi }{3}}{6862929920}+\frac{7209 \sin \frac{14 \varphi }{3}}{7046758400}-\frac{1215 \sin \frac{20 \varphi }{3}}{1165918208}+\frac{388557 \sin \frac{26 \varphi }{3}}{1430491955200}
\nonumber
\end{eqnarray}

\begin{eqnarray}
\phi _{10,1,0}&=&-\frac{27 \sqrt{3} \cos \frac{2 \varphi }{3}}{4289331200}-\frac{81 \sin \frac{2 \varphi }{3}}{4289331200}\nonumber
\end{eqnarray}

\normalsize
\subsection{Second singular exponent ($j=2$)}\mbox{}
\tiny
\begin{eqnarray}
\phi _{0,2,0} & = & \sin \frac{4 \varphi }{3}-\frac{1}{\sqrt{3}} \cos \frac{4 \varphi }{3}
\nonumber\\
\phi _{0,2,1} & = & \frac{1}{4 \sqrt{3}}\cos \frac{\varphi }{3}-\frac{1}{4}\sin \frac{\varphi }{3}-\frac{1}{28} \sin \frac{7\varphi }{3}+\frac{1}{28 \sqrt{3}} \cos \frac{7\varphi }{3}
\nonumber\\
\phi _{0,2,2} & = & -\frac{\sqrt{3}}{32}\cos \frac{2 \varphi }{3}-\frac{\sqrt{3}}{28} \cos \frac{4 \varphi }{3}-\frac{3}{32} \sin \frac{2 \varphi }{3}+\frac{3}{28} \sin \frac{4 \varphi }{3}
\nonumber\\
\phi _{0,2,3} & = & \frac{23 \sqrt{3}}{896} \cos \frac{\varphi }{3}+\frac{5}{128 \sqrt{3}} \cos \frac{5 \varphi }{3}+\frac{17 \sqrt{3}}{1120} \cos \frac{7 \varphi }{3}-\frac{69}{896} \sin \frac{\varphi }{3}+
\frac{5}{128} \sin \frac{5 \varphi }{3}
\nonumber\\ & - &
\frac{51}{1120}\sin \frac{7 \varphi }{3}-\frac{1}{280} \sin \frac{13\varphi }{3}+\frac{1}{280 \sqrt{3}} \cos \frac{13\varphi }{3}
\nonumber\\
\phi _{0,2,4} & = & -\frac{85}{1792 \sqrt{3}} \cos \frac{2 \varphi }{3}-\frac{17 \sqrt{3}}{1024} \cos \frac{4 \varphi }{3}-\frac{35}{2048 \sqrt{3}} \cos \frac{8 \varphi }{3}-\frac{35}{1664 \sqrt{3}} \cos \frac{10 \varphi }{3}
\nonumber\\ & - &
\frac{85}{1792} \sin \frac{2 \varphi }{3}+\frac{51}{1024} \sin \frac{4 \varphi }{3}+-\frac{35}{2048} \sin \frac{8 \varphi }{3}+\frac{35}{1664} \sin \frac{10 \varphi }{3}
\nonumber\\
\phi _{0,2,5} & = & \frac{391 \sqrt{3}}{28672} \cos \frac{\varphi }{3}+\frac{75 \sqrt{3}}{8192} \cos \frac{5 \varphi }{3}+\frac{531 \sqrt{3}}{53248} \cos \frac{7 \varphi }{3}+\frac{21 \sqrt{3}}{8192} \cos \frac{11 \varphi }{3}+\frac{1015}{106496 \sqrt{3}} \cos \frac{13 \varphi }{3}
\nonumber\\ & - &
\frac{1173}{28672} \sin \frac{\varphi }{3}+\frac{225}{8192} \sin \frac{5 \varphi }{3}-\frac{1593}{53248} \sin \frac{7 \varphi }{3}+\frac{63}{8192} \sin \frac{11 \varphi }{3}-\frac{1015}{106496} \sin \frac{13 \varphi }{3}
\nonumber\\ & - &
\frac{319}{745472} \sin \frac{19\varphi }{3}+\frac{319}{745472 \sqrt{3}} \cos \frac{19\varphi }{3}
\nonumber\\
\phi _{0,2,6} & = & -\frac{323 \sqrt{3}}{32768} \cos \frac{2 \varphi }{3}-\frac{3797 \sqrt{3}}{372736} \cos \frac{4 \varphi }{3}-\frac{21 \sqrt{3}}{4096} \cos \frac{8 \varphi }{3}-\frac{19459 \sqrt{3}}{3407872} \cos \frac{10 \varphi }{3}-\frac{77 \sqrt{3}}{65536} \cos \frac{14 \varphi }{3}
\nonumber\\ & - &
\frac{5307 \sqrt{3}}{3540992} \cos \frac{16 \varphi }{3}-\frac{969}{32768} \sin \frac{2 \varphi }{3}+\frac{11391}{372736} \sin \frac{4 \varphi }{3}-\frac{63}{4096} \sin \frac{8 \varphi }{3}+\frac{58377}{3407872} \sin \frac{10 \varphi }{3}
\nonumber\\ & - &
\frac{231}{65536} \sin \frac{14 \varphi }{3}+\frac{15921}{3540992} \sin \frac{16 \varphi }{3}
\nonumber\\
\phi _{0,2,7} & = & \frac{15023 \sqrt{3}}{1703936} \cos \frac{\varphi }{3}+\frac{867 \sqrt{3}}{131072} \cos \frac{5 \varphi }{3}+\frac{2002421}{95420416 \sqrt{3}} \cos \frac{7 \varphi }{3}
\nonumber\\ & + &
\frac{737 \sqrt{3}}{262144} \cos \frac{11 \varphi }{3}+\frac{5786399 \sqrt{3}}{1812987904} \cos \frac{13 \varphi }{3}+\frac{143 \sqrt{3}}{262144} \cos \frac{17 \varphi }{3}+\frac{217587 \sqrt{3}}{311607296} \cos \frac{19 \varphi }{3}
\nonumber\\ & - &
\frac{45069}{1703936}\sin \frac{\varphi }{3}+\frac{2601}{131072} \sin \frac{5 \varphi }{3}-\frac{2002421}{95420416} \sin \frac{7 \varphi }{3}+\frac{2211}{262144} \sin \frac{11 \varphi }{3}
\nonumber\\ & - &
\frac{17359197}{1812987904} \sin \frac{13 \varphi }{3}+\frac{429}{262144} \sin \frac{17 \varphi }{3}-\frac{652761}{311607296} \sin \frac{19 \varphi }{3}-\frac{4139}{68485120} \sin \frac{25\varphi }{3}
\nonumber\\ & + &
\frac{4139}{68485120 \sqrt{3}} \cos \frac{25\varphi }{3}
\nonumber\\
\phi _{0,2,8} & = & -\frac{7341 \sqrt{3}}{1064960} \cos \frac{2 \varphi }{3}-\frac{4616423}{218103808 \sqrt{3}} \cos \frac{4 \varphi }{3}-\frac{44319 \sqrt{3}}{10485760} \cos \frac{8 \varphi }{3}-\frac{49508365}{3625975808 \sqrt{3}} \cos \frac{10 \varphi }{3}
\nonumber\\ & - &
\frac{5577 \sqrt{3}}{3670016} \cos \frac{14 \varphi }{3}-\frac{279019063 \sqrt{3}}{159542935552} \cos \frac{16 \varphi }{3}-\frac{2145 \sqrt{3}}{8388608} \cos \frac{20 \varphi }{3}-\frac{11327483}{11331174400 \sqrt{3}} \cos \frac{22 \varphi }{3}
\nonumber\\ & - &
\frac{22023}{1064960} \sin \frac{2 \varphi }{3}+\frac{4616423}{218103808} \sin \frac{4 \varphi }{3}-\frac{132957}{10485760} \sin \frac{8 \varphi }{3}+\frac{49508365}{3625975808} \sin \frac{10 \varphi }{3}
\nonumber\\ & - &
\frac{16731}{3670016} \sin \frac{14 \varphi }{3}+\frac{837057189}{159542935552} \sin \frac{16 \varphi }{3}-\frac{6435}{8388608} \sin \frac{20 \varphi }{3}+\frac{11327483}{11331174400} \sin \frac{22 \varphi }{3}
\nonumber\\
\phi _{0,2,9} & = & \frac{5494833 \sqrt{3}}{872415232} \cos \frac{\varphi }{3}+\frac{42185 \sqrt{3}}{8388608} \cos \frac{5 \varphi }{3}+\frac{86863845 \sqrt{3}}{16575889408} \cos \frac{7 \varphi }{3}+\frac{21879 \sqrt{3}}{8388608} \cos \frac{11 \varphi }{3}
\nonumber\\ & + &
\frac{1818260415 \sqrt{3}}{638171742208} \cos \frac{13 \varphi }{3}+\frac{190905 \sqrt{3}}{234881024} \cos \frac{17 \varphi }{3}+\frac{15107532969 \sqrt{3}}{15954293555200} \cos \frac{19 \varphi }{3}
\nonumber\\ & + &
\frac{12155}{33554432 \sqrt{3}} \cos \frac{23 \varphi }{3}+\frac{600356599}{1269091532800 \sqrt{3}} \cos \frac{25 \varphi }{3}-\frac{16484499}{872415232} \sin \frac{\varphi }{3}+\frac{126555}{8388608} \sin \frac{5 \varphi }{3}
\nonumber\\ & - &
\frac{260591535}{16575889408} \sin \frac{7 \varphi }{3}+\frac{65637}{8388608} \sin \frac{11 \varphi }{3}-\frac{5454781245}{638171742208} \sin \frac{13 \varphi }{3}+\frac{572715}{234881024} \sin \frac{17 \varphi }{3}
\nonumber\\ & - &
\frac{45322598907}{15954293555200} \sin \frac{19 \varphi }{3}+\frac{12155}{33554432} \sin \frac{23 \varphi }{3}-\frac{600356599}{1269091532800}\sin \frac{25 \varphi }{3}
\nonumber\\ & - &
\frac{134768503}{13960006860800} \sin \frac{31\varphi }{3}+\frac{134768503}{13960006860800 \sqrt{3}} \cos \frac{31\varphi }{3}
\nonumber\\
\phi _{0,2,10} & = & -\frac{18025733 \sqrt{3} \cos \frac{2 \varphi }{3}}{3489660928}-\frac{87168519 \sqrt{3} \cos \frac{4 \varphi }{3}}{16575889408}-\frac{102421 \sqrt{3} \cos \frac{8 \varphi }{3}}{29360128}-\frac{1344919515 \sqrt{3} \cos \frac{10 \varphi }{3}}{364669566976}
\nonumber\\ & - &
\frac{736593 \sqrt{3} \cos \frac{14 \varphi }{3}}{469762048}-\frac{172871033 \sqrt{3} \cos \frac{16 \varphi }{3}}{99714334720}-\frac{303875 \cos \frac{20 \varphi }{3}}{234881024 \sqrt{3}}-\frac{906362008781 \sqrt{3} \cos \frac{22 \varphi }{3}}{1786880878182400}
\nonumber\\ & - &
\frac{46189 \cos \frac{26 \varphi }{3}}{268435456 \sqrt{3}}-\frac{7032308149 \cos \frac{28 \varphi }{3}}{30911443763200 \sqrt{3}}-\frac{54077199 \sin \frac{2 \varphi }{3}}{3489660928}+\frac{261505557 \sin \frac{4 \varphi }{3}}{16575889408}-\frac{307263 \sin \frac{8 \varphi }{3}}{29360128}
\nonumber\\ & + &
\frac{4034758545 \sin \frac{10 \varphi }{3}}{364669566976}-\frac{2209779 \sin \frac{14 \varphi }{3}}{469762048}+\frac{518613099 \sin \frac{16 \varphi }{3}}{99714334720}-\frac{303875 \sin \frac{20 \varphi }{3}}{234881024}
\nonumber\\ & + &
\frac{2719086026343 \sin \frac{22 \varphi }{3}}{1786880878182400}-\frac{46189 \sin \frac{26 \varphi }{3}}{268435456}+\frac{7032308149 \sin \frac{28 \varphi }{3}}{30911443763200}
\nonumber
\end{eqnarray}

\begin{eqnarray}
\phi _{2,2,0}&=&\frac{1}{28} \sqrt{3} \cos \frac{4 \varphi }{3}-\frac{3}{28} \sin \frac{4 \varphi }{3}
\nonumber\\
\phi _{2,2,1}&=&-\frac{3}{112} \sqrt{3} \cos \frac{\varphi }{3}-\frac{1}{35} \sqrt{3} \cos \frac{7 \varphi }{3}-\frac{97 \sqrt{3} \cos \frac{13 \varphi }{3}}{7280}+\frac{9}{112} \sin \frac{\varphi }{3}+\frac{3}{35} \sin \frac{7 \varphi }{3}+\frac{291 \sin \frac{13 \varphi }{3}}{7280}
\nonumber\\
\phi _{2,2,2}&=&\frac{15}{896} \sqrt{3} \cos \frac{2 \varphi }{3}+\frac{141 \sqrt{3} \cos \frac{4 \varphi }{3}}{4480}+\frac{123 \sqrt{3} \cos \frac{10 \varphi }{3}}{5824}+\frac{45}{896} \sin \frac{2 \varphi }{3}-\frac{423 \sin \frac{4 \varphi }{3}}{4480}-\frac{369 \sin \frac{10 \varphi }{3}}{5824}
\nonumber\\
\phi _{2,2,3}&=&-\frac{93 \sqrt{3} \cos \frac{\varphi }{3}}{3584}-\frac{5}{512} \sqrt{3} \cos \frac{5 \varphi }{3}-\frac{633 \sqrt{3} \cos \frac{7 \varphi }{3}}{23296}
\nonumber\\
&-&\frac{21839 \sqrt{3} \cos \frac{13 \varphi }{3}}{1863680}-\frac{68787 \sqrt{3} \cos \frac{19 \varphi }{3}}{35409920}+\frac{279 \sin \frac{\varphi }{3}}{3584}-\frac{15}{512} \sin \frac{5 \varphi }{3}+\frac{1899 \sin \frac{7 \varphi }{3}}{23296}+\frac{65517 \sin \frac{13 \varphi }{3}}{1863680}+\frac{206361 \sin \frac{19 \varphi }{3}}{35409920}
\nonumber\\
\phi _{2,2,4}&=&\frac{77 \sqrt{3} \cos \frac{2 \varphi }{3}}{4096}+\frac{50497 \sqrt{3} \cos \frac{4 \varphi }{3}}{1863680}+\frac{45 \sqrt{3} \cos \frac{8 \varphi }{3}}{8192}
\nonumber\\
&+&\frac{3207 \sqrt{3} \cos \frac{10 \varphi }{3}}{163840}+\frac{751 \sqrt{3} \cos \frac{16 \varphi }{3}}{110656}+\frac{231 \sin \frac{2 \varphi }{3}}{4096}-\frac{151491 \sin \frac{4 \varphi }{3}}{1863680}+\frac{135 \sin \frac{8 \varphi }{3}}{8192}-\frac{9621 \sin \frac{10 \varphi }{3}}{163840}-\frac{2253 \sin \frac{16 \varphi }{3}}{110656}
\nonumber\\
\phi _{2,2,5}&=&-\frac{49509 \sqrt{3} \cos \frac{\varphi }{3}}{2129920}-\frac{207 \sqrt{3} \cos \frac{5 \varphi }{3}}{16384}-\frac{21369 \sqrt{3} \cos \frac{7 \varphi }{3}}{917504}-\frac{99 \sqrt{3} \cos \frac{11 \varphi }{3}}{32768}-\frac{1137957 \sqrt{3} \cos \frac{13 \varphi }{3}}{87162880}-\frac{22170977 \sqrt{3} \cos \frac{19 \varphi }{3}}{6232145920}
\nonumber\\
&-&\frac{9624067 \sqrt{3} \cos \frac{25 \varphi }{3}}{31160729600}+\frac{148527 \sin \frac{\varphi }{3}}{2129920}-\frac{621 \sin \frac{5 \varphi }{3}}{16384}+\frac{64107 \sin \frac{7 \varphi }{3}}{917504}-\frac{297 \sin \frac{11 \varphi }{3}}{32768}+\frac{3413871 \sin \frac{13 \varphi }{3}}{87162880}+\frac{66512931 \sin \frac{19 \varphi }{3}}{6232145920}
\nonumber\\
&+&\frac{28872201 \sin \frac{25 \varphi }{3}}{31160729600}
\nonumber\\
\phi _{2,2,6}&=&\frac{5931 \sqrt{3} \cos \frac{2 \varphi }{3}}{327680}+\frac{2762703 \sqrt{3} \cos \frac{4 \varphi }{3}}{119275520}+\frac{10593 \sqrt{3} \cos \frac{8 \varphi }{3}}{1310720}+\frac{32529369 \sqrt{3} \cos \frac{10 \varphi }{3}}{1812987904}+\frac{429 \sqrt{3} \cos \frac{14 \varphi }{3}}{262144}
\nonumber\\
&+&\frac{816384879 \sqrt{3} \cos \frac{16 \varphi }{3}}{99714334720}+\frac{54482031 \sqrt{3} \cos \frac{22 \varphi }{3}}{28327936000}+\frac{17793 \sin \frac{2 \varphi }{3}}{327680}-\frac{8288109 \sin \frac{4 \varphi }{3}}{119275520}+\frac{31779 \sin \frac{8 \varphi }{3}}{1310720}-\frac{97588107 \sin \frac{10 \varphi }{3}}{1812987904}+\frac{1287 \sin \frac{14 \varphi }{3}}{262144}
\nonumber\\
&-&\frac{2449154637 \sin \frac{16 \varphi }{3}}{99714334720}-\frac{163446093 \sin \frac{22 \varphi }{3}}{28327936000}
\nonumber\\
\phi _{2,2,7}&=&-\frac{78308991 \sqrt{3} \cos \frac{\varphi }{3}}{3816816640}-\frac{449031 \sqrt{3} \cos \frac{5 \varphi }{3}}{34078720}-\frac{1470654989 \sqrt{3} \cos \frac{7 \varphi }{3}}{72519516160}-\frac{26169 \sqrt{3} \cos \frac{11 \varphi }{3}}{5242880}-\frac{1279456777 \sqrt{3} \cos \frac{13 \varphi }{3}}{99714334720}-\frac{6435 \sqrt{3} \cos \frac{17 \varphi }{3}}{7340032}
\nonumber\\
&-&\frac{4498216551 \sqrt{3} \cos \frac{19 \varphi }{3}}{906493952000}-\frac{34493067223 \sqrt{3} \cos \frac{25 \varphi }{3}}{34900017152000}-\frac{58410950531 \sqrt{3} \cos \frac{31 \varphi }{3}}{1081900531712000}+\frac{234926973 \sin \frac{\varphi }{3}}{3816816640}-\frac{1347093 \sin \frac{5 \varphi }{3}}{34078720}
\nonumber\\
&+&\frac{4411964967 \sin \frac{7 \varphi }{3}}{72519516160}-\frac{78507 \sin \frac{11 \varphi }{3}}{5242880}+\frac{3838370331 \sin \frac{13 \varphi }{3}}{99714334720}-\frac{19305 \sin \frac{17 \varphi }{3}}{7340032}+\frac{13494649653 \sin \frac{19 \varphi }{3}}{906493952000}+\frac{103479201669 \sin \frac{25 \varphi }{3}}{34900017152000}
\nonumber\\
&+&\frac{175232851593 \sin \frac{31 \varphi }{3}}{1081900531712000}
\nonumber\\
\phi _{2,2,8}&=&\frac{102271257 \sqrt{3} \cos \frac{2 \varphi }{3}}{6106906624}+\frac{641736027 \sqrt{3} \cos \frac{4 \varphi }{3}}{31876710400}+\frac{2680227 \sqrt{3} \cos \frac{8 \varphi }{3}}{293601280}+\frac{10376443733 \sqrt{3} \cos \frac{10 \varphi }{3}}{638171742208}+\frac{176319 \sqrt{3} \cos \frac{14 \varphi }{3}}{58720256}
\nonumber\\
&+&\frac{36582771103 \sqrt{3} \cos \frac{16 \varphi }{3}}{4198498304000}+\frac{109395 \sqrt{3} \cos \frac{20 \varphi }{3}}{234881024}+\frac{3259134726429 \sqrt{3} \cos \frac{22 \varphi }{3}}{1116800548864000}+\frac{10028501487 \sqrt{3} \cos \frac{28 \varphi }{3}}{19319652352000}+\frac{306813771 \sin \frac{2 \varphi }{3}}{6106906624}
\nonumber\\
&-&\frac{1925208081 \sin \frac{4 \varphi }{3}}{31876710400}+\frac{8040681 \sin \frac{8 \varphi }{3}}{293601280}-\frac{31129331199 \sin \frac{10 \varphi }{3}}{638171742208}+\frac{528957 \sin \frac{14 \varphi }{3}}{58720256}-\frac{109748313309 \sin \frac{16 \varphi }{3}}{4198498304000}+\frac{328185 \sin \frac{20 \varphi }{3}}{234881024}
\nonumber\\
&-&\frac{9777404179287 \sin \frac{22 \varphi }{3}}{1116800548864000}-\frac{30085504461 \sin \frac{28 \varphi }{3}}{19319652352000}
\nonumber
\end{eqnarray}
\newpage

\begin{eqnarray}
\phi _{4,2,0}&=&-\frac{3 \sqrt{3} \cos \frac{4 \varphi }{3}}{2240}+\frac{9 \sin \frac{4 \varphi }{3}}{2240}
\nonumber\\
\phi _{4,2,1}&=&\frac{3 \sqrt{3} \cos \frac{\varphi }{3}}{1792}+\frac{45 \sqrt{3} \cos \frac{7 \varphi }{3}}{23296}+\frac{291 \sqrt{3} \cos \frac{13 \varphi }{3}}{465920}-\frac{1737 \sqrt{3} \cos \frac{19 \varphi }{3}}{8852480}-\frac{9 \sin \frac{\varphi }{3}}{1792}-\frac{135 \sin \frac{7 \varphi }{3}}{23296}-\frac{873 \sin \frac{13 \varphi }{3}}{465920}+\frac{5211 \sin \frac{19 \varphi }{3}}{8852480}
\nonumber\\
\phi _{4,2,2}&=&-\frac{3 \sqrt{3} \cos \frac{2 \varphi }{3}}{2048}-\frac{141 \sqrt{3} \cos \frac{4 \varphi }{3}}{46592}-\frac{17181 \sqrt{3} \cos \frac{10 \varphi }{3}}{7454720}-\frac{111 \sqrt{3} \cos \frac{16 \varphi }{3}}{221312}-\frac{9 \sin \frac{2 \varphi }{3}}{2048}+\frac{423 \sin \frac{4 \varphi }{3}}{46592}+\frac{51543 \sin \frac{10 \varphi }{3}}{7454720}+\frac{333 \sin \frac{16 \varphi }{3}}{221312}
\nonumber\\
\phi _{4,2,3}&=&\frac{339 \sqrt{3} \cos \frac{\varphi }{3}}{106496}+\frac{9 \sqrt{3} \cos \frac{5 \varphi }{3}}{8192}+\frac{22401 \sqrt{3} \cos \frac{7 \varphi }{3}}{5963776}+\frac{1161297 \sqrt{3} \cos \frac{13 \varphi }{3}}{566558720}+\frac{1288377 \sqrt{3} \cos \frac{19 \varphi }{3}}{3116072960}-\frac{731631 \sqrt{3} \cos \frac{25 \varphi }{3}}{15580364800}-\frac{1017 \sin \frac{\varphi }{3}}{106496}
\nonumber\\
&+&\frac{27 \sin \frac{5 \varphi }{3}}{8192}-\frac{67203 \sin \frac{7 \varphi }{3}}{5963776}-\frac{3483891 \sin \frac{13 \varphi }{3}}{566558720}-\frac{3865131 \sin \frac{19 \varphi }{3}}{3116072960}+\frac{2194893 \sin \frac{25 \varphi }{3}}{15580364800}
\nonumber\\
\phi _{4,2,4}&=&-\frac{297 \sqrt{3} \cos \frac{2 \varphi }{3}}{106496}-\frac{75477 \sqrt{3} \cos \frac{4 \varphi }{3}}{17039360}-\frac{99 \sqrt{3} \cos \frac{8 \varphi }{3}}{131072}-\frac{210081 \sqrt{3} \cos \frac{10 \varphi }{3}}{56655872}-\frac{11511609 \sqrt{3} \cos \frac{16 \varphi }{3}}{7122452480}-\frac{48771 \sqrt{3} \cos \frac{22 \varphi }{3}}{186368000}-\frac{891 \sin \frac{2 \varphi }{3}}{106496}
\nonumber\\
&+&\frac{226431 \sin \frac{4 \varphi }{3}}{17039360}-\frac{297 \sin \frac{8 \varphi }{3}}{131072}+\frac{630243 \sin \frac{10 \varphi }{3}}{56655872}+\frac{34534827 \sin \frac{16 \varphi }{3}}{7122452480}+\frac{146313 \sin \frac{22 \varphi }{3}}{186368000}
\nonumber\\
\phi _{4,2,5}&=&\frac{595161 \sqrt{3} \cos \frac{\varphi }{3}}{136314880}+\frac{14949 \sqrt{3} \cos \frac{5 \varphi }{3}}{6815744}+\frac{12637593 \sqrt{3} \cos \frac{7 \varphi }{3}}{2589982720}+\frac{1287 \sqrt{3} \cos \frac{11 \varphi }{3}}{2621440}+\frac{638879919 \sqrt{3} \cos \frac{13 \varphi }{3}}{199428669440}+\frac{88487757 \sqrt{3} \cos \frac{19 \varphi }{3}}{76703334400}
\nonumber\\
&+&\frac{2921545173 \sqrt{3} \cos \frac{25 \varphi }{3}}{17450008576000}-\frac{1085706519 \sqrt{3} \cos \frac{31 \varphi }{3}}{108190053171200}-\frac{1785483 \sin \frac{\varphi }{3}}{136314880}+\frac{44847 \sin \frac{5 \varphi }{3}}{6815744}-\frac{37912779 \sin \frac{7 \varphi }{3}}{2589982720}+\frac{3861 \sin \frac{11 \varphi }{3}}{2621440}
\nonumber\\
&-&\frac{1916639757 \sin \frac{13 \varphi }{3}}{199428669440}-\frac{265463271 \sin \frac{19 \varphi }{3}}{76703334400}-\frac{8764635519 \sin \frac{25 \varphi }{3}}{17450008576000}+\frac{3257119557 \sin \frac{31 \varphi }{3}}{108190053171200}
\nonumber\\
\phi _{4,2,6}&=&-\frac{8363223 \sqrt{3} \cos \frac{2 \varphi }{3}}{2181038080}-\frac{48435399 \sqrt{3} \cos \frac{4 \varphi }{3}}{9064939520}-\frac{1683 \sqrt{3} \cos \frac{8 \varphi }{3}}{1048576}-\frac{212164605 \sqrt{3} \cos \frac{10 \varphi }{3}}{45583695872}-\frac{1287 \sqrt{3} \cos \frac{14 \varphi }{3}}{4194304}-\frac{25127548689 \sqrt{3} \cos \frac{16 \varphi }{3}}{9971433472000}
\nonumber\\
&-&\frac{652815981 \sqrt{3} \cos \frac{22 \varphi }{3}}{839699660800}-\frac{1069062327 \sqrt{3} \cos \frac{28 \varphi }{3}}{11039801344000}-\frac{25089669 \sin \frac{2 \varphi }{3}}{2181038080}+\frac{145306197 \sin \frac{4 \varphi }{3}}{9064939520}-\frac{5049 \sin \frac{8 \varphi }{3}}{1048576}+\frac{636493815 \sin \frac{10 \varphi }{3}}{45583695872}
\nonumber\\
&-&\frac{3861 \sin \frac{14 \varphi }{3}}{4194304}+\frac{75382646067 \sin \frac{16 \varphi }{3}}{9971433472000}+\frac{1958447943 \sin \frac{22 \varphi }{3}}{839699660800}+\frac{3207186981 \sin \frac{28 \varphi }{3}}{11039801344000}
\nonumber
\end{eqnarray}

\begin{eqnarray}
\phi _{6,2,0}&=&\frac{3 \sqrt{3} \cos \frac{4 \varphi }{3}}{116480}-\frac{9 \sin \frac{4 \varphi }{3}}{116480}
\nonumber\\
\phi _{6,2,1}&=&-\frac{3 \sqrt{3} \cos \frac{\varphi }{3}}{66560}-\frac{3 \sqrt{3} \cos \frac{7 \varphi }{3}}{57344}-\frac{873 \sqrt{3} \cos \frac{13 \varphi }{3}}{70819840}+\frac{5211 \sqrt{3} \cos \frac{19 \varphi }{3}}{779018240}+\frac{9 \sin \frac{\varphi }{3}}{66560}+\frac{9 \sin \frac{7 \varphi }{3}}{57344}+\frac{2619 \sin \frac{13 \varphi }{3}}{70819840}-\frac{15633 \sin \frac{19 \varphi }{3}}{779018240}
\nonumber\\
\phi _{6,2,2}&=&\frac{27 \sqrt{3} \cos \frac{2 \varphi }{3}}{532480}+\frac{459 \sqrt{3} \cos \frac{4 \varphi }{3}}{4259840}+\frac{9153 \sqrt{3} \cos \frac{10 \varphi }{3}}{113311744}+\frac{1503 \sqrt{3} \cos \frac{16 \varphi }{3}}{97377280}-\frac{5211 \sqrt{3} \cos \frac{22 \varphi }{3}}{856064000}+\frac{81 \sin \frac{2 \varphi }{3}}{532480}-\frac{1377 \sin \frac{4 \varphi }{3}}{4259840}-\frac{27459 \sin \frac{10 \varphi }{3}}{113311744}
\nonumber\\
&-&\frac{4509 \sin \frac{16 \varphi }{3}}{97377280}+\frac{15633 \sin \frac{22 \varphi }{3}}{856064000}
\nonumber\\
\phi _{6,2,3}&=&-\frac{297 \sqrt{3} \cos \frac{\varphi }{3}}{2129920}-\frac{99 \sqrt{3} \cos \frac{5 \varphi }{3}}{2129920}-\frac{2871 \sqrt{3} \cos \frac{7 \varphi }{3}}{17039360}-\frac{1167711 \sqrt{3} \cos \frac{13 \varphi }{3}}{12464291840}-\frac{1231299 \sqrt{3} \cos \frac{19 \varphi }{3}}{77901824000}+\frac{23785659 \sqrt{3} \cos \frac{25 \varphi }{3}}{4362502144000}
\nonumber\\
&+&\frac{57409293 \sqrt{3} \cos \frac{31 \varphi }{3}}{67618783232000}+\frac{891 \sin \frac{\varphi }{3}}{2129920}-\frac{297 \sin \frac{5 \varphi }{3}}{2129920}+\frac{8613 \sin \frac{7 \varphi }{3}}{17039360}+\frac{3503133 \sin \frac{13 \varphi }{3}}{12464291840}+\frac{3693897 \sin \frac{19 \varphi }{3}}{77901824000}-\frac{71356977 \sin \frac{25 \varphi }{3}}{4362502144000}
\nonumber\\
&-&\frac{172227879 \sin \frac{31 \varphi }{3}}{67618783232000}
\nonumber\\
\phi _{6,2,4}&=&\frac{19701 \sqrt{3} \cos \frac{2 \varphi }{3}}{136314880}+\frac{767421 \sqrt{3} \cos \frac{4 \varphi }{3}}{3237478400}+\frac{99 \sqrt{3} \cos \frac{8 \varphi }{3}}{2621440}+\frac{293115 \sqrt{3} \cos \frac{10 \varphi }{3}}{1424490496}+\frac{35369211 \sqrt{3} \cos \frac{16 \varphi }{3}}{383516672000}+\frac{65479527 \sqrt{3} \cos \frac{22 \varphi }{3}}{4985716736000}
\nonumber\\
&-&\frac{306636441 \sqrt{3} \cos \frac{28 \varphi }{3}}{77278609408000}+\frac{59103 \sin \frac{2 \varphi }{3}}{136314880}-\frac{2302263 \sin \frac{4 \varphi }{3}}{3237478400}+\frac{297 \sin \frac{8 \varphi }{3}}{2621440}-\frac{879345 \sin \frac{10 \varphi }{3}}{1424490496}-\frac{106107633 \sin \frac{16 \varphi }{3}}{383516672000}-\frac{196438581 \sin \frac{22 \varphi }{3}}{4985716736000}
\nonumber\\
&+&\frac{919909323 \sin \frac{28 \varphi }{3}}{77278609408000}
\nonumber
\end{eqnarray}

\begin{eqnarray}
\phi _{8,2,0}&=&-\frac{9 \sqrt{3} \cos \frac{4 \varphi }{3}}{29818880}+\frac{27 \sin \frac{4 \varphi }{3}}{29818880}
\nonumber\\
\phi _{8,2,1}&=&\frac{81 \sqrt{3} \cos \frac{\varphi }{3}}{119275520}+\frac{1773 \sqrt{3} \cos \frac{7 \varphi }{3}}{2266234880}+\frac{873 \sqrt{3} \cos \frac{13 \varphi }{3}}{6232145920}-\frac{15633 \sqrt{3} \cos \frac{19 \varphi }{3}}{155803648000}-\frac{243 \sin \frac{\varphi }{3}}{119275520}-\frac{5319 \sin \frac{7 \varphi }{3}}{2266234880}-\frac{2619 \sin \frac{13 \varphi }{3}}{6232145920}+\frac{46899 \sin \frac{19 \varphi }{3}}{155803648000}
\nonumber\\
\phi _{8,2,2}&=&-\frac{891 \sqrt{3} \cos \frac{2 \varphi }{3}}{954204160}-\frac{81 \sqrt{3} \cos \frac{4 \varphi }{3}}{40468480}-\frac{2619 \sqrt{3} \cos \frac{10 \varphi }{3}}{1812987904}-\frac{4023 \sqrt{3} \cos \frac{16 \varphi }{3}}{19475456000}+\frac{1516401 \sqrt{3} \cos \frac{22 \varphi }{3}}{8725004288000}-\frac{2673 \sin \frac{2 \varphi }{3}}{954204160}+\frac{243 \sin \frac{4 \varphi }{3}}{40468480}
\nonumber\\
&+&\frac{7857 \sin \frac{10 \varphi }{3}}{1812987904}+\frac{12069 \sin \frac{16 \varphi }{3}}{19475456000}-\frac{4549203 \sin \frac{22 \varphi }{3}}{8725004288000}
\nonumber
\end{eqnarray}

\begin{eqnarray}
\phi _{10,2,0}&=&\frac{27 \sqrt{3} \cos \frac{4 \varphi }{3}}{11331174400}-\frac{81 \sin \frac{4 \varphi }{3}}{11331174400}
\nonumber
\end{eqnarray}

\normalsize
\subsection{Higher order exponents ($j=3,\ldots,11$)}\mbox{}
\tiny
\begin{eqnarray}
\phi _{0,3,0} & = &\cos  2\varphi
\nonumber\\
\phi _{0,3,1}& = &-\frac{\cos  \varphi }{4}-\frac{1}{12} \cos  3 \varphi
\nonumber\\
\phi _{0,3,2}& = &\frac{3}{32}+\frac{1}{8} \cos  2 \varphi
\nonumber\\
\phi _{0,3,3}& = &-\frac{\cos  \varphi }{8}-\frac{7}{128} \cos  3 \varphi -\frac{1}{128} \cos  5 \varphi
\nonumber\\
\phi _{0,3,4}& = &\frac{5}{96}+\frac{77 \cos  2 \varphi }{1024}+\frac{17}{640} \cos  4 \varphi
\nonumber\\
\phi _{0,3,5}& = &-\frac{315 \cos  \varphi }{4096}-\frac{897 \cos  3 \varphi }{20480}-\frac{187 \cos  5 \varphi }{15360}-\frac{7 \cos  7 \varphi }{7680}
\nonumber\\
\phi _{0,3,6}& = &\frac{2205}{65536}+\frac{433 \cos  2 \varphi }{8192}+\frac{4037 \cos  4 \varphi }{163840}+\frac{843 \cos  6 \varphi }{143360}
\nonumber\\
\phi _{0,3,7}& = &-\frac{875 \cos  \varphi }{16384}-\frac{4521 \cos  3 \varphi }{131072}-\frac{62597 \cos  5 \varphi }{4587520}-\frac{2529 \cos  7 \varphi }{917504}-\frac{349 \cos  9 \varphi }{2752512}
\nonumber\\
\phi _{0,3,8}& = &\frac{1575}{65536}+\frac{248941 \cos  2 \varphi }{6291456}+\frac{59729 \cos  4 \varphi }{2752512}+\frac{218139 \cos  6 \varphi }{29360128}+\frac{1835 \cos  8 \varphi }{1376256}
\nonumber\\
\phi _{0,3,9}& = &-\frac{1670823 \cos  \varphi }{41943040}-\frac{233745 \cos  3 \varphi }{8388608}-\frac{1558145 \cos  5 \varphi }{117440512}-\frac{1412281 \cos  7 \varphi }{352321536}
\nonumber\\ & - &
\frac{6973 \cos  9 \varphi }{11010048}-\frac{2201 \cos  11 \varphi }{110100480}
\nonumber
\end{eqnarray}

\begin{eqnarray}
\phi _{2,3,0}& = &-\frac{1}{12} \cos  2 \varphi
\nonumber\\
\phi _{2,3,1}& = &\frac{\cos  \varphi }{16}+\frac{7}{96} \cos  3 \varphi +\frac{1}{32} \cos  5 \varphi
\nonumber\\
\phi _{2,3,2}& = &-\frac{5}{128}-\frac{61}{768} \cos  2 \varphi -\frac{53}{960} \cos  4 \varphi
\nonumber\\
\phi _{2,3,3}& = &\frac{45 \cos  \varphi }{512}+\frac{137 \cos  3 \varphi }{1920}+\frac{487 \cos  5 \varphi }{15360}+\frac{71 \cos  7 \varphi }{15360}
\nonumber\\
\phi _{2,3,4}& = &-\frac{385}{8192}-\frac{15437 \cos  2 \varphi }{184320}-\frac{19523 \cos  4 \varphi }{368640}-\frac{143 \cos  6 \varphi }{7680}
\nonumber\\
\phi _{2,3,5}& = &\frac{67949 \cos  \varphi }{737280}+\frac{309119 \cos  3 \varphi }{4423680}+\frac{372203 \cos  5 \varphi }{10321920}+\frac{2447 \cos  7 \varphi }{245760}+\frac{2365 \cos  9 \varphi }{3096576}
\nonumber\\
\phi _{2,3,6}& = &-\frac{61733}{1310720}-\frac{648079 \cos  2 \varphi }{7864320}-\frac{2200529 \cos  4 \varphi }{41287680}-\frac{844089 \cos  6 \varphi }{36700160}-\frac{168421 \cos  8 \varphi }{30965760}
\nonumber\\
\phi _{2,3,7}& = &\frac{2806877 \cos  \varphi }{31457280}+\frac{6434681 \cos  3 \varphi }{94371840}+\frac{10053427 \cos  5 \varphi }{264241152}+\frac{55988167 \cos  7 \varphi }{3963617280}+\frac{2806673 \cos  9 \varphi }{990904320}+\frac{301339 \cos  11 \varphi }{2179989504}
\nonumber
\end{eqnarray}

\begin{eqnarray}
\phi _{4,3,0}& = &\frac{1}{384} \cos  2 \varphi
\nonumber\\
\phi _{4,3,1}& = &-\frac{5 \cos  \varphi }{1536}-\frac{31 \cos  3 \varphi }{7680}-\frac{1}{768} \cos  5 \varphi +\frac{3 \cos  7 \varphi }{8960}
\nonumber\\
\phi _{4,3,2}& = &\frac{35}{12288}+\frac{29 \cos  2 \varphi }{4608}+\frac{923 \cos  4 \varphi }{184320}+\frac{29 \cos  6 \varphi }{26880}
\nonumber\\
\phi _{4,3,3}& = &-\frac{161 \cos  \varphi }{18432}-\frac{3611 \cos  3 \varphi }{442368}-\frac{23633 \cos  5 \varphi }{5160960}-\frac{3077 \cos  7 \varphi }{3440640}+\frac{9103 \cos  9 \varphi }{92897280}
\nonumber\\
\phi _{4,3,4}& = &\frac{189}{32768}+\frac{391699 \cos  2 \varphi }{35389440}+\frac{17161 \cos  4 \varphi }{2064384}+\frac{67513 \cos  6 \varphi }{18350080}+\frac{1325 \cos  8 \varphi }{2322432}
\nonumber\\
\phi _{4,3,5}& = &-\frac{73213 \cos  \varphi }{5242880}-\frac{1681349 \cos  3 \varphi }{141557760}-\frac{138697 \cos  5 \varphi }{18874368}-\frac{5281387 \cos  7 \varphi }{1981808640}-\frac{683371 \cos  9 \varphi }{1857945600}+\frac{334729 \cos  11 \varphi }{13624934400}
\nonumber
\end{eqnarray}

\begin{eqnarray}
\phi _{6,3,0}& = &-\frac{\cos  2 \varphi }{23040}
\nonumber\\
\phi _{6,3,1}& = &\frac{7 \cos  \varphi }{92160}+\frac{13 \cos  3 \varphi }{138240}+\frac{\cos  5 \varphi }{43008}-\frac{3 \cos  7 \varphi }{286720}+\frac{43 \cos  9 \varphi }{23224320}
\nonumber\\
\phi _{6,3,2}& = &-\frac{7}{81920}-\frac{427 \cos  2 \varphi }{2211840}-\frac{787 \cos  4 \varphi }{5160960}-\frac{9 \cos  6 \varphi }{286720}+\frac{53 \cos  8 \varphi }{5806080}
\nonumber\\
\phi _{6,3,3}& = &\frac{161 \cos  \varphi }{491520}+\frac{2801 \cos  3 \varphi }{8847360}+\frac{3035 \cos  5 \varphi }{16515072}+\frac{2767 \cos  7 \varphi }{82575360}-\frac{1571 \cos  9 \varphi }{176947200}-\frac{7361 \cos  11 \varphi }{4541644800}
\nonumber
\end{eqnarray}

\begin{eqnarray}
\phi _{8,3,0}& = &\frac{\cos  2 \varphi }{2211840}
\nonumber\\
\phi_{8,3,1}& = &-\frac{\cos  \varphi }{983040}-\frac{11 \cos  3 \varphi }{8847360}-\frac{\cos  5 \varphi }{4128768}+\frac{\cos  7 \varphi }{6881280}-\frac{43 \cos  9 \varphi }{928972800}+\frac{\cos  11 \varphi }{162201600}
\nonumber
\end{eqnarray}

\begin{eqnarray}
\phi _{0,4,0}& = &\sin \frac{8}{3} \varphi +\frac{1}{\sqrt{3}} \cos \frac{8}{3} \varphi
\nonumber\\
\phi _{0,4,1}& = &-\frac{1}{4 \sqrt{3}}\cos \frac{5 \varphi }{3}-\frac{1}{4}\sin \frac{5 \varphi }{3}-\frac{5}{44} \sin \frac{11\varphi }{3}-\frac{5}{44 \sqrt{3}} \cos \frac{11\varphi }{3}
\nonumber\\
\phi _{0,4,2}& = &\frac{\sqrt{3}}{32}\cos \frac{2 \varphi }{3}+\frac{3}{32}\sin \frac{2 \varphi }{3}+\frac{\sqrt{3}}{22}\cos \frac{8 \varphi }{3}+\frac{3}{22}\sin \frac{8 \varphi }{3}
\nonumber\\
\phi _{0,4,3}& = &-\frac{5}{128 \sqrt{3}}\cos \frac{\varphi }{3}+\frac{5}{128}\sin \frac{\varphi }{3}-\frac{43 \sqrt{3}}{1408}\cos \frac{5 \varphi }{3}-\frac{129}{1408}\sin \frac{5 \varphi }{3}
\nonumber\\ & - &
\frac{25 \sqrt{3}}{1232}\cos \frac{11 \varphi }{3}-\frac{75}{1232}\sin \frac{11 \varphi }{3}-\frac{25}{2464}\sin \frac{17\varphi }{3}-\frac{25}{2464 \sqrt{3}}\cos \frac{17\varphi }{3}
\nonumber\\
\phi _{0,4,4}& = &\frac{155}{2816 \sqrt{3}}\cos \frac{2 \varphi }{3}+\frac{155}{2816}\sin \frac{2 \varphi }{3}+\frac{35}{2048 \sqrt{3}}\cos \frac{4 \varphi }{3}-\frac{35}{2048}\sin \frac{4 \varphi }{3}+\frac{1675 \sqrt{3}}{78848}\cos \frac{8 \varphi }{3}
\nonumber\\ & + &
\frac{5025}{78848}\sin \frac{8 \varphi }{3}+\frac{725}{23936 \sqrt{3}}\cos \frac{14 \varphi }{3}+\frac{725}{23936}\sin \frac{14 \varphi }{3}
\nonumber\\
\phi _{0,4,5}& = &-\frac{945 \sqrt{3}}{90112}\cos \frac{\varphi }{3}+\frac{2835}{90112}\sin \frac{\varphi }{3}-\frac{5375 \sqrt{3}}{315392}\cos \frac{5 \varphi }{3}-\frac{16125 }{315392}\sin \frac{5 \varphi }{3}-\frac{21 \sqrt{3}}{8192}\cos \frac{7 \varphi }{3}
\nonumber\\ & + &
\frac{63}{8192}\sin \frac{7 \varphi }{3}-\frac{72225 \sqrt{3}}{5361664}\cos \frac{11 \varphi }{3}-\frac{216675}{5361664}\sin \frac{11 \varphi }{3}-\frac{5365}{382976 \sqrt{3}}\cos \frac{17 \varphi }{3}
\nonumber\\ & - &
\frac{5365}{382976}\sin \frac{17 \varphi }{3}-\frac{3133}{2680832} \sin \frac{23\varphi }{3}-\frac{3133}{2680832 \sqrt{3}}\cos \frac{23\varphi }{3}
\nonumber\\
\phi _{0,4,6}& = &\frac{4375 \sqrt{3}}{360448}\cos \frac{2 \varphi }{3}+\frac{13125}{360448}\sin \frac{2 \varphi }{3}+\frac{525 \sqrt{3}}{90112}\cos \frac{4 \varphi }{3}-\frac{1575}{90112}\sin \frac{4 \varphi }{3}+\frac{5125 \sqrt{3}}{382976}\cos \frac{8 \varphi }{3}
\nonumber\\ & + &
\frac{15375}{382976}\sin \frac{8 \varphi }{3}+\frac{77 \sqrt{3}}{65536}\cos \frac{10 \varphi }{3}-\frac{231}{65536}\sin \frac{10 \varphi }{3}+\frac{684685 \sqrt{3}}{85786624}\cos \frac{14 \varphi }{3}
\nonumber\\ & + &
\frac{2054055}{85786624}\sin \frac{14 \varphi }{3}+\frac{140577 \sqrt{3}}{61659136}\cos \frac{20 \varphi }{3}+\frac{421731}{61659136}\sin \frac{20 \varphi }{3}
\nonumber\\
\phi _{0,4,7}& = &-\frac{11625 \sqrt{3} \cos \frac{\varphi }{3}}{1441792}-\frac{25375 \sqrt{3} \cos \frac{5 \varphi }{3}}{2228224}-\frac{833 \sqrt{3} \cos \frac{7 \varphi }{3}}{262144}-\frac{9794825 \cos \frac{11 \varphi }{3}}{343146496 \sqrt{3}}
\nonumber\\ & - &
\frac{143 \sqrt{3} \cos \frac{13 \varphi }{3}}{262144}-\frac{36189151 \sqrt{3} \cos \frac{17 \varphi }{3}}{7892369408}-\frac{984039 \sqrt{3} \cos \frac{23 \varphi }{3}}{916078592}-\frac{129905 \cos \frac{29 \varphi }{3}}{801568768 \sqrt{3}}
\nonumber\\ & + &
\frac{34875 \sin \frac{\varphi }{3}}{1441792}-\frac{76125 \sin \frac{5 \varphi }{3}}{2228224}+\frac{2499 \sin \frac{7 \varphi }{3}}{262144}-\frac{9794825 \sin \frac{11 \varphi }{3}}{343146496}
\nonumber\\ & + &
\frac{429 \sin \frac{13 \varphi }{3}}{262144}-\frac{108567453 \sin \frac{17 \varphi }{3}}{7892369408}-\frac{2952117 \sin \frac{23 \varphi }{3}}{916078592}-\frac{129905 \sin \frac{29 \varphi }{3}}{801568768}
\nonumber\\
\phi _{0,4,8}& = &\frac{53925 \sqrt{3} \cos \frac{2 \varphi }{3}}{6127616}+\frac{10725 \sqrt{3} \cos \frac{4 \varphi }{3}}{2097152}+\frac{1305175 \cos \frac{8 \varphi }{3}}{46137344 \sqrt{3}}+\frac{897 \sqrt{3} \cos \frac{10 \varphi }{3}}{524288}
\nonumber\\ & + &
\frac{301061615 \cos \frac{14 \varphi }{3}}{15784738816 \sqrt{3}}+\frac{2145 \sqrt{3} \cos \frac{16 \varphi }{3}}{8388608}+\frac{300557929 \sqrt{3} \cos \frac{20 \varphi }{3}}{117258059776}+\frac{89746385 \cos \frac{26 \varphi }{3}}{57219678208 \sqrt{3}}
\nonumber\\ & + &
\frac{161775 \sin \frac{2 \varphi }{3}}{6127616}-\frac{32175 \sin \frac{4 \varphi }{3}}{2097152}+\frac{1305175 \sin \frac{8 \varphi }{3}}{46137344}-\frac{2691 \sin \frac{10 \varphi }{3}}{524288}+\frac{301061615 \sin \frac{14 \varphi }{3}}{15784738816}
\nonumber\\ & - &
\frac{6435 \sin \frac{16 \varphi }{3}}{8388608}+\frac{901673787 \sin \frac{20 \varphi }{3}}{117258059776}+\frac{89746385 \sin \frac{26 \varphi }{3}}{57219678208}
\nonumber
\end{eqnarray}

\begin{eqnarray}
\phi _{2,4,0}& = &-\frac{1}{44} \sqrt{3} \cos \frac{8 \varphi }{3}-\frac{3}{44} \sin \frac{8 \varphi }{3}
\nonumber\\
\phi _{2,4,1}& = &\frac{3}{176} \sqrt{3} \cos \frac{5 \varphi }{3}+\frac{13}{616} \sqrt{3} \cos \frac{11 \varphi }{3}+\frac{9}{176} \sin \frac{5 \varphi }{3}+\frac{39}{616} \sin \frac{11 \varphi }{3}
\nonumber\\
\phi _{2,4,2}& = &-\frac{15 \sqrt{3} \cos \frac{2 \varphi }{3}}{1408}-\frac{225 \sqrt{3} \cos \frac{8 \varphi }{3}}{9856}-\frac{37 \sqrt{3} \cos \frac{14 \varphi }{3}}{2618}-\frac{45 \sin \frac{2 \varphi }{3}}{1408}-\frac{675 \sin \frac{8 \varphi }{3}}{9856}-\frac{111 \sin \frac{14 \varphi }{3}}{2618}
\nonumber\\
\phi_{2,4,3}& = &\frac{35 \sqrt{3} \cos \frac{\varphi }{3}}{5632}+\frac{105 \sqrt{3} \cos \frac{5 \varphi }{3}}{5632}+\frac{13599 \sqrt{3} \cos \frac{11 \varphi }{3}}{670208}+\frac{5751 \sqrt{3} \cos \frac{17 \varphi }{3}}{670208}-\frac{105 \sin \frac{\varphi }{3}}{5632}+\frac{315 \sin \frac{5 \varphi }{3}}{5632}+\frac{40797 \sin \frac{11 \varphi }{3}}{670208}+\frac{17253 \sin \frac{17 \varphi }{3}}{670208}
\nonumber\\
\phi _{2,4,4}& = &-\frac{55 \sqrt{3} \cos \frac{2 \varphi }{3}}{4096}-\frac{315 \sqrt{3} \cos \frac{4 \varphi }{3}}{90112}-\frac{110427 \sqrt{3} \cos \frac{8 \varphi }{3}}{5361664}-\frac{163539 \sqrt{3} \cos \frac{14 \varphi }{3}}{10723328}-\frac{302101 \sqrt{3} \cos \frac{20 \varphi }{3}}{61659136}-\frac{165 \sin \frac{2 \varphi }{3}}{4096}+\frac{945 \sin \frac{4 \varphi }{3}}{90112}
\nonumber\\
&-&\frac{331281 \sin \frac{8 \varphi }{3}}{5361664}-\frac{490617 \sin \frac{14 \varphi }{3}}{10723328}-\frac{906303 \sin \frac{20 \varphi }{3}}{61659136}
\nonumber\\
\phi _{2,4,5}& = &\frac{405 \sqrt{3} \cos \frac{\varphi }{3}}{45056}+\frac{53 \sqrt{3} \cos \frac{5 \varphi }{3}}{2992}+\frac{63 \sqrt{3} \cos \frac{7 \varphi }{3}}{32768}+\frac{990783 \sqrt{3} \cos \frac{11 \varphi }{3}}{53616640}+\frac{10275141 \sqrt{3} \cos \frac{17 \varphi }{3}}{986546176}+\frac{155323 \sqrt{3} \cos \frac{23 \varphi }{3}}{57254912}-\frac{1215 \sin \frac{\varphi }{3}}{45056}
\nonumber\\
&+&\frac{159 \sin \frac{5 \varphi }{3}}{2992}-\frac{189 \sin \frac{7 \varphi }{3}}{32768}+\frac{2972349 \sin \frac{11 \varphi }{3}}{53616640}+\frac{30825423 \sin \frac{17 \varphi }{3}}{986546176}+\frac{465969 \sin \frac{23 \varphi }{3}}{57254912}
\nonumber\\
\phi _{2,4,6}& = &-\frac{337743 \sqrt{3} \cos \frac{2 \varphi }{3}}{24510464}-\frac{1503 \sqrt{3} \cos \frac{4 \varphi }{3}}{262144}-\frac{15848761 \sqrt{3} \cos \frac{8 \varphi }{3}}{857866240}-\frac{273 \sqrt{3} \cos \frac{10 \varphi }{3}}{262144}-\frac{579169167 \sqrt{3} \cos \frac{14 \varphi }{3}}{39461847040}-\frac{687402189 \sqrt{3} \cos \frac{20 \varphi }{3}}{102600802304}
\nonumber\\
&-&\frac{38987053 \sqrt{3} \cos \frac{26 \varphi }{3}}{26566279168}-\frac{1013229 \sin \frac{2 \varphi }{3}}{24510464}+\frac{4509 \sin \frac{4 \varphi }{3}}{262144}-\frac{47546283 \sin \frac{8 \varphi }{3}}{857866240}+\frac{819 \sin \frac{10 \varphi }{3}}{262144}-\frac{1737507501 \sin \frac{14 \varphi }{3}}{39461847040}-\frac{2062206567 \sin \frac{20 \varphi }{3}}{102600802304}
\nonumber\\
&-&\frac{116961159 \sin \frac{26 \varphi }{3}}{26566279168}
\nonumber
\end{eqnarray}

\begin{eqnarray}
\phi _{4,4,0}& = &\frac{3 \sqrt{3} \cos \frac{8 \varphi }{3}}{4928}+\frac{9 \sin \frac{8 \varphi }{3}}{4928}
\nonumber\\
\phi _{4,4,1}& = &-\frac{15 \sqrt{3} \cos \frac{5 \varphi }{3}}{19712}-\frac{333 \sqrt{3} \cos \frac{11 \varphi }{3}}{335104}+\frac{597 \sqrt{3} \cos \frac{23 \varphi }{3}}{1926848}-\frac{45 \sin \frac{5 \varphi }{3}}{19712}-\frac{999 \sin \frac{11 \varphi }{3}}{335104}+\frac{1791 \sin \frac{23 \varphi }{3}}{1926848}
\nonumber\\
\phi _{4,4,2}& = &\frac{15 \sqrt{3} \cos \frac{2 \varphi }{3}}{22528}+\frac{1035 \sqrt{3} \cos \frac{8 \varphi }{3}}{670208}+\frac{5517 \sqrt{3} \cos \frac{14 \varphi }{3}}{5361664}-\frac{597 \sqrt{3} \cos \frac{20 \varphi }{3}}{7707392}+\frac{45 \sin \frac{2 \varphi }{3}}{22528}+\frac{3105 \sin \frac{8 \varphi }{3}}{670208}+\frac{16551 \sin \frac{14 \varphi }{3}}{5361664}-\frac{1791 \sin \frac{20 \varphi }{3}}{7707392}
\nonumber\\
\phi _{4,4,3}& = &-\frac{45 \sqrt{3} \cos \frac{\varphi }{3}}{90112}-\frac{225 \sqrt{3} \cos \frac{5 \varphi }{3}}{139264}-\frac{41163 \sqrt{3} \cos \frac{11 \varphi }{3}}{21446656}-\frac{414027 \sqrt{3} \cos \frac{17 \varphi }{3}}{493273088}+\frac{4179 \sqrt{3} \cos \frac{23 \varphi }{3}}{114509824}+\frac{375033 \sqrt{3} \cos \frac{29 \varphi }{3}}{4226453504}+\frac{135 \sin \frac{\varphi }{3}}{90112}
\nonumber\\
&-&\frac{675 \sin \frac{5 \varphi }{3}}{139264}-\frac{123489 \sin \frac{11 \varphi }{3}}{21446656}-\frac{1242081 \sin \frac{17 \varphi }{3}}{493273088}+\frac{12537 \sin \frac{23 \varphi }{3}}{114509824}+\frac{1125099 \sin \frac{29 \varphi }{3}}{4226453504}
\nonumber\\
\phi _{4,4,4}& = &\frac{135 \sqrt{3} \cos \frac{2 \varphi }{3}}{95744}+\frac{45 \sqrt{3} \cos \frac{4 \varphi }{3}}{131072}+\frac{143487 \sqrt{3} \cos \frac{8 \varphi }{3}}{61276160}+\frac{926991 \sqrt{3} \cos \frac{14 \varphi }{3}}{493273088}+\frac{32002167 \sqrt{3} \cos \frac{20 \varphi }{3}}{51300401152}-\frac{51351 \sqrt{3} \cos \frac{26 \varphi }{3}}{1300447232}+\frac{405 \sin \frac{2 \varphi }{3}}{95744}
\nonumber\\
&-&\frac{135 \sin \frac{4 \varphi }{3}}{131072}+\frac{430461 \sin \frac{8 \varphi }{3}}{61276160}+\frac{2780973 \sin \frac{14 \varphi }{3}}{493273088}+\frac{96006501 \sin \frac{20 \varphi }{3}}{51300401152}-\frac{154053 \sin \frac{26 \varphi }{3}}{1300447232}
\nonumber\\
\phi _{4,4,5}& = &-\frac{9855 \sqrt{3} \cos \frac{\varphi }{3}}{8912896}-\frac{572103 \sqrt{3} \cos \frac{5 \varphi }{3}}{245104640}-\frac{117 \sqrt{3} \cos \frac{7 \varphi }{3}}{524288}-\frac{2968203 \sqrt{3} \cos \frac{11 \varphi }{3}}{1127481344}-\frac{329537211 \sqrt{3} \cos \frac{17 \varphi }{3}}{205201604608}-\frac{2526434613 \sqrt{3} \cos \frac{23 \varphi }{3}}{5950846533632}
\nonumber\\
&+&\frac{3132411 \sqrt{3} \cos \frac{29 \varphi }{3}}{166457245696}+\frac{81227253 \sqrt{3} \cos \frac{35 \varphi }{3}}{3769767034880}+\frac{29565 \sin \frac{\varphi }{3}}{8912896}-\frac{1716309 \sin \frac{5 \varphi }{3}}{245104640}+\frac{351 \sin \frac{7 \varphi }{3}}{524288}-\frac{8904609 \sin \frac{11 \varphi }{3}}{1127481344}-\frac{988611633 \sin \frac{17 \varphi }{3}}{205201604608}
\nonumber\\
&-&\frac{7579303839 \sin \frac{23 \varphi }{3}}{5950846533632}+\frac{9397233 \sin \frac{29 \varphi }{3}}{166457245696}+\frac{243681759 \sin \frac{35 \varphi }{3}}{3769767034880}
\nonumber
\end{eqnarray}

\begin{eqnarray}
\phi _{6,4,0}& = &-\frac{3 \sqrt{3} \cos \frac{8 \varphi }{3}}{335104}-\frac{9 \sin \frac{8 \varphi }{3}}{335104}
\nonumber\\
\phi _{6,4,1}& = &\frac{3 \sqrt{3} \cos \frac{5 \varphi }{3}}{191488}+\frac{39 \sqrt{3} \cos \frac{11 \varphi }{3}}{1914880}-\frac{1791 \sqrt{3} \cos \frac{23 \varphi }{3}}{200392192}+\frac{9 \sin \frac{5 \varphi }{3}}{191488}+\frac{117 \sin \frac{11 \varphi }{3}}{1914880}-\frac{5373 \sin \frac{23 \varphi }{3}}{200392192}
\nonumber\\
\phi _{6,4,2}& = &-\frac{27 \sqrt{3} \cos \frac{2 \varphi }{3}}{1531904}-\frac{639 \sqrt{3} \cos \frac{8 \varphi }{3}}{15319040}-\frac{68463 \sqrt{3} \cos \frac{14 \varphi }{3}}{2466365440}+\frac{5373 \sqrt{3} \cos \frac{20 \varphi }{3}}{801568768}+\frac{191637 \sqrt{3} \cos \frac{26 \varphi }{3}}{23245494272}-\frac{81 \sin \frac{2 \varphi }{3}}{1531904}-\frac{1917 \sin \frac{8 \varphi }{3}}{15319040}
\nonumber\\
&-&\frac{205389 \sin \frac{14 \varphi }{3}}{2466365440}+\frac{16119 \sin \frac{20 \varphi }{3}}{801568768}+\frac{574911 \sin \frac{26 \varphi }{3}}{23245494272}
\nonumber\\
\phi _{6,4,3}& = &\frac{9 \sqrt{3} \cos \frac{\varphi }{3}}{557056}+\frac{1647 \sqrt{3} \cos \frac{5 \varphi }{3}}{30638080}+\frac{11493 \sqrt{3} \cos \frac{11 \varphi }{3}}{176168960}+\frac{1622793 \sqrt{3} \cos \frac{17 \varphi }{3}}{64125501440}-\frac{26865 \sqrt{3} \cos \frac{23 \varphi }{3}}{2905686784}-\frac{46559115 \sqrt{3} \cos \frac{29 \varphi }{3}}{5950846533632}-\frac{27 \sin \frac{\varphi }{3}}{557056}
\nonumber\\
&+&\frac{4941 \sin \frac{5 \varphi }{3}}{30638080}+\frac{34479 \sin \frac{11 \varphi }{3}}{176168960}+\frac{4868379 \sin \frac{17 \varphi }{3}}{64125501440}-\frac{80595 \sin \frac{23 \varphi }{3}}{2905686784}-\frac{139677345 \sin \frac{29 \varphi }{3}}{5950846533632}
\nonumber
\end{eqnarray}

\begin{eqnarray}
\phi _{8,4,0}& = &\frac{9 \sqrt{3} \cos \frac{8 \varphi }{3}}{107233280}+\frac{27 \sin \frac{8 \varphi }{3}}{107233280}
\nonumber\\
\phi _{8,4,1}& = &-\frac{81 \sqrt{3} \cos \frac{5 \varphi }{3}}{428933120}-\frac{477 \sqrt{3} \cos \frac{11 \varphi }{3}}{1973092352}+\frac{5373 \sqrt{3} \cos \frac{23 \varphi }{3}}{46490988544}-\frac{243 \sin \frac{5 \varphi }{3}}{428933120}-\frac{1431 \sin \frac{11 \varphi }{3}}{1973092352}+\frac{16119 \sin \frac{23 \varphi }{3}}{46490988544}
\nonumber
\end{eqnarray}


\begin{eqnarray}
\phi _{0,5,0}& = &\sin \frac{10 \varphi }{3}-\frac{1}{\sqrt{3}} \cos \frac{10 \varphi }{3}
\nonumber\\
\phi _{0,5,1}& = &\frac{\cos \frac{7 \varphi }{3}}{4 \sqrt{3}}+\frac{7 \cos \frac{13 \varphi }{3}}{52 \sqrt{3}}-\frac{1}{4} \sin \frac{7 \varphi }{3}-\frac{7}{52} \sin \frac{13 \varphi }{3}
\nonumber\\
\phi _{0,5,2}& = &-\frac{1}{32} \sqrt{3} \cos \frac{4 \varphi }{3}-\frac{5}{104} \sqrt{3} \cos \frac{10 \varphi }{3}+\frac{3}{32} \sin \frac{4 \varphi }{3}+\frac{15}{104} \sin \frac{10 \varphi }{3}
\nonumber\\
\phi _{0,5,3}& = &\frac{5 \cos \frac{\varphi }{3}}{128 \sqrt{3}}+\frac{53 \sqrt{3} \cos \frac{7 \varphi }{3}}{1664}+\frac{145 \sqrt{3} \cos \frac{13 \varphi }{3}}{6656}+\frac{77 \cos \frac{19 \varphi }{3}}{6656 \sqrt{3}}-\frac{5}{128} \sin \frac{\varphi }{3}-\frac{159 \sin \frac{7 \varphi }{3}}{1664}
\nonumber\\ & - &
\frac{435 \sin \frac{13 \varphi }{3}}{6656}-\frac{77 \sin \frac{19 \varphi }{3}}{6656}
\nonumber\\
\phi _{0,5,4}& = &-\frac{35 \cos \frac{2 \varphi }{3}}{2048 \sqrt{3}}-\frac{95 \cos \frac{4 \varphi }{3}}{1664 \sqrt{3}}-\frac{2407 \sqrt{3} \cos \frac{10 \varphi }{3}}{106496}-\frac{1043 \cos \frac{16 \varphi }{3}}{31616 \sqrt{3}}-\frac{35 \sin \frac{2 \varphi }{3}}{2048}+\frac{95 \sin \frac{4 \varphi }{3}}{1664}
\nonumber\\ & + &
\frac{7221 \sin \frac{10 \varphi }{3}}{106496}+\frac{1043 \sin \frac{16 \varphi }{3}}{31616}
\nonumber\\
\phi _{0,5,5}& = &\frac{1155 \sqrt{3} \cos \frac{\varphi }{3}}{106496}+\frac{21 \sqrt{3} \cos \frac{5 \varphi }{3}}{8192}+\frac{7685 \sqrt{3} \cos \frac{7 \varphi }{3}}{425984}+\frac{117621 \sqrt{3} \cos \frac{13 \varphi }{3}}{8093696}+\frac{42763 \cos \frac{19 \varphi }{3}}{2782208 \sqrt{3}}
\nonumber\\ & + &
\frac{14749 \cos \frac{25 \varphi }{3}}{11128832 \sqrt{3}}-\frac{3465 \sin \frac{\varphi }{3}}{106496}+\frac{63 \sin \frac{5 \varphi }{3}}{8192}-\frac{23055 \sin \frac{7 \varphi }{3}}{425984}-\frac{352863 \sin \frac{13 \varphi }{3}}{8093696}
\nonumber\\ & - &
\frac{42763 \sin \frac{19 \varphi }{3}}{2782208}-\frac{14749 \sin \frac{25 \varphi }{3}}{11128832}
\nonumber\\
\phi _{0,5,6}& = &-\frac{1281 \sqrt{3} \cos \frac{2 \varphi }{3}}{212992}-\frac{43645 \sqrt{3} \cos \frac{4 \varphi }{3}}{3407872}-\frac{77 \sqrt{3} \cos \frac{8 \varphi }{3}}{65536}-\frac{17885 \sqrt{3} \cos \frac{10 \varphi }{3}}{1245184}-\frac{562807 \sqrt{3} \cos \frac{16 \varphi }{3}}{64749568}
\nonumber\\ & - &
\frac{254793 \sqrt{3} \cos \frac{22 \varphi }{3}}{101171200}-\frac{3843 \sin \frac{2 \varphi }{3}}{212992}+\frac{130935 \sin \frac{4 \varphi }{3}}{3407872}-\frac{231 \sin \frac{8 \varphi }{3}}{65536}+\frac{53655 \sin \frac{10 \varphi }{3}}{1245184}
\nonumber\\ & + &
\frac{1688421 \sin \frac{16 \varphi }{3}}{64749568}+\frac{764379 \sin \frac{22 \varphi }{3}}{101171200}
\nonumber\\
\phi _{0,5,7}& = &\frac{57855 \sqrt{3} \cos \frac{\varphi }{3}}{6815744}+\frac{11165 \sqrt{3} \cos \frac{5 \varphi }{3}}{3407872}+\frac{1577555 \sqrt{3} \cos \frac{7 \varphi }{3}}{129499136}+\frac{143 \sqrt{3} \cos \frac{11 \varphi }{3}}{262144}
\nonumber\\ & + &
\frac{88150615 \cos \frac{13 \varphi }{3}}{2848980992 \sqrt{3}}+\frac{32566961 \sqrt{3} \cos \frac{19 \varphi }{3}}{6474956800}+\frac{1929147 \sqrt{3} \cos \frac{25 \varphi }{3}}{1618739200}+\frac{3279997 \cos \frac{31 \varphi }{3}}{17806131200 \sqrt{3}}
\nonumber\\ & - &
\frac{173565 \sin \frac{\varphi }{3}}{6815744}+\frac{33495 \sin \frac{5 \varphi }{3}}{3407872}-\frac{4732665 \sin \frac{7 \varphi }{3}}{129499136}+\frac{429 \sin \frac{11 \varphi }{3}}{262144}-\frac{88150615 \sin \frac{13 \varphi }{3}}{2848980992}
\nonumber\\ & - &
\frac{97700883 \sin \frac{19 \varphi }{3}}{6474956800}-\frac{5787441 \sin \frac{25 \varphi }{3}}{1618739200}-\frac{3279997 \sin \frac{31 \varphi }{3}}{17806131200}
\nonumber\\
\phi _{0,5,8}& = &-\frac{1172325 \sqrt{3} \cos \frac{2 \varphi }{3}}{218103808}-\frac{1216131 \sqrt{3} \cos \frac{4 \varphi }{3}}{129499136}-\frac{231 \sqrt{3} \cos \frac{8 \varphi }{3}}{131072}-\frac{107170105 \cos \frac{10 \varphi }{3}}{3506438144 \sqrt{3}}
\nonumber\\ & - &
\frac{2145 \sqrt{3} \cos \frac{14 \varphi }{3}}{8388608}-\frac{2318071 \cos \frac{16 \varphi }{3}}{111288320 \sqrt{3}}-\frac{292877947 \sqrt{3} \cos \frac{22 \varphi }{3}}{103599308800}-\frac{964432279 \cos \frac{28 \varphi }{3}}{551990067200 \sqrt{3}}
\nonumber\\ & - &
\frac{3516975 \sin \frac{2 \varphi }{3}}{218103808}+\frac{3648393 \sin \frac{4 \varphi }{3}}{129499136}-\frac{693 \sin \frac{8 \varphi }{3}}{131072}+\frac{107170105 \sin \frac{10 \varphi }{3}}{3506438144}-\frac{6435 \sin \frac{14 \varphi }{3}}{8388608}
\nonumber\\ & + &
\frac{2318071 \sin \frac{16 \varphi }{3}}{111288320}+\frac{878633841 \sin \frac{22 \varphi }{3}}{103599308800}+\frac{964432279 \sin \frac{28 \varphi }{3}}{551990067200}
\nonumber
\end{eqnarray}

\begin{eqnarray}
\phi _{0,6,0}& = &\cos  4\varphi
\nonumber\\
\phi _{0,6,1}& = &-\frac{1}{4} \cos  3 \varphi -\frac{3}{20} \cos  5 \varphi
\nonumber\\
\phi _{0,6,2}& = &\frac{3}{32} \cos  2 \varphi +\frac{3}{20} \cos  4 \varphi
\nonumber\\
\phi _{0,6,3}& = &-\frac{5 \cos  \varphi }{128}-\frac{63}{640} \cos  3 \varphi -\frac{11}{160} \cos  5 \varphi -\frac{1}{80} \cos  7 \varphi
\nonumber\\
\phi _{0,6,4}& = &\frac{35}{2048}+\frac{15}{256} \cos  2 \varphi +\frac{363 \cos  4 \varphi }{5120}+\frac{157 \cos  6 \varphi }{4480}
\nonumber\\
\phi _{0,6,5}& = &-\frac{21 \cos  \varphi }{512}-\frac{231 \cos  3\varphi }{4096}-\frac{6603 \cos  5 \varphi }{143360}-\frac{471 \cos  7 \varphi }{28672}-\frac{41 \cos  9 \varphi }{28672}
\nonumber\\
\phi _{0,6,6}& = &\frac{189}{10240}+\frac{2849 \cos  2 \varphi }{65536}+\frac{1301 \cos  4 \varphi }{28672}+\frac{25461 \cos  6 \varphi }{917504}+\frac{349 \cos  8 \varphi }{43008}
\nonumber\\
\phi _{0,6,7}& = &-\frac{47817 \cos  \varphi }{1310720}-\frac{10495 \cos  3 \varphi }{262144}-\frac{120495 \cos  5 \varphi }{3670016}-\frac{59261 \cos  7 \varphi }{3670016}
\nonumber\\ & - &
\frac{6631 \cos  9 \varphi }{1720320}-\frac{229 \cos  11 \varphi }{1146880}
\nonumber
\end{eqnarray}

\begin{eqnarray}
\phi _{0,7,0}& = & \sin \frac{14 \varphi }{3}+\frac{1}{\sqrt{3}} \cos \frac{14 \varphi }{3}
\nonumber\\
\phi _{0,7,1}& = & -\frac{\cos \frac{11 \varphi }{3}}{4 \sqrt{3}}-\frac{11 \cos \frac{17 \varphi }{3}}{68 \sqrt{3}}-\frac{1}{4} \sin \frac{11 \varphi }{3}-\frac{11}{68} \sin \frac{17 \varphi }{3}
\nonumber\\
\phi _{0,7,2}& = & \frac{1}{32} \sqrt{3} \cos \frac{8 \varphi }{3}+\frac{7}{136} \sqrt{3} \cos \frac{14 \varphi }{3}+\frac{3}{32} \sin \frac{8 \varphi }{3}+\frac{21}{136} \sin \frac{14 \varphi }{3}
\nonumber\\
\phi _{0,7,3}& = & -\frac{5 \cos \frac{5 \varphi }{3}}{128 \sqrt{3}}-\frac{73 \sqrt{3} \cos \frac{11 \varphi }{3}}{2176}-\frac{259 \sqrt{3} \cos \frac{17 \varphi }{3}}{10880}-\frac{143 \cos \frac{23 \varphi }{3}}{10880 \sqrt{3}}
\nonumber\\ & - &
\frac{5}{128} \sin \frac{5 \varphi }{3}-\frac{219 \sin \frac{11 \varphi }{3}}{2176}-\frac{777 \sin \frac{17 \varphi }{3}}{10880}-\frac{143 \sin \frac{23 \varphi }{3}}{10880}
\nonumber\\
\phi _{0,7,4}& = & \frac{35 \cos \frac{2 \varphi }{3}}{2048 \sqrt{3}}+\frac{65 \cos \frac{8 \varphi }{3}}{1088 \sqrt{3}}+\frac{851 \sqrt{3} \cos \frac{14 \varphi }{3}}{34816}+\frac{1835 \cos \frac{20 \varphi }{3}}{50048 \sqrt{3}}+\frac{35 \sin \frac{2 \varphi }{3}}{2048}
\nonumber\\ & + &
\frac{65 \sin \frac{8 \varphi }{3}}{1088}+\frac{2553 \sin \frac{14 \varphi }{3}}{34816}+\frac{1835 \sin \frac{20 \varphi }{3}}{50048}
\nonumber\\
\phi _{0,7,5}& = & -\frac{21 \sqrt{3} \cos \frac{\varphi }{3}}{8192}-\frac{1575 \sqrt{3} \cos \frac{5 \varphi }{3}}{139264}-\frac{2701 \sqrt{3} \cos \frac{11 \varphi }{3}}{139264}-\frac{51273 \sqrt{3} \cos \frac{17 \varphi }{3}}{3203072}
\nonumber\\ & - &
\frac{89915 \cos \frac{23 \varphi }{3}}{5204992 \sqrt{3}}-\frac{979 \cos \frac{29 \varphi }{3}}{650624 \sqrt{3}}+\frac{63 \sin \frac{\varphi }{3}}{8192}-\frac{4725 \sin \frac{5 \varphi }{3}}{139264}-\frac{8103 \sin \frac{11 \varphi }{3}}{139264}
\nonumber\\ & - &
\frac{153819 \sin \frac{17 \varphi }{3}}{3203072}-\frac{89915 \sin \frac{23 \varphi }{3}}{5204992}-\frac{979 \sin \frac{29 \varphi }{3}}{650624}
\nonumber\\
\phi _{0,7,6}& = & \frac{1743 \sqrt{3} \cos \frac{2 \varphi }{3}}{278528}+\frac{77 \sqrt{3} \cos \frac{4 \varphi }{3}}{65536}+\frac{15281 \sqrt{3} \cos \frac{8 \varphi }{3}}{1114112}+\frac{100799 \sqrt{3} \cos \frac{14 \varphi }{3}}{6406144}
\nonumber\\ & + &
\frac{3231697 \sqrt{3} \cos \frac{20 \varphi }{3}}{333119488}+\frac{132627 \sqrt{3} \cos \frac{26 \varphi }{3}}{46444544}+\frac{5229 \sin \frac{2 \varphi }{3}}{278528}-\frac{231 \sin \frac{4 \varphi }{3}}{65536}
\nonumber\\ & + &
\frac{45843 \sin \frac{8 \varphi }{3}}{1114112}+\frac{302397 \sin \frac{14 \varphi }{3}}{6406144}+\frac{9695091 \sin \frac{20 \varphi }{3}}{333119488}+\frac{397881 \sin \frac{26 \varphi }{3}}{46444544}
\nonumber
\end{eqnarray}

\begin{eqnarray}
\phi _{0,8,0}& = & \sin \frac{16 \varphi }{3}-\frac{1}{\sqrt{3}} \cos \frac{16 \varphi }{3}
\nonumber\\
\phi _{0,8,1}& = & \frac{\cos \frac{13 \varphi }{3}}{4 \sqrt{3}}+\frac{13 \cos \frac{19 \varphi }{3}}{76 \sqrt{3}}-\frac{1}{4} \sin \frac{13 \varphi }{3}-\frac{13}{76} \sin \frac{19 \varphi }{3}
\nonumber\\
\phi _{0,8,2}& = & -\frac{1}{32} \sqrt{3} \cos \frac{10 \varphi }{3}-\frac{1}{19} \sqrt{3} \cos \frac{16 \varphi }{3}+\frac{3}{32} \sin \frac{10 \varphi }{3}+\frac{3}{19} \sin \frac{16 \varphi }{3}
\nonumber\\
\phi _{0,8,3}& = & \frac{5 \cos \frac{7 \varphi }{3}}{128 \sqrt{3}}+\frac{83 \sqrt{3} \cos \frac{13 \varphi }{3}}{2432}+\frac{41 \sqrt{3} \cos \frac{19 \varphi }{3}}{1672}+\frac{91 \cos \frac{25 \varphi }{3}}{6688 \sqrt{3}}-\frac{5}{128} \sin \frac{7 \varphi }{3}
\nonumber\\ & - &
\frac{249 \sin \frac{13 \varphi }{3}}{2432}-\frac{123 \sin \frac{19 \varphi }{3}}{1672}-\frac{91 \sin \frac{25 \varphi }{3}}{6688}
\nonumber\\
\phi _{0,8,4}& = & -\frac{35 \cos \frac{4 \varphi }{3}}{2048 \sqrt{3}}-\frac{295 \cos \frac{10 \varphi }{3}}{4864 \sqrt{3}}-\frac{5371 \sqrt{3} \cos \frac{16 \varphi }{3}}{214016}-\frac{2309 \cos \frac{22 \varphi }{3}}{60800 \sqrt{3}}+\frac{35 \sin \frac{4 \varphi }{3}}{2048}
\nonumber\\ & + &
\frac{295 \sin \frac{10 \varphi }{3}}{4864}+\frac{16113 \sin \frac{16 \varphi }{3}}{214016}+\frac{2309 \sin \frac{22 \varphi }{3}}{60800}
\nonumber\\
\phi _{0,8,5}& = & \frac{21 \sqrt{3} \cos \frac{\varphi }{3}}{8192}+\frac{1785 \sqrt{3} \cos \frac{7 \varphi }{3}}{155648}+\frac{17015 \sqrt{3} \cos \frac{13 \varphi }{3}}{856064}+\frac{354033 \sqrt{3} \cos \frac{19 \varphi }{3}}{21401600}
\nonumber\\ & + &
\frac{122377 \cos \frac{25 \varphi }{3}}{6809600 \sqrt{3}}+\frac{116909 \cos \frac{31 \varphi }{3}}{74905600 \sqrt{3}}-\frac{63 \sin \frac{\varphi }{3}}{8192}-\frac{5355 \sin \frac{7 \varphi }{3}}{155648}-\frac{51045 \sin \frac{13 \varphi }{3}}{856064}
\nonumber\\ & - &
\frac{1062099 \sin \frac{19 \varphi }{3}}{21401600}-\frac{122377 \sin \frac{25 \varphi }{3}}{6809600}-\frac{116909 \sin \frac{31 \varphi }{3}}{74905600}
\nonumber\\
\phi _{0,8,6}& = & -\frac{77 \sqrt{3} \cos \frac{2 \varphi }{3}}{65536}-\frac{987 \sqrt{3} \cos \frac{4 \varphi }{3}}{155648}-\frac{96145 \sqrt{3} \cos \frac{10 \varphi }{3}}{6848512}-\frac{4343 \sqrt{3} \cos \frac{16 \varphi }{3}}{267520}
\nonumber\\ & - &
\frac{24148337 \sqrt{3} \cos \frac{22 \varphi }{3}}{2396979200}-\frac{988653 \sqrt{3} \cos \frac{28 \varphi }{3}}{331724800}-\frac{231 \sin \frac{2 \varphi }{3}}{65536}+\frac{2961 \sin \frac{4 \varphi }{3}}{155648}
\nonumber\\ & + &
\frac{288435 \sin \frac{10 \varphi }{3}}{6848512}+\frac{13029 \sin \frac{16 \varphi }{3}}{267520}+\frac{72445011 \sin \frac{22 \varphi }{3}}{2396979200}+\frac{2965959 \sin \frac{28 \varphi }{3}}{331724800}
\nonumber
\end{eqnarray}

\begin{eqnarray}
\phi _{0,9,0}& = & \cos  6 \varphi
\nonumber\\
\phi _{0,9,1}& = & -\frac{1}{4} \cos  5 \varphi -\frac{5}{28} \cos  7 \varphi
\nonumber\\
\phi _{0,9,2}& = & \frac{3}{32} \cos  4 \varphi +\frac{9}{56} \cos  6\varphi
\nonumber\\
\phi _{0,9,3}& = & -\frac{5}{128} \cos  3 \varphi -\frac{93}{896} \cos  5 \varphi -\frac{135 \cos  7 \varphi }{1792}-\frac{25 \cos  9 \varphi }{1792}
\nonumber\\
\phi _{0,9,4}& = & \frac{35 \cos  2 \varphi }{2048}+\frac{55}{896} \cos  4 \varphi +\frac{315 \cos  6 \varphi }{4096}+\frac{5}{128} \cos  8 \varphi \nonumber\\ \phi _{0,9,5}& = & -\frac{63 \cos  \varphi }{8192}-\frac{285 \cos  3 \varphi }{8192}-\frac{6975 \cos  5 \varphi }{114688}-\frac{835 \cos  7 \varphi }{16384}-\frac{19 \cos  9 \varphi }{1024}-\frac{23 \cos  11 \varphi }{14336}
\nonumber
\end{eqnarray}

\begin{eqnarray}
\phi _{0,10,0}& = & \sin \frac{20 \varphi }{3}+\frac{1}{\sqrt{3}} \cos \frac{20 \varphi }{3}
\nonumber\\
\phi _{0,10,1}& = & -\frac{\cos \frac{17 \varphi }{3}}{4 \sqrt{3}}-\frac{17 \cos \frac{23 \varphi }{3}}{92 \sqrt{3}}-\frac{1}{4} \sin \frac{17 \varphi }{3}-\frac{17}{92} \sin \frac{23 \varphi }{3}
\nonumber\\
\phi _{0,10,2}& = & \frac{1}{32} \sqrt{3} \cos \frac{14 \varphi }{3}+\frac{5}{92} \sqrt{3} \cos \frac{20 \varphi }{3}+\frac{3}{32} \sin \frac{14 \varphi }{3}+\frac{15}{92} \sin \frac{20 \varphi }{3}
\nonumber\\
\phi _{0,10,3}& = & -\frac{5 \cos \frac{11 \varphi }{3}}{128 \sqrt{3}}-\frac{103 \sqrt{3} \cos \frac{17 \varphi }{3}}{2944}-\frac{245 \sqrt{3} \cos \frac{23 \varphi }{3}}{9568}-\frac{17 \cos \frac{29 \varphi }{3}}{1196 \sqrt{3}}-\frac{5}{128} \sin \frac{11 \varphi }{3}
\nonumber\\ & - &
\frac{309 \sin \frac{17 \varphi }{3}}{2944}-\frac{735 \sin \frac{23 \varphi }{3}}{9568}-\frac{17 \sin \frac{29 \varphi }{3}}{1196}
\nonumber\\
\phi _{0,10,4}& = & \frac{35 \cos \frac{8 \varphi }{3}}{2048 \sqrt{3}}+\frac{365 \cos \frac{14 \varphi }{3}}{5888 \sqrt{3}}+\frac{7987 \sqrt{3} \cos \frac{20 \varphi }{3}}{306176}+\frac{3413 \cos \frac{26 \varphi }{3}}{85376 \sqrt{3}}+\frac{35 \sin \frac{8 \varphi }{3}}{2048}
\nonumber\\ & + &
\frac{365 \sin \frac{14 \varphi }{3}}{5888}+\frac{23961 \sin \frac{20 \varphi }{3}}{306176}+\frac{3413 \sin \frac{26 \varphi }{3}}{85376}
\nonumber
\end{eqnarray}

\begin{eqnarray}
\phi _{0,11,0}& = & \sin \frac{22 \varphi }{3}-\frac{1}{\sqrt{3}} \cos \frac{22 \varphi }{3}
\nonumber\\
\phi _{0,11,1}& = & \frac{\cos \frac{19 \varphi }{3}}{4 \sqrt{3}}+\frac{19 \cos \frac{25 \varphi }{3}}{100 \sqrt{3}}-\frac{1}{4} \sin \frac{19 \varphi }{3}-\frac{19}{100} \sin \frac{25 \varphi }{3}
\nonumber\\
\phi _{0,11,2}& = & -\frac{1}{32} \sqrt{3} \cos \frac{16 \varphi }{3}-\frac{11}{200} \sqrt{3} \cos \frac{22 \varphi }{3}+\frac{3}{32} \sin \frac{16 \varphi }{3}+\frac{33}{200} \sin \frac{22 \varphi }{3}
\nonumber\\
\phi _{0,11,3}& = & \frac{5 \cos \frac{13 \varphi }{3}}{128 \sqrt{3}}+\frac{113 \sqrt{3} \cos \frac{19 \varphi }{3}}{3200}+\frac{583 \sqrt{3} \cos \frac{25 \varphi }{3}}{22400}+\frac{323 \cos \frac{31 \varphi }{3}}{22400 \sqrt{3}}-\frac{5}{128} \sin \frac{13 \varphi }{3}
\nonumber\\ & - &
\frac{339 \sin \frac{19 \varphi }{3}}{3200}-\frac{1749 \sin \frac{25 \varphi }{3}}{22400}-\frac{323 \sin \frac{31 \varphi }{3}}{22400}
\nonumber\\
\phi _{0,11,4}& = & -\frac{35 \cos \frac{10 \varphi }{3}}{2048 \sqrt{3}}-\frac{\cos \frac{16 \varphi }{3}}{16 \sqrt{3}}-\frac{9487 \sqrt{3} \cos \frac{22 \varphi }{3}}{358400}-\frac{4043 \cos \frac{28 \varphi }{3}}{99200 \sqrt{3}}+\frac{35 \sin \frac{10 \varphi }{3}}{2048}
\nonumber\\ & + &
\frac{1}{16} \sin \frac{16 \varphi }{3}+\frac{28461 \sin \frac{22 \varphi }{3}}{358400}+\frac{4043 \sin \frac{28 \varphi }{3}}{99200}
\nonumber
\end{eqnarray}

\begin{eqnarray}
\phi _{0,12,0}& = & \cos  8 \varphi
\nonumber\\
\phi _{0,12,1}& = & -\frac{1}{4} \cos  7 \varphi -\frac{7}{36} \cos  9 \varphi
\nonumber\\
\phi _{0,12,2}& = & \frac{3}{32} \cos  6 \varphi +\frac{1}{6} \cos  8 \varphi
\nonumber\\
\phi _{0,12,3}& = & -\frac{5}{128} \cos  5 \varphi -\frac{41}{384} \cos  7 \varphi -\frac{19}{240} \cos  9 \varphi -\frac{7}{480} \cos  11 \varphi
\nonumber
\end{eqnarray}

\begin{eqnarray}
\phi _{0,13,0}& = & \sin \frac{26\varphi }{3}+\frac{1}{\sqrt{3}} \cos \frac{26\varphi }{3}
\nonumber\\
\phi _{0,13,1}& = & -\frac{\cos \frac{23 \varphi }{3}}{4 \sqrt{3}}-\frac{23 \cos \frac{29 \varphi }{3}}{116 \sqrt{3}}-\frac{1}{4} \sin \frac{23 \varphi }{3}-\frac{23}{116} \sin \frac{29 \varphi }{3}
\nonumber\\
\phi _{0,13,2}& = & \frac{1}{32} \sqrt{3} \cos \frac{20 \varphi }{3}+\frac{13}{232} \sqrt{3} \cos \frac{26 \varphi }{3}+\frac{3}{32} \sin \frac{20 \varphi }{3}+\frac{39}{232} \sin \frac{26 \varphi }{3}
\nonumber
\end{eqnarray}

\begin{eqnarray}
\phi _{0,14,0}& = & \sin \frac{28 \varphi }{3}-\frac{1}{\sqrt{3}} \cos \frac{28 \varphi }{3}
\nonumber\\
\phi _{0,14,1}& = & \frac{\cos \frac{25 \varphi }{3}}{4 \sqrt{3}}+\frac{25 \cos \frac{31 \varphi }{3}}{124 \sqrt{3}}-\frac{1}{4} \sin \frac{25 \varphi }{3}-\frac{25}{124} \sin \frac{31 \varphi }{3}
\nonumber\\
\phi _{0,14,2}& = & -\frac{1}{32} \sqrt{3} \cos \frac{22 \varphi }{3}-\frac{7}{124} \sqrt{3} \cos \frac{28 \varphi }{3}+\frac{3}{32} \sin \frac{22 \varphi }{3}+\frac{21}{124} \sin \frac{28 \varphi }{3}
\nonumber
\end{eqnarray}

\begin{eqnarray}
\phi_{0,15,0} & = &  \cos  10 \varphi
\nonumber\\
\phi _{0,15,1}& = & -\frac{1}{4} \cos  9 \varphi -\frac{9}{44} \cos  11 \varphi
\nonumber
\end{eqnarray}

\begin{eqnarray}
\phi _{0,16,0}& = & \sin \frac{32 \varphi }{3}+\frac{1}{\sqrt{3}} \cos \frac{32 \varphi }{3}
\nonumber
\end{eqnarray}

\begin{eqnarray}
\phi _{0,17,0}& = & \sin \frac{34\varphi }{3}-\frac{1}{\sqrt{3}} \cos \frac{34\varphi }{3}
\nonumber
\end{eqnarray}
\normalsize


\section{Dual functions and dual shadows for $90^\circ$ V-notch $-\pi \le \varphi \le \pi/2$}
\subsection{First singular exponent ($j=1$)}\mbox{}

\scriptsize
\begin{eqnarray}
\psi _{0,1,0}& = &\sin \frac{2\varphi }{3}+\frac{1}{\sqrt{3}} \cos \frac{2\varphi }{3}
\nonumber\\
\psi _{0,1,1}& = &-\frac{1}{4 \sqrt{3}}\cos \frac{5 \varphi }{3}-\frac{1}{4}\sin \frac{5 \varphi }{3}-\frac{5}{4} \sin \frac{\varphi }{3}+\frac{5}{4 \sqrt{3}} \cos \frac{\varphi }{3}
\nonumber\\
\psi _{0,1,2}& = &-\frac{\sqrt{3}}{8}\cos \frac{2 \varphi }{3}-\frac{3}{8} \sin \frac{2 \varphi }{3}+\frac{\sqrt{3}}{32}\cos \frac{8 \varphi }{3}+\frac{3}{32}\sin \frac{8 \varphi }{3}
\nonumber\\
\psi _{0,1,3}& = &\frac{5 \sqrt{3}}{128}\cos \frac{\varphi }{3}-\frac{15}{128}\sin \frac{\varphi }{3}+\frac{7 \sqrt{3}}{128}\cos \frac{5 \varphi }{3}+\frac{21}{128}\sin \frac{5 \varphi }{3}-\frac{5}{128 \sqrt{3}}\cos \frac{11 \varphi }{3}-\frac{5}{128}\sin \frac{11 \varphi }{3}
\nonumber\\
&+&\frac{25}{128 \sqrt{3}} \cos \frac{7\varphi }{3}-\frac{25}{128} \sin \frac{7\varphi }{3}
\nonumber\\
\psi _{0,1,4}& = &-\frac{65 \sqrt{3}}{2048}\cos \frac{2 \varphi }{3}-\frac{195}{2048}\sin \frac{2 \varphi }{3}-\frac{85}{896 \sqrt{3}}\cos \frac{4 \varphi }{3}+\frac{85}{896}\sin \frac{4 \varphi }{3}-\frac{5}{64 \sqrt{3}}\cos \frac{8 \varphi }{3}-\frac{5}{64}\sin \frac{8 \varphi }{3}+\frac{35}{2048 \sqrt{3}}\cos \frac{14 \varphi }{3}
\nonumber\\
&+&\frac{35}{2048}\sin \frac{14 \varphi }{3}
\nonumber
\end{eqnarray}

\begin{eqnarray}
\psi _{2,1,0}& = &-\frac{1}{4} \sqrt{3} \cos \frac{2 \varphi }{3}-\frac{3}{4} \sin \frac{2 \varphi }{3}
\nonumber\\
\psi _{2,1,1}& = &\frac{1}{32} \sqrt{3} \cos \frac{\varphi }{3}+\frac{3}{16} \sqrt{3} \cos \frac{5 \varphi }{3}+\frac{31}{224} \sqrt{3} \cos \frac{7 \varphi }{3}-\frac{3}{32} \sin \frac{\varphi }{3}+\frac{9}{16} \sin \frac{5 \varphi }{3}-\frac{93}{224} \sin \frac{7 \varphi }{3}
\nonumber\\
\psi _{2,1,2}& = &-\frac{15}{256} \sqrt{3} \cos \frac{2 \varphi }{3}-\frac{15}{448} \sqrt{3} \cos \frac{4 \varphi }{3}-\frac{15}{128} \sqrt{3} \cos \frac{8 \varphi }{3}-\frac{45}{256} \sin \frac{2 \varphi }{3}+\frac{45}{448} \sin \frac{4 \varphi }{3}-\frac{45}{128} \sin \frac{8 \varphi }{3}
\nonumber\\
\psi _{2,1,3}& = &\frac{3}{112} \sqrt{3} \cos \frac{\varphi }{3}+\frac{15}{256} \sqrt{3} \cos \frac{5 \varphi }{3}+\frac{167 \sqrt{3} \cos \frac{7 \varphi }{3}}{35840}+\frac{35}{512} \sqrt{3} \cos \frac{11 \varphi }{3}+\frac{16659 \sqrt{3} \cos \frac{13 \varphi }{3}}{465920}-\frac{9}{112} \sin \frac{\varphi }{3}+\frac{45}{256} \sin \frac{5 \varphi }{3}
\nonumber\\
&-&\frac{501 \sin \frac{7 \varphi }{3}}{35840}+\frac{105}{512} \sin \frac{11 \varphi }{3}-\frac{49977 \sin \frac{13 \varphi }{3}}{465920}
\nonumber\\
\psi _{2,1,4}& = &-\frac{241 \sqrt{3} \cos \frac{2 \varphi }{3}}{8192}-\frac{1857 \sqrt{3} \cos \frac{4 \varphi }{3}}{286720}-\frac{385 \sqrt{3} \cos \frac{8 \varphi }{3}}{8192}-\frac{341 \sqrt{3} \cos \frac{10 \varphi }{3}}{46592}-\frac{315 \sqrt{3} \cos \frac{14 \varphi }{3}}{8192}-\frac{723 \sin \frac{2 \varphi }{3}}{8192}+\frac{5571 \sin \frac{4 \varphi }{3}}{286720}
\nonumber\\
&-&\frac{1155 \sin \frac{8 \varphi }{3}}{8192}+\frac{1023 \sin \frac{10 \varphi }{3}}{46592}-\frac{945 \sin \frac{14 \varphi }{3}}{8192}
\nonumber
\end{eqnarray}

\begin{eqnarray}
\psi _{4,1,0}& = &\frac{3}{128} \sqrt{3} \cos \frac{2 \varphi }{3}+\frac{9}{128} \sin \frac{2 \varphi }{3}
\nonumber\\
\psi _{4,1,1}& = &-\frac{9}{512} \sqrt{3} \cos \frac{\varphi }{3}-\frac{15}{512} \sqrt{3} \cos \frac{5 \varphi }{3}-\frac{93 \sqrt{3} \cos \frac{7 \varphi }{3}}{8960}+\frac{9 \sqrt{3} \cos \frac{13 \varphi }{3}}{1280}+\frac{27}{512} \sin \frac{\varphi }{3}-\frac{45}{512} \sin \frac{5 \varphi }{3}+\frac{279 \sin \frac{7 \varphi }{3}}{8960}
\nonumber\\
&-&\frac{27 \sin \frac{13 \varphi }{3}}{1280}
\nonumber\\
\psi _{4,1,2}& = &\frac{15}{512} \sqrt{3} \cos \frac{2 \varphi }{3}+\frac{2601 \sqrt{3} \cos \frac{4 \varphi }{3}}{143360}+\frac{105 \sqrt{3} \cos \frac{8 \varphi }{3}}{4096}+\frac{159 \sqrt{3} \cos \frac{10 \varphi }{3}}{23296}+\frac{45}{512} \sin \frac{2 \varphi }{3}-\frac{7803 \sin \frac{4 \varphi }{3}}{143360}+\frac{315 \sin \frac{8 \varphi }{3}}{4096}
\nonumber\\
&-&\frac{477 \sin \frac{10 \varphi }{3}}{23296}
\nonumber\\
\psi _{4,1,3}& = &-\frac{3051 \sqrt{3} \cos \frac{\varphi }{3}}{114688}-\frac{525 \sqrt{3} \cos \frac{5 \varphi }{3}}{16384}-\frac{3051 \sqrt{3} \cos \frac{7 \varphi }{3}}{212992}-\frac{315 \sqrt{3} \cos \frac{11 \varphi }{3}}{16384}-\frac{190557 \sqrt{3} \cos \frac{13 \varphi }{3}}{29818880}+\frac{1248279 \sqrt{3} \cos \frac{19 \varphi }{3}}{566558720}
\nonumber\\
&+&\frac{9153 \sin \frac{\varphi }{3}}{114688}-\frac{1575 \sin \frac{5 \varphi }{3}}{16384}+\frac{9153 \sin \frac{7 \varphi }{3}}{212992}-\frac{945 \sin \frac{11 \varphi }{3}}{16384}+\frac{571671 \sin \frac{13 \varphi }{3}}{29818880}-\frac{3744837 \sin \frac{19 \varphi }{3}}{566558720}
\nonumber\\
\psi _{4,1,4}& = &\frac{10317 \sqrt{3} \cos \frac{2 \varphi }{3}}{327680}+\frac{9219 \sqrt{3} \cos \frac{4 \varphi }{3}}{425984}+\frac{945 \sqrt{3} \cos \frac{8 \varphi }{3}}{32768}+\frac{418479 \sqrt{3} \cos \frac{10 \varphi }{3}}{36700160}+\frac{3465 \sqrt{3} \cos \frac{14 \varphi }{3}}{262144}+\frac{3825 \sqrt{3} \cos \frac{16 \varphi }{3}}{1011712}
\nonumber\\
&+&\frac{30951 \sin \frac{2 \varphi }{3}}{327680}-\frac{27657 \sin \frac{4 \varphi }{3}}{425984}+\frac{2835 \sin \frac{8 \varphi }{3}}{32768}-\frac{1255437 \sin \frac{10 \varphi }{3}}{36700160}+\frac{10395 \sin \frac{14 \varphi }{3}}{262144}-\frac{11475 \sin \frac{16 \varphi }{3}}{1011712}
\nonumber
\end{eqnarray}
\normalsize
\newpage

\subsection{Higher order singular exponent ($j=2,4$)}\mbox{}

\scriptsize
\begin{eqnarray}
\psi _{0,2,0}& = &\sin \frac{4\varphi }{3}-\frac{1}{\sqrt{3}} \cos \frac{4\varphi }{3}
\nonumber\\
\psi _{0,2,1}& = &\frac{1}{4 \sqrt{3}}\cos \frac{7\varphi }{3}-\frac{1}{4}\sin \frac{7\varphi }{3}-\frac{7}{4}\sin \frac{\varphi }{3}+\frac{7}{4 \sqrt{3}} \cos \frac{\varphi }{3}
\nonumber\\
\psi _{0,2,2}& = &-\frac{\sqrt{3}}{4}\cos \frac{4 \varphi }{3}+\frac{3}{4}\sin \frac{4 \varphi }{3}-\frac{\sqrt{3}}{32}\cos \frac{10 \varphi }{3}+\frac{3}{32}\sin \frac{10 \varphi }{3}
\nonumber\\
\psi _{0,2,3}& = &\frac{\sqrt{3}}{32}\cos \frac{\varphi }{3}-\frac{3}{32}\sin \frac{\varphi }{3}+\frac{17 \sqrt{3}}{128}\cos \frac{7 \varphi }{3}-\frac{51}{128}\sin \frac{7 \varphi }{3}+\frac{5}{128 \sqrt{3}}\cos \frac{13 \varphi }{3}-\frac{5}{128}\sin \frac{13 \varphi }{3}+\frac{7}{16} \sin \frac{5\varphi }{3}+\frac{7}{16 \sqrt{3}} \cos \frac{5\varphi }{3}
\nonumber\\
\psi _{0,2,4}& = &-\frac{91}{640 \sqrt{3}} \cos \frac{2 \varphi }{3}-\frac{29 \sqrt{3}}{1024} \cos \frac{4 \varphi }{3}-\frac{55}{256 \sqrt{3}} \cos \frac{10 \varphi }{3}-\frac{35}{2048 \sqrt{3}} \cos \frac{16 \varphi }{3}-\frac{91}{640}\sin \frac{2 \varphi }{3}+\frac{87}{1024} \sin \frac{4 \varphi }{3}+\frac{55}{256} \sin \frac{10 \varphi }{3}
\nonumber\\
&+&\frac{35}{2048} \sin \frac{16 \varphi }{3}
\nonumber
\end{eqnarray}

\begin{eqnarray}
\psi _{0,4,0}& = &\sin \frac{8}{3} \varphi +\frac{1}{\sqrt{3}} \cos \frac{8}{3} \varphi
\nonumber\\
\psi _{0,4,1}& = &-\frac{1}{4 \sqrt{3}}\cos \frac{11 \varphi }{3}-\frac{1}{4}\sin \frac{11 \varphi }{3}-\frac{11}{20} \sin \frac{5\varphi }{3}-\frac{11}{20 \sqrt{3}} \cos \frac{5\varphi }{3}
\nonumber\\
\psi _{0,4,2}& = & \frac{1}{10} \sqrt{3} \cos\frac{8 \varphi }{3}+\frac{1}{32} \sqrt{3} \cos\frac{14 \varphi }{3}+\frac{3}{10} \sin\frac{8 \varphi }{3}+\frac{3}{32} \sin\frac{14 \varphi }{3}
\nonumber\\
\psi _{0,4,3}& = &\frac{11 \cos\frac{\varphi }{3}}{160 \sqrt{3}}-\frac{7}{80} \sqrt{3} \cos\frac{5 \varphi }{3}-\frac{37}{640} \sqrt{3} \cos\frac{11 \varphi }{3}-\frac{5 \cos\frac{17 \varphi }{3}}{128 \sqrt{3}}-\frac{11}{160} \sin\frac{\varphi }{3}-\frac{21}{80} \sin\frac{5 \varphi }{3}-\frac{111}{640} \sin\frac{11 \varphi }{3}-\frac{5}{128} \sin\frac{17 \varphi }{3}
\nonumber\\
\psi _{0,4,4}& = &-\frac{53 \cos\frac{2 \varphi }{3}}{640 \sqrt{3}}+\frac{427 \sqrt{3} \cos\frac{8 \varphi }{3}}{5120}+\frac{25 \cos\frac{14 \varphi }{3}}{256 \sqrt{3}}+\frac{35 \cos\frac{20 \varphi }{3}}{2048 \sqrt{3}}-\frac{53}{640} \sin\frac{2 \varphi }{3}+\frac{1281 \sin\frac{8 \varphi }{3}}{5120}+\frac{25}{256} \sin\frac{14 \varphi }{3}+\frac{35 \sin\frac{20 \varphi }{3}}{2048}
\nonumber
\end{eqnarray}

\normalsize
\bigskip
\section{Authors' addresses}

\end{document}